\documentclass{article}
\usepackage{setspace}
\usepackage{titletoc}

\PassOptionsToPackage{numbers, compress}{natbib}

\usepackage[final]{neurips_2022}

\usepackage[utf8]{inputenc} 
\usepackage[T1]{fontenc}    
\usepackage{url}            
\usepackage{booktabs}       
\usepackage{amsfonts}       
\usepackage{nicefrac}       
\usepackage{microtype}      
\usepackage{xcolor,wrapfig}         

\newcommand{\alert}[1]{{#1}}
\newcommand{\alertt}[1]{{#1}}
\newcommand{\alerttt}[1]{{#1}}
\newcommand{\Alert}[1]{{ #1}}
\newcommand{\Alertt}[1]{{#1}}
\usepackage{xspace}
\usepackage{graphicx}
\usepackage{enumitem}
\usepackage[noend,algo2e,ruled,linesnumbered]{algorithm2e}

\usepackage[breaklinks=true,pdfstartview=FitH,colorlinks=false,pagebackref=true]{hyperref}
\usepackage{footnotebackref}
\usepackage{epsfig,epstopdf}
\usepackage{subcaption,multirow}
\usepackage{rotating}

\usepackage{amsmath}
\usepackage{amssymb}
\usepackage{mathtools}
\usepackage{amsthm}

\usepackage[capitalize,noabbrev,nameinlink]{cleveref}
\crefformat{equation}{\textup{#2(#1)#3}}
\crefrangeformat{equation}{\textup{#3(#1)#4--#5(#2)#6}}
\crefmultiformat{equation}{\textup{#2(#1)#3}}{ and \textup{#2(#1)#3}}
{, \textup{#2(#1)#3}}{, and \textup{#2(#1)#3}}
\crefrangemultiformat{equation}{\textup{#3(#1)#4--#5(#2)#6}}%
{ and \textup{#3(#1)#4--#5(#2)#6}}{, \textup{#3(#1)#4--#5(#2)#6}}{, and \textup{#3(#1)#4--#5(#2)#6}}

\Crefformat{equation}{#2Equation~\textup{(#1)}#3}
\Crefrangeformat{equation}{Equations~\textup{#3(#1)#4--#5(#2)#6}}
\Crefmultiformat{equation}{Equations~\textup{#2(#1)#3}}{ and \textup{#2(#1)#3}}
{, \textup{#2(#1)#3}}{, and \textup{#2(#1)#3}}
\Crefrangemultiformat{equation}{Equations~\textup{#3(#1)#4--#5(#2)#6}}%
{ and \textup{#3(#1)#4--#5(#2)#6}}{, \textup{#3(#1)#4--#5(#2)#6}}{, and \textup{#3(#1)#4--#5(#2)#6}}

\theoremstyle{plain}
\newtheorem{theorem}{Theorem}[section]

\theoremstyle{definition}

\theoremstyle{remark}

\usepackage[textsize=tiny]{todonotes}
\newcommand*{\lmss}{\fontfamily{lmss}\selectfont}

\newcommand{\T}{\mathsf{T}}

\newcommand{\I}{\mathcal{I}}

\newcommand{\J}{\mathcal{J}}

\newcommand{\N}{\mathbb{N}}

\newcommand\inner[2]{\left\langle #1, #2 \right\rangle}
\newcommand{\Tpg}{T_{\rm PG}^{\lambda}}

 \DeclareMathOperator*{\argmin}{arg\,min}

 \DeclareMathOperator*{\dist}{dist}
 
 \DeclarePairedDelimiter{\ceil}{\lceil}{\rceil}

 \newcommand{\norm}[1]{{\left\|{#1}\right\|}}

\newcommand{\bftab}{\fontseries{b}\selectfont}
\def\news{{\lmss news20}\xspace}
\def\webspam{{\lmss webspam}\xspace}
\def\rcv{{\lmss rcv1.binary}\xspace}
\def\gisette{{\lmss gisette\_scale}\xspace}
\def\duke{{\lmss duke}\xspace}
\def\leu{{\lmss leukemia}\xspace}
\def\colon{{\lmss colon-cancer}\xspace}
\def\El{{\lmss E2006-log1p}\xspace}
\def\Et{{\lmss E2006-tfidf}\xspace}

\title{Accelerated Projected Gradient Algorithms for Sparsity Constrained
Optimization Problems}
\author{%
  Jan Harold Alcantara\\
  Academia Sinica\\
  Taipei, Taiwan\\
  \texttt{jan.harold.alcantara@gmail.com}
  \And
  Ching-pei Lee\\
  Academia Sinica\\
  Taipei, Taiwan\\
  \texttt{leechingpei@gmail.com}
}
\begin{document}

\setlist{nolistsep}
\maketitle

\begin{abstract}
We consider the projected gradient algorithm for the nonconvex best subset selection
problem that minimizes a given empirical loss function under an $\ell_0$-norm
constraint.
Through decomposing the feasible set of the given sparsity constraint
as a finite union of linear subspaces, we present two acceleration
schemes with global convergence guarantees, one by same-space
extrapolation and the other by subspace identification.
The former fully utilizes the problem structure to greatly accelerate
the optimization speed with only negligible additional cost.
The latter leads to a two-stage meta-algorithm that first uses
classical projected gradient iterations to identify the correct
subspace containing an optimal solution, and then switches to
a highly-efficient smooth optimization method in the identified
subspace to attain superlinear convergence.
Experiments demonstrate that the proposed accelerated algorithms are
magnitudes faster than their non-accelerated counterparts as well as
the state of the art.
\end{abstract}

\section{Introduction}
We consider the sparsity-constrained optimization problem in $\Re^n$:
	\begin{equation}
	\label{eq:problem}
	\min\nolimits _{w \in A_s} f(w),
	\end{equation}
where $f$ is convex with $L$-Lipschitz continuous gradient, 
$s \in \N$,
and $A_s$ is the sparsity set given by
	\begin{equation}
		A_s \coloneqq \{ w \in \Re^n : \|w\|_0 \leq s\}, 
		\label{eq:As}
	\end{equation}
where $\|w\|_0$ denotes the $\ell_0$-norm that indicates the number of 
nonzero components in $w$. 
We further assume that $f$ is lower-bounded on $A_s$.

A classical problem that fits in the framework of 
\cref{eq:problem} is the best subset 
selection problem in linear regression \citep{BeaKM67a,HocL67a}. Given a 
response vector $y\in \Re^{m}$ 
and a design 
matrix of explanatory variables $X\in \Re^{m\times n}$, traditional linear regression minimizes a 
least squares (LS) loss 
function 
	\begin{equation}
		f(w) = \| y - Xw\|^2/2.
		\label{eq:ls_loss}
	\end{equation}
However, due to either high dimensionality in terms 
of the number of features $n$ or having significantly fewer instances $m$ than 
features $n$ (i.e., $m\ll n$), we often seek a linear model that selects only a 
subset of the 
explanatory variables that will best predict the outcome $y$. Towards
this goal, we can solve
\cref{eq:problem} with 
$f$ given by \cref{eq:ls_loss} to fit the training data while simultaneously
selecting the best-$s$
features. Indeed, such a sparse linear regression problem is 
fundamental in many 
scientific applications, such as high-dimensional statistical learning and 
signal processing \citep{JK17}. 
The loss in \cref{eq:ls_loss} can be generalized to the following
linear empirical risk to cover various tasks in machine learning beyond regression
\begin{equation}
	f(w) = g(Xw), \quad g(z) = \sum\nolimits_{i=1}^m g_i (z_i),
	\label{eq:erm}
\end{equation}
where $g$ is convex.
Such a problem structure makes evaluations of
the objective and its derivatives highly efficient, and such efficient
computation is a key motivation for our algorithms for
\cref{eq:problem}.

\paragraph{Related Works.} The discontinuous cardinality constraint in 
\cref{eq:problem} makes
the problem difficult to solve.
To make the optimization problem easier, a popular approach is to
slightly sacrifice the quality of the solution (either not strictly
satisfying the sparsity level constraint or the prediction performance
is deteriorated)
to use continuous surrogate functions for
the $\ell_0$-norm, which lead to a continuous nonlinear programming
problem, where abundant algorithms are at our disposal.
For instance, using a convex penalty surrogate such as the 
$\ell_1$-norm in the case of LASSO \citep{Tibshirani96}, the problem 
\cref{eq:problem} can be relaxed into a convex (unconstrained) one that can be
efficiently solved by many algorithms.
Other algorithms based on continuous nonconvex
relaxations such as the use of smoothly clipped absolute deviation
\citep{FanLi01} and the minimax concave penalty \citep{Zhang10}
regularizers are also popular in scenarios with a higher level of noise
and outliers in the data. 
However, for applications in which enforcing the constraints or
getting the best prediction performance is of utmost importance, solving
the original problem \cref{eq:problem} is inevitable. \alert{(For a detailed
review, we refer the 
	interested reader to \cite[Section 1]{BKM16}.)}
Unfortunately, methods for \cref{eq:problem} are not as well-studied as
those for the surrogate problems.
Moreover, existing methods are indeed
still preliminary and too slow to be useful in large-scale problems
often faced in modern machine learning tasks.

In view of the present unsatisfactory status for scenarios that
simultaneously involve high-volume data and need to get the best prediction performance,
this work proposes efficient algorithms to directly solve \cref{eq:problem}
with large-scale data.
\alertt{To our knowledge, all the most popular algorithms that directly tackles
\cref{eq:problem} without the use of surrogates involve using the
well-known projected gradient (PG) algorithm, at least as a major component}
\alert{\citep{BT11,BKM16,Blu12,BluDav09,BahRB13a}}.\footnote{\cite{GotTT18a} proposed 
an 
algorithm for a similar
optimization problem that minimizes $f(w) + C\norm{w}_0$ for some
$C > 0$. But whether it is equivalent to \cref{eq:problem} is unclear
because both problems are nonconvex, and for any prespecified sparsity
level $s$, it is hard to find $C$ that leads to a solution $w^*$
with $\norm{w^*}_0 = s$.}
	\cite{BT11} proved 
linear convergence of the objective value with the LS loss
function \cref{eq:ls_loss} for the iterates generated by PG under a
scalable restricted isometry property, which also served as their tool
to accelerate PG. However, given any problem instance, it is  hard, if not
computationally impossible, to verify whether the said property holds.
On the other hand, \cite{BKM16} established
global subsequential convergence to a stationary point for the iterates of PG 
on \cref{eq:problem} without the need for such isometry conditions, and their 
results 
are valid for general loss functions $f$ beyond \cref{eq:ls_loss}.
While some theoretical guarantees are 
known, the practicality of PG for solving 
\cref{eq:problem} remains a big problem in real-world applications
as its empirical convergence speed tends to be slow.
\alertt{The PG approach is called iterative hard thresholding (IHT) in studies of
compressed sensing \citep{BluDav09} that mainly focuses on the
LS case.
To accelerate IHT, several approaches that alternates between a PG
step and a subspace optimization step are also proposed
\cite{Blu12,BahRB13a}, but such methods mainly focus on the
LS case and statistical properties, while their convergence speed is
less studied from an optimization perspective.}
Recently, ``acceleration'' approaches for PG on general nonconvex regularized problems have
been studied in \cite{NIPS2015_f7664060,WCP18}. While their proposed
algorithms are also applicable to \cref{eq:problem}, the obtained convergence 
speed for nonconvex problems is
not faster than that of PG.

\Alertt{This work is inspired by our earlier work \cite{JHA22a}, which
considered a much broader class of problems without requiring convexity nor 
differentiability assumptions
for $f$, and hence obtained only much weaker convergence results,
with barely any convergence rates, for such general problems.}

\paragraph{Contributions.} In this work, we revisit the PG
algorithm for solving the general problem
\cref{eq:problem} and propose two acceleration schemes by leveraging the 
combinatorial nature of $\ell_0$-norm.
In particular, we decompose the feasible set $A_s$ as the finite union
of $s$-dimensional linear subspaces, each representing a subset of the
coordinates $\{1,\dotsc,n\}$, as
detailed in \cref{eq:As_decomp} of \cref{sec:pgm}. Such subspaces are utilized 
in devising techniques to efficiently accelerate
PG.
Our first acceleration scheme is based on a same-space extrapolation
technique such that we
conduct extrapolation only when two consecutive iterates $w_{k-1}$ and
$w_k$ lie in the same 
subspace, and the step size for this
extrapolation is determined by a spectral initialization combined with
backtracking 
to ensure sufficient function decrease.
This 
is motivated by the observation that for \cref{eq:erm},
objective and derivatives at the extrapolated point can be inferred
efficiently through a linear combination of $X w_{k-1}$ and $X w_{k}$.
The second acceleration technique starts with plain PG, and
when consecutive iterates stay in the same subspace, it
begins to alternate between a full PG step and a truncated Newton step
in the subspace to obtain superlinear convergence with extremely low
computational cost.
Our main contributions are as follows:
\begin{enumerate}[leftmargin=*]
	\item We prove that PG for \cref{eq:problem} is globally 
	convergent to a local optimum with a local linear rate, improving upon the sublinear results of \citet{BKM16}. We emphasize 
	that our framework, like 
	\citep{BKM16}, is applicable to 
	general loss functions $f$ satisfying the convexity and smoothness requirements, and therefore 
	covers not only the classical sparse regression problem but also
	many other ones encompassed by the empirical risk minimization
	(ERM) framework.
	
	\item By decomposing $A_s$ as the union of linear subspaces,
		we further show that PG is provably capable of identifying a subspace containing a 
	local optimum of \cref{eq:problem}.
	By exploiting this property,
	we propose 
	two acceleration strategies with practical implementation and convergence
	guarantees
	for the general problem class \cref{eq:problem}. 
	Our acceleration provides both computational and theoretical
	advantages for convergence, and can in particular obtain
	superlinear convergence.
	
	\item In comparison with existing acceleration methods for
		nonconvex problems \citep{NIPS2015_f7664060,WCP18}, this work
		provides new acceleration schemes with faster theoretical speeds
		(see \cref{thm:rate2,thm:newton}), and beyond being applied to
		the classical PG algorithm,
		those schemes can
	also easily be combined with existing accelerated PG approaches to
	further make them converge even faster.
	\item Numerical experiments exemplify the significant improvement
		in both iterations and running time brought by our
	acceleration methods, in particular over the projected gradient
	algorithm by \cite{BKM16} as well as the accelerated proximal
	gradient method for nonconvex problems proposed by
	\cite{NIPS2015_f7664060}.
\end{enumerate}

This work is organized as follows. We review the projected gradient algorithm 
and prove its local linear convergence and subspace identification
for arbitrary smooth loss functions in \cref{sec:pgm}. In \cref{sec:accelerate}, 
we propose the acceleration schemes devised through
decomposing the constraint set in \cref{eq:problem} into
subspaces of $\Re^n$.
Experiments in \cref{sec:exp} then illustrate the effectiveness of
the proposed acceleration techniques, and \cref{sec:conclusion} concludes this
work.
All proofs, details of the experiment settings, and additional
experiments are in the appendices.

\section{Projected Gradient Algorithm}\label{sec:pgm}
The \emph{projected gradient algorithm} for solving \cref{eq:problem} is given 
by the iterations
	\begin{equation}
	w ^{k+1} \in  \Tpg (w^k) \coloneqq P_{A_s}(w ^k - \lambda \nabla f(w 
	^k)), \label{eq:pgm}
	\end{equation}
where $P_{A_s}(w )$ denotes the projection of $w $ onto $A_s$,
which is set-valued because of the nonconvexity of $A_s$.
When $f$ is given by \cref{eq:ls_loss},
global linear convergence of this algorithm under a restricted
isometry condition is established in \citep{BT11}.
For a general convex $f$ with $L$-Lipschitz continuous
gradients, that is, 
	\begin{equation}
	\label{eq:L-lipschitz}
	\|\nabla f(w) - \nabla f(w') \| \leq L \|w - w '\| \quad \forall 
	w, w' \in \Re^n,
	\end{equation}
the global subsequential convergence of \cref{eq:pgm} is proved in 
\cite{BKM16}, but neither global nor local rates of convergence is provided.
In this section, we present an alternative proof of global convergence and
more importantly establish its local linear convergence.

A useful observation that we will utilize in the proofs of
our coming convergence results is that the nonconvex set
$A_s$ given by \cref{eq:As} can be decomposed as 
a finite union of subspaces in $\Re^n$:
	\begin{equation}
	\label{eq:As_decomp}
		A_s = \bigcup\nolimits _{J\in \J _s} A_{J}, \quad A_{J} \coloneqq {\rm
		span} \{ e_j: j\in J\}, \quad
	 \J_s \coloneqq \left\{ J\subseteq \{1,2,\dots,n\}: |J| = s\right\},
	\end{equation}
where $e_j$ is the $j$th standard unit vector in $\Re^n$. \alert{Throughout 
this paper, we assume that $\lambda\in (0,L^{-1})$.}

\begin{theorem}
	\label{thm:pgm}
	\alert{Let $\{w ^k\}$ be a sequence 
	generated by \cref{eq:pgm}.}  Then:
	\begin{enumerate}[leftmargin=*]
		\item[(a)] (Subsequential convergence) 
			\alert{Either $\{f(w^k)\}$ is strictly decreasing, or there exists 
			$N>0$ 
			such that $w^k = w^N$ for all $k\geq N$. In addition,} any 
			accumulation point $w^*$ of $\{w^k\}$ satisfies
			$w^* \in P_{A_s} (w^* - \lambda \nabla f(w^*))$, and is
			hence a stationary point of \cref{eq:problem}.
		\item[(b)] (Subspace identification and full convergence) There 
		exists $N \in \N$ such that
		\begin{equation}
		\{w^k\}_{k=N}^{\infty}\subseteq \bigcup\nolimits
		_{J\in\I_{w^*}} A_J,\quad
		\quad \I_{w^*} \coloneqq \{ J\in \J_s: w^* \in
		A_J\}.
		\label{eq:identification}
		\end{equation} 
		whenever $w^k\to w^*$. In particular, if $\Tpg (w ^*)$ is a singleton 
		for an accumulation point $w^*$ of $\{w^k\}$, then
		$w^*$ is a local minimum for \cref{eq:problem}, $w ^k\to w ^*$, and 
			\cref{eq:identification} holds. 
	\item[(c)] ($Q$-linear convergence) If $\Tpg (w ^*)$ is a singleton 
	for an accumulation point $w^*$ and $w \mapsto w -
	\lambda \nabla f(w)$ is a contraction over $A_J$ for all $J\in
	\I_{w^*}$, then $\{w^k\}$ converges to $w^*$ at a $Q$-linear rate.
	In other words, there is $N_2 \in \N$ and $\gamma \in [0,1)$ such that
		\begin{equation}
		\label{eq:Q-linear}
		\norm{w^{k+1} - w^*} \leq \gamma \norm{w^{k} - w^*}, \quad \forall
		k \geq N_2.
		\end{equation}
	\end{enumerate}
\end{theorem}

	It is well-known that an optimal solution of \cref{eq:problem} is
	also a
	stationary point of it \cite[Theorem 2.2]{BE13}, and therefore (a) proves the 
	global subsequential convergence of PG to candidate solutions of 
	\cref{eq:problem}. Consider $z^* \coloneqq w^* - \lambda \nabla f(w^*)$,
	and let $\tau$ be a permutation of $\{1,\dotsc,n\}$ such that
	$z^*_{\tau(1)} \geq z^*_{\tau(2)} \geq \cdots \geq z^*_{\tau(n)}$.
	The requirement of $\Tpg (w ^*)$ being a singleton in
	\cref{thm:pgm} (b) then simply means the mild condition of
	$z^*_{\tau(s)} > z^*_{\tau(s+1)}$, which is almost always true in
	practice.
	The requirement for (c) can be fulfilled when $f$ confined to $A_J$ is
	strongly convex, even if $f$ itself is not.
	This often holds true in practice when $f$ is of the form
	\cref{eq:erm} and we restrict $s$ in \cref{eq:problem} to be
	smaller than the number of data instances $m$, and is thus also mild.
	The existence of a stationary point can be guaranteed when $\{w^k\}$
	is a bounded sequence, often guaranteed when $f$ is coercive on
	$A_J$ for each $J \in \J_s$. 

In comparison to existing results in \cite{BKM16,HA13a,BolSTV18a},
parts (b) and (c) of \cref{thm:pgm} are new.
In particular, part (b) provides a full convergence result that usually
requires stronger regularity assumptions like the
Kurdyka-{\L}ojasiewicz (KL) condition \cite{HA13a,BolSTV18a} \alert{(see also 
\cref{eq:KL2})} \alertt{that requries the objective value to decrease
	proportionally with the minimum-norm subgradient in a neighborhood
	of the accumulation point, but we just need the very mild
	singleton condition right at the
accumulation point only}. Part (c) gives a local 
linear convergence for the PG iterates even if the
problem is nonconvex, while the rates in \cite{BolSTV18a} requires a
KL condition and the rate is measured in the objective value.

The following result further provides rates of convergence of the 
objective values even without the conventional KL assumption.
The first rate below follows from \cite{LeeW19a}.

\begin{theorem}
	\label{thm:sublinear}
	\alert{Let $\{w ^k\}$ be a sequence generated by \cref{eq:pgm}}.
	If $w^k\to w^*$, such as when $\Tpg (w^*)$ is a singleton at an 
	accumulation 
	point $w^*$ of \cref{eq:pgm}, then
		\begin{equation}
		\label{eq:sublinear}
		f(w^k) - f(w^*) = \alert{o}(k^{-1}).
		\end{equation}
		Moreover, under the hypothesis of \cref{thm:pgm} (c), the
		objective converges to $f(w^*)$ R-linearly, i.e., 		
			\begin{equation}
			f(w^k) - f(w^*) = O(\exp(-k)).
			\label{eq:linear}
			\end{equation}
\end{theorem}

By using \cref{thm:pgm}, we can also easily get
rates faster than \cref{eq:sublinear} under a version of the KL
condition that is easier to understand and verify than those assumed
in existing works.
In particular, existing analyses require the KL condition to hold in a
neighborhood in $\Re^n$ of an
accumulation point, but we just need it
to hold around $w^*$ within $A_J$ for the restriction $f|_{A_J}$ for each $J \in 
\I_{w^*}$.
These results are postponed to \cref{thm:rate2} in the next section as
the PG method is a special case of our acceleration framework.

\section{Accelerated methods}
\label{sec:accelerate}

The main focus of this work is the proposal in this section of new
techniques with solid convergence guarantees to accelerate the PG
algorithm presented in the preceding section.
Our techniques fully exploit the subspace identification property
described by the inclusion \cref{eq:identification}, as well as the
problem structure of \cref{eq:erm} to devise efficient algorithms.

We emphasize that the two acceleration strategies described below can
be combined together, and they are also widely applicable such that
they can be employed to other existing algorithms for
\cref{eq:problem} as long as such algorithms have a property similar to
\cref{eq:identification}.

\subsection{Acceleration by extrapolation}
Traditional extrapolation techniques are found in the realm of convex
optimization to accelerate algorithms 
\citep{BeckTeboulle09,Nesterov83} with guaranteed convergence
improvements, but were often only adopted as heuristics in the
nonconvex setting, until some recent works showed that theoretical
convergence can also be achieved \citep{NIPS2015_f7664060,WCP18}.
However, unlike the convex case, these
extrapolation strategies for nonconvex problems do not lead to faster
convergence speed nor an intuitive reason for doing so.
An extrapolation step proceeds by choosing a 
\emph{positive} stepsize along the 
direction determined by two consecutive iterates. That is, given two iterates 
$w^{k-1}$ and $w^k$, an 
intermediate point $z^k \coloneqq w^k + t_k (w^k-w^{k-1})$ for some stepsize 
$t_k \geq 0$ is first calculated 
before applying the original algorithmic map ($\Tpg$ in our
case).\footnote{it is clear that if $t_k \equiv 0$, we reduce back to
the original algorithm.}

Another popular acceleration scheme for gradient algorithms is the
spectral approach pioneered by \cite{JB88a}.
They take the differences of the gradients and of the iterates in two
consecutive iterations to estimate the curvature at the
current point, and use it to decide the step size for updating along
the reversed gradient direction.
It has been shown in \cite{WriNF09a} that equipping this step size
with a backtracking procedure leads to significantly faster
convergence for proximal gradient on regularized optimization
problems, which includes our PG for \cref{eq:problem} as a special case.

To describe our proposed double acceleration procedure that combines
extrapolation and spectral techniques, we first observe that all PG iterates 
lie on $A_s$, 
and that $A_s$ can be finitely decomposed as \cref{eq:As_decomp}. When two 
consecutive iterates lie 
on the same convex subspace $A_J$ for some $J\in \J_s$, within these two iterations, we
are actually conducting convex optimization.
In this case, an extrapolation step within $A_J$ is reasonable because it will 
not
violate the constraint, and acceleration can be expected from the
improved rates of accelerated proximal gradient on convex problems in
\cite{BeckTeboulle09,Nes13a}.
Judging from \cref{thm:pgm} (b),
the corresponding $J$ is also a candidate index set that belongs to
$\I_{w^*}$, so extrapolation within $A_J$ makes further sense.
We set $t_k = 0$ to skip the extrapolation step if $d^k$ is not a
descent direction for $f$ at $w^k$.
Otherwise, we start from some $\hat t_k > 0$ decided by the curvature
information of $f$, and then execute a backtracking linesearch along
$d^k \coloneqq w^k - w^{k-1}$ to set
$t_k = \eta^i \hat t_k$ for the smallest integer $i\geq 0$ that provides
sufficient descent
\begin{equation}
	f (w^k + t_k d^k) \leq f (w^k) - \sigma t_k^2 \|d^k\|^2,
	\label{eq:suffdecrease}
\end{equation}
given parameters $\eta,\sigma \in (0,1)$.
We then apply \cref{eq:pgm} to $z^k = w^k +
t_k d^k$ to obtain $w^{k+1}$.

For the spectral initialization $\hat t_k$ for accelerating the
convergence,
instead of directly using approaches of \cite{JB88a,WriNF09a} that takes
the reversed gradient as the update direction,
we need to devise a different mechanism as
our direction $d^k$ is not directly related to the gradient.
We observe that for the stepsize
\begin{equation}
	\label{eq:BB}
	\alpha_k \coloneqq \inner{s^k}{s^k}/\inner{s^k}{r^k}, \quad
	s^k \coloneqq w^{k} - w^{k-1}, \quad r^k \coloneqq \nabla f(w^k) -
	\nabla f(w^{k-1})
\end{equation}
used in \cite{JB88a}, the final update $-\alpha_k\nabla f(w^k)$ is actually
the minimizer of the following subproblem
\begin{equation}
	\label{eq:quad}
	\min_{d\in\Re^n}\quad \inner{\nabla f(w^k)}{d} +  \norm{d}^2/ (2
	\alpha_k).
\end{equation}
By juxtaposing the above quadratic problem and the upper bound provided by 
the descent lemma \cite[Lemma 5.7]{Beck17},
we can view $\alpha_k^{-1}$ as an estimate of the local Lipschitz
parameter that could be much smaller than $L$ but still guarantee
descent of the objective.
We thus follow this idea to decide $\hat t_k$ using such curvature
estimate and the descent lemma by
\begin{equation}
	\label{eq:BB_extrapolation}
	\hat t_k = \argmin_{t\geq 0}\; \inner{\nabla f(w^k)}{t d^k} +
	\norm{t d^k}^2 / (2\alpha_k)\quad \Leftrightarrow \quad
	\hat t_k = -\inner{\alpha_k \nabla f(w^k)}{
	d^k}/ \norm{d^k}^2.
\end{equation}
Another interpretation of \cref{eq:BB} is that $\alpha_k^{-1} I$
also serves as an estimate of $\nabla^2 f(w^k)$,\footnote{As $\nabla f$ is Lipschitz
continuous, it is differentiable almost everywhere.
Here, we denote by $\nabla ^2 f(w^k)$ a generalized Hessian of $f$
at $w$, which is well-defined for $f$ with Lipschitz continuous
gradient \citep{HUJ84a}.}
and the objective in \cref{eq:quad} is a low-cost approximation of the
second-order Taylor expansion of $f$.
However, we notice that for problems in the form of \cref{eq:erm} and
with $d^k \in A_J$, the exact second-order Taylor expansion
\begin{equation}
	\label{eq:taylor2}
f(w^k + t d^k) \approx f(w^k) + t\inner{\nabla f(w^k)}{d^k} + 
t^2\inner{ \nabla ^2 f(w^k) d^k} { d^k }/2
\end{equation}
can be calculated efficiently.
In particular, for \cref{eq:erm}
and any $d^k \in A_J$, we get from $X d^k = X_{:,J} d^k_J$:
\begin{equation}
\label{eq:gradhess}
\begin{aligned}
	\nabla f(w^k)^{\top} d^k &= \nabla g\left( (X w^k) \right)^{\top}
	\left( X_{:,J} d^k_J \right),\\
	\inner{\nabla^2 f(w^k) d^k}{d^k}  &= \inner{(X_{:,J} d^k_J)}{\nabla^2
		g\left( (Xw^k) \right) (X_{:,J} d^k_J)},
	\end{aligned}
\end{equation}
which can be calculated in $O(ms)$ time by computing $X_{:,J} d^k_J$ first.
This $O(ms)$ cost is much cheaper than
the $O(mn)$ one for evaluating the full gradient of $f$ needed in the
PG step,
so our extrapolation plus spectral techniques has only negligible
cost.
Moreover, for our case of $d^k = w^k - w^{k-1}$,
we can further reduce the cost of calculate $X_{:,J} d^k_J$ and thus
\cref{eq:gradhess} to $O(m)$ by recycling intermediate computational
results needed in evaluating $f(w^k)$ through $X_{:,J} d^k_J = X w^k - X
w^{k-1}$.
With such tricks for efficient computation,
we therefore consider the more accurate approximation
to let $\hat t_k$ be the scalar that minimizes the quadratic function
on the right-hand side of \cref{eq:taylor2} for problems in the form
\cref{eq:erm}. That is, we use
	\begin{equation}
		\label{eq:tk}
		\hat t_k  \coloneqq - \inner{\nabla f(w^k)}{d^k}/ \inner{ 
		\nabla ^2 f(w^k) d^k} { d^k }.
	\end{equation}
Finally, for both \cref{eq:tk} and \cref{eq:BB_extrapolation},
we safeguard $\hat t_k$ by
\begin{equation}
	\hat t_k \leftarrow P_{[c_k \alpha_{\min}
	, c_k \alpha_{\max} ] } \left( \hat t_k \right)
	\label{eq:hattk}
\end{equation}
for some fixed $\alpha_{\max} \geq \alpha_{\min} > 0$,
where
\begin{equation}
 c_k \coloneqq \norm{(\nabla f(w^k))_{J} }/ (\zeta_k \norm{d^k}),
 \quad
	\zeta_k \coloneqq - \inner{d^k}{\nabla f(w^k)} /
	(\norm{d^k} \norm{(\nabla f(w^k))_{J} }) \in (0,1].
	\label{eq:cosine}
\end{equation}
We also note that the low cost of evaluating $X d^k$ is also the key
to making the backtracking in \cref{eq:suffdecrease} practical,
as each $f(w^k + \eta^i \hat t_k d^k)$ can be calculated in $O(m)$
time through linear combinations of $X w^k$ and $X d^k$.
The above procedure is summarized in
\cref{alg:extrapolation} \alert{with global convergence guaranteed 
by \cref{thm:global_acce}. In \cref{thm:rate2}, we establish its 
full convergence as well as its
convergence rates under a KL condition at $w^*$: there exists 
neighborhood $U \subset \Re^n$ of
$w^*$, $\theta \in [0,1]$, and $\kappa > 0$ such that for every $J
\in \I_{w^*}$, 
\begin{equation}
\left(f(w) - f(w^*)\right)^{{\theta}} \leq 
{\kappa } \norm{(\nabla
	f(w))_{J}},\quad \forall w\in A_J \cap U.
\label{eq:KL2}
\end{equation}
We denote
by $n_k$ the number of successful 
extrapolation 
steps in the
first $k$ iterations of \cref{alg:extrapolation}.}

\begin{theorem}
	\label{thm:global_acce}
	Under the hypotheses of \cref{thm:pgm}, any
	accumulation point of a sequence generated by \cref{alg:extrapolation}
	is a stationary point.
\end{theorem}
\begin{theorem}
\label{thm:rate2}
\alert{Consider either \cref{eq:pgm} or
\cref{alg:extrapolation} with }$\eta,
\sigma,{\epsilon}\in (0,1)$, and $\alpha_{\max} \geq \alpha_{\min} > 0$,
and suppose that there is an accumulation
point $w^*$ of the iterates at which
the KL condition holds. Then $w^k\to w^*$. Moreover, the following rates
hold:
\begin{enumerate}[leftmargin=*]
	\item[(a)] If \alert{$\theta \in (1/2,1)$}: $f(w^k) - f(w^*) =
	O((k + n_k)^{-1/\alert{(2\theta -1)}})$.
\item[(b)] If \alert{$\theta \in (0, 1/2]$}: $f(w^k) - f(w^*) = O(\exp(-(k
		+ n_k)))$.
	\item[\alert{(c)}] \alert{If \alert{$\theta = 0$}:
there 
is $k_0 \geq 0$ such that $f(w^k) = f(w^*)$ for all $k
	\geq k_0$}.
\end{enumerate}
\end{theorem}
\alert{We stress that convexity of $f$ is not required in 
	\cref{thm:global_acce,thm:rate2} \Alertt{except the second half of the
		last item of \cref{thm:rate2}}.} 
There are several advantages of the proposed extrapolation strategy
over existing ones in \cite{NIPS2015_f7664060,WCP18}.
The most obvious one is the faster rates in \cref{thm:rate2} over
PG such that each successful extrapolation step in our method
contributes to the convergence speed, while existing methods only
provide the same convergence speed as PG.
Next, existing strategies only use prespecified step sizes without information
from the given problem nor the current progress,
and they only restrict such step sizes to be within $[0,1]$.
Our method, on the other hand, fully takes advantage of the function
curvature and can allow for arbitrarily large step sizes to better
decrease the objective.
In fact, we often observe $t_k \gg 1$ in our numerical
experiments.
Moreover, our acceleration techniques utilize the nature of
\cref{eq:As_decomp} and \cref{eq:erm} to obtain very efficient
implementation for ERM problems such that the per-iteration cost of
\cref{alg:extrapolation} is almost the same as that of PG, while the approach of
\cite{NIPS2015_f7664060} requires evaluating $f$ and $\nabla f$ at two
points per iteration, and thus has twice the per-iteration cost.

A finite termination result similar to \cref{thm:rate2} (c) is presented in 
\cite{LiuY17a}
under a H\"olderian error bound that is closely related to the KL
condition, but their result requires convexity of both the smooth term
and the regularizer, so it is not applicable to \cref{eq:problem} that
involves a nonconvex constraint.


\begin{algorithm2e}[tb]
\DontPrintSemicolon
\caption{Accelerated projected gradient algorithm by extrapolation (APG)}
\label{alg:extrapolation}
		Given an initial vector
		$w^0\in\Re^n$ and parameters ${\epsilon},\eta,\sigma\in (0,1)$,
		$\alpha_{\max} \geq \alpha_{\min} > 0$, $\lambda \in (0,1/L)$.

		\For{$k=0,1,2,\dotsc$}{
			\If{$k>0$; $w^{k-1}$ and $w^k$ activate the same $A_J$;
			and ${\zeta_k \geq \epsilon}$}{
			\label{line:if}
				$d^k \leftarrow w^k - w^{k-1}$, and compute $\hat t_k$
				from \cref{eq:hattk} with either \cref{eq:BB_extrapolation} or \cref{eq:tk}

			\For{$i=0,1,\dotsc$}{
				$t_k \leftarrow \eta^i \hat t_k$

				\lIf{\cref{eq:suffdecrease} is satisfied}{$z^k
				\leftarrow w^k + t_k d^k$, and break}
			}
		
		}
		\lElse{
			$z^k \leftarrow w^k$
		}
		$w^{k+1} \leftarrow\Tpg (z^k)$
	}
\end{algorithm2e}

\subsection{Subspace Identification}
\label{sec:identify}

In line with the above discussion, we interpret 
\cref{eq:identification} 
as a theoretical property guaranteeing that the iterates of the projected gradient algorithm 
\cref{eq:pgm}
will eventually identify the subspaces $A_J$ that contain a candidate solution 
$w^*$ after a finite number of iterations. Consequently, the task of minimizing
$f$ over the nonconvex set $A_s$ can be reduced to a convex optimization 
problem of minimizing $f$ over $A_J$. Motivated by this, we present a 
two-stage algorithm described in \cref{alg:identification} that switches to 
a high-order method for smooth convex optimization after a candidate
piece $A_J$ is identified to obtain even faster convergence.
Since $\nabla f$ is assumed to be 
Lipschitz continuous, the generalized Hessian of $f$ exists everywhere
\citep{HUJ84a}, so we may employ a semismooth Newton
(SSN) method \citep{QiS93a} with backtracking linesearch to get a faster
convergence speed with low cost (details in \cref{subsec:ssn}). In particular, we 
reduce the computation costs by considering the restriction of $f$ on the
subspace $A_J$ by 
treating
the coordinates not in $J$ as non-variables so that the problem
considered is indeed smooth and convex.
As we cannot know a priori whether $I_{w^*}$ is indeed identified,
we adopt the approach implemented in 
\cite{LeeW12a,YSL20a-nourl,LeeCP20-nourl} to consider it
identified when $w^k$ activates the same $A_J$ for long enough
consecutive iterations.
To further safeguard that we are not optimizing over a wrong subspace,
we also incorporate the idea of
\cite{Wri12a,BarIM20a,YSL20a-nourl,LeeCP20-nourl}
to periodically alternate to a PG step \cref{eq:pgm} after switching
to the SSN stage.
A detailed description of this two-stage algorithm is in
\cref{alg:identification}.

In the following theorem, we show that superlinear convergence can be
obtained for \cref{alg:identification} even if we take only one SSN
step every time between two steps of \cref{eq:pgm}, using a simplified
setting of twice-differentiability.
\alert{For our next theorem}, we need to introduce some additional notations.
Given any $w \in A_J$, we use $f_{J}(w_J) \coloneqq f(w)$ to
denote the function of considering only the coordinates of $w$ in
$J$ as variables and treating the remaining as constant zeros.
\alert{We assume that the conditions of \cref{thm:pgm} (b) hold with $w^* \in 
A_s$, and that $f$ is
twice-differentiable around a neighborhood $U$ of $w^*$ with
$\nabla^2 f_J$ Lipschitz continuous in $U$ and $\nabla^2 f_J(w^*)$
positive definite for all $J \in \I_{w^*}$.}
\begin{theorem}
\alert{Suppose that starting after $k
\geq N$ and $P_{A_s}(w^k) \subset U$}, we conduct $t$ Newton steps
between every two steps of \cref{eq:pgm} for $t \geq 1$:
	\begin{equation}
	\label{eq:newton}
	w^{k,0} \in  P_{A_s}(w^k), \,
	\begin{cases}
	J &\in  \I_{w^{k,0}},\\
	w^{k,j+1}_{i} &= 0, \quad \forall i \notin J, \quad j=1,\dotsc,t-1,\\
	w^{k,j+1}_J &= w^{k,j}_J -
	\nabla^2 f_J(w^{k,j}_J)^{-1} \nabla f_J(w^{k,j}_J),\\
	\end{cases}
	\,
	w^{k+1} \in \Tpg (w^{k,t}).
	\end{equation}
Then $w^k \rightarrow w^*$ at a $Q$-quadratic rate.
\label{thm:newton}
\end{theorem}

In practice, the linear system for obtaining the SSN step is only
solved inexactly via a (preconditioned) conjugate gradient (PCG) method,
and with suitable stopping conditions for PCG and proper algorithmic
modifications such as those in \cite{MCY19a,MorYZZ20a},
superlinear convergence can still be obtained easily.
Interested readers are referred to \cref{subsec:ssn} for a more detailed
description of our implementation.

\begin{algorithm2e}[tb]
\DontPrintSemicolon
\caption{Accelerated projected gradient algorithm by subspace identification (PG+)}
\label{alg:identification}

		Given an initial vector $w^0\in\Re^n$ and 
		$S,t\in \mathbb{N}$. Set Unchanged $\leftarrow 0$.

		\For{$k=0,1,2,\dotsc$}{
		\If{$k > 0$, and $w^{k-1}$ and $w^k$ activate the same component of $A_s$}{
			Let $J\in \J_s$ correspond to the activated component

			Unchanged $\leftarrow$ Unchanged $+1$
		}
		\lElse{
			Unchanged $\leftarrow 0$
		}
		
		\If{Unchanged $\geq S$}
		{
			$y^k \leftarrow P_{A_J}(w^k)$ and
			use $t$ steps of SSN described in \cref{subsec:ssn},
			starting from $y^k$, to find $z^k$ that approximately
			minimizes $f|_{A_J}$

			\lIf{SSN fails}{
				$z^k \leftarrow w^k$ and Unchanged $\leftarrow 0$.
			}
		}
		\lElse{
			$z^k \leftarrow w^k$
		}
		$w^{k+1} \leftarrow\Tpg (z^k)$
	}
\end{algorithm2e}

\section{Experiments}
\label{sec:exp}
In this section, we conduct numerical experiments to demonstrate the accelerated techniques presented 
in \cref{sec:accelerate}. We employ \cref{alg:extrapolation} (APG) with \cref{eq:tk} to 
accelerate PG, and 
further accelerate APG by incorporating subspace identification described in 
\cref{alg:identification}, which we denote by APG+.\footnote{That is, if $Unchanged<S$ in 
\cref{alg:identification}, we calculate $z^k$ as in \cref{alg:extrapolation}} 
Comparisons with the extrapolated PG algorithm of \citet{NIPS2015_f7664060}, 
which we denote by PG-LL, are also presented. 
PG-LL is a state-of-the-art approach for nonconvex regularized
optimization and thus suitable for \cref{eq:problem}.
For $f$ in \cref{eq:problem}, we consider both LS
\cref{eq:ls_loss} and logistic regression (LR)
\begin{equation}
	f(w) = \sum\nolimits_{i=1}^m \log\left( 1 + \exp\left( - y_i x_i^{\top}
	w \right) \right) + \mu\norm{w}^2/2,
	\label{eq:lr}
\end{equation}
where $(x_i,y_i) \in \Re^n \times \{-1,1\}$, $i=1,\dotsc, m$, are
the training instances, and $\mu>0$ is a small regularization parameter added to make the logistic loss 
coercive.

\begin{figure}[tbh] 
	\centering 
	Logistic regression
	\begin{tabular}{ccc}
		\begin{subfigure}[b]{0.29\textwidth} 
			\includegraphics[width=.86\linewidth]{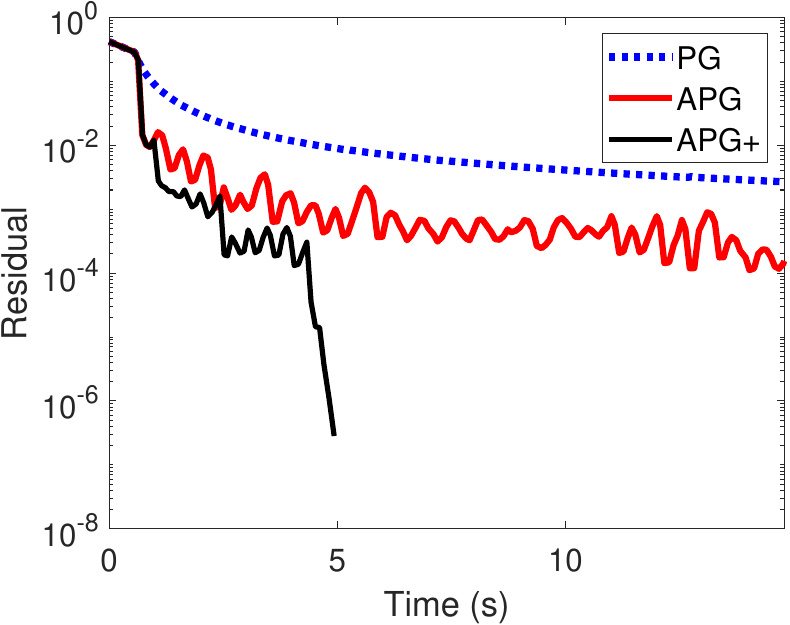}
			\caption{\news, $s=\ceil {0.01m}$} 
		\end{subfigure}& 
		\begin{subfigure}[b]{0.29\textwidth} 
		\includegraphics[width=.86\linewidth]{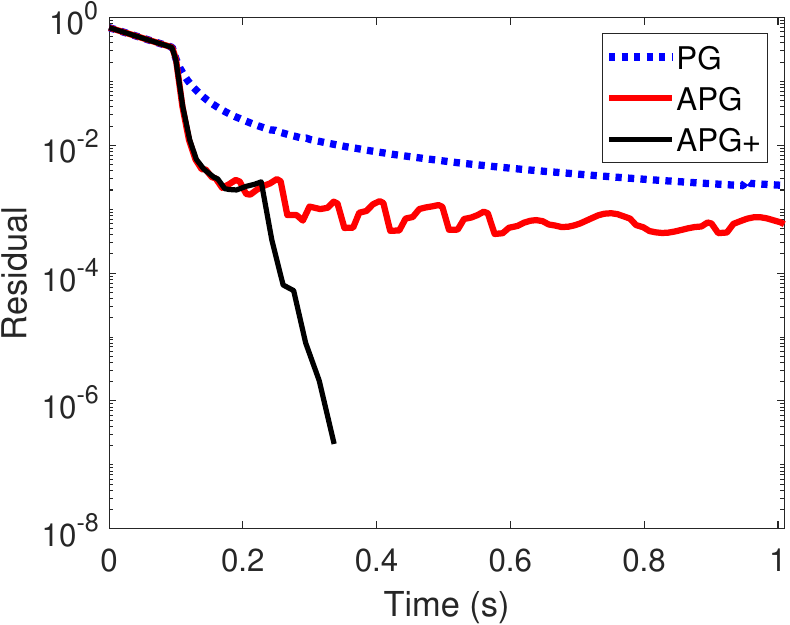}
		 
		\caption{\rcv, $s=\ceil {0.01m}$} 
		\end{subfigure}& 
		\begin{subfigure}[b]{0.29\textwidth} 
		\includegraphics[width=.86\linewidth]{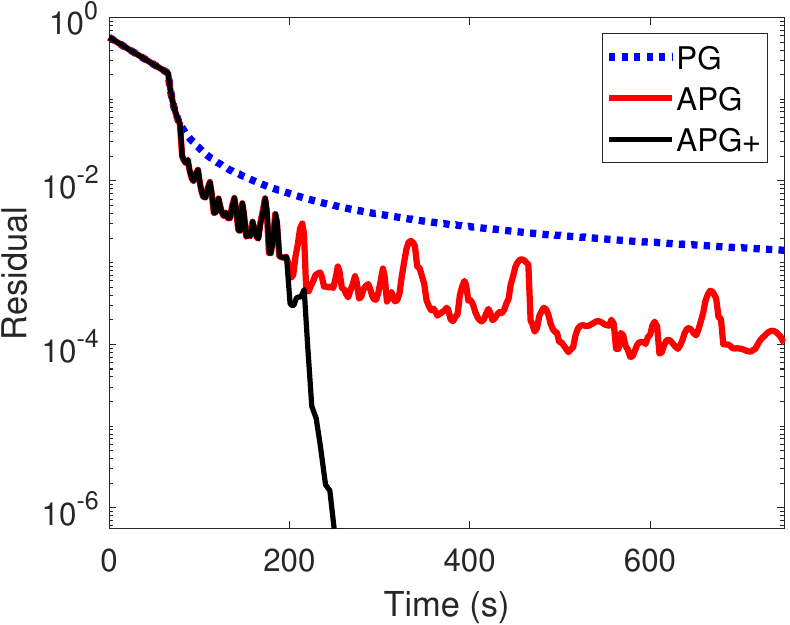}
		 
		\caption{\webspam, $s=\ceil {0.001m}$} 
		\end{subfigure}
	\end{tabular} 
	\smallskip

	Least square

	\begin{tabular}{cc}
		\begin{subfigure}[b]{0.29\textwidth} 
			\includegraphics[width=.86\linewidth]{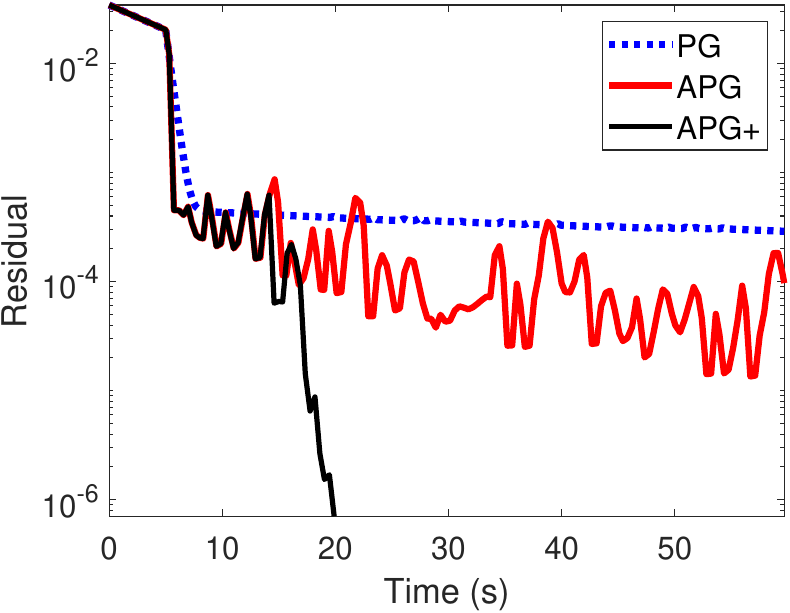}
			\caption{\El, $s=\ceil {0.01m}$} 
		\end{subfigure}& 
		\begin{subfigure}[b]{0.29\textwidth} 
			\includegraphics[width=.86\linewidth]{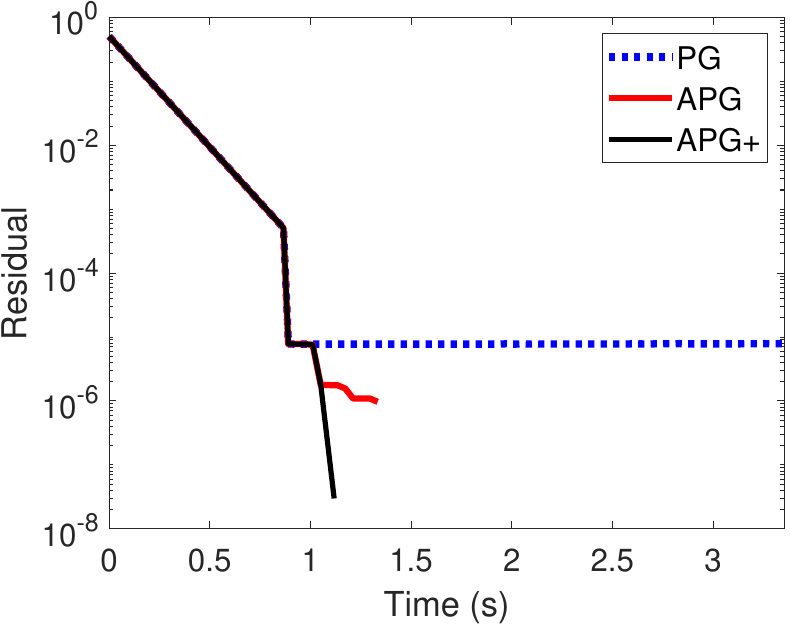}
			\caption{\Et, $s=\ceil {0.01m}$} 
		\end{subfigure}
	\end{tabular} 
	\caption{Experiment on sparse regularized LR and
		LS. We present time v.s. residual in
	\cref{eq:opt}.}
	\label{fig:main} 
\end{figure}

The algorithms are implemented in MATLAB and tested 
with public datasets in \Alert{\cref{tbl:data,tbl:data2} in 
\cref{app:expdetails}. All algorithms compared start from $w^0=0$ and 
terminate when
the first-order optimality condition 
\begin{equation}
\text{Residual}(w) \coloneqq \norm{w - P_{A_s}\left( w - \lambda \nabla f\left( 
w
	\right) \right)}/ (1+\norm{w} + \lambda \norm{\nabla f \left( w
	\right)}) < \hat \epsilon
\label{eq:opt}
\end{equation}
is met for some given $\hat \epsilon > 0$. 	More setting and parameter 
details of our experiments are in \cref{app:expdetails}.

\paragraph{Comparisons of algorithms for large datasets.} To fit the practical 
scenario of using \cref{eq:problem}, we specifically selected
	high-dimensional datasets with $n$ larger than $m$.
	We conduct experiments with various $s$ to widely test
	the performance under different scenarios.
	In particular, we consider $s \in \{ \ceil{0.01m}, \ceil{0.05m},
	\ceil{0.1m}\}$ on all data except for the 
	largest dataset \webspam, for which we set $s \in \{ \ceil{0.001m}, 
	\ceil{0.005m}, \ceil{0.01m}\}$.	}
	The results of the experiment with the smallest $s$ are summarized in
\cref{fig:main}, and results of the other two settings of $s$ are in
\cref{sec:moreexp}.


\begin{table}[h] 
	\caption{Comparison of algorithms for \cref{eq:problem} to meet
	\cref{eq:opt} with $\hat \epsilon = 10^{-6}$, with
	\cref{eq:lr} and \cref{eq:ls_loss} and with sparsity levels $s_1=\ceil{0.01m}$ 
	and $s_2 = \ceil{0.05m}$ for all datasets
		except \webspam where $s_1 = \ceil{0.001m}$ and $s_2 = \ceil{0.005m}$.
		CPU: CPU time in seconds. GE: number of gradient
		evaluations. In one iteration, PG, APG, and APG+ needs one
		gradient evaluation , while PG-LL and PG-LL+ needs two.
		CG: number of Hessian-vector products in the PCG
		procedure for obtaining SSN steps.
		PA: prediction accuracy (for \cref{eq:lr}).
		MSE: mean-squared error (for \cref{eq:ls_loss}).
		Time with $*$ indicates that the algorithm is terminated after
		running $10000$ iterations without satisfying \cref{eq:opt}.}
	\label{tbl:exp}
	\begin{tabular}{@{}l@{}lrrrrrrrr@{}}
		\toprule 
		\multirow{3}*{Dataset} & \multirow{3}*{Method} & 
		\multicolumn{4}{c}{$s_1$} & 
		\multicolumn{4}{c}{$s_2$}\\
		\cmidrule(lr){3-6} 
		\cmidrule(lr){7-10} 
		& & \multicolumn{1}{c}{CPU} & \multicolumn{1}{c}{GE} & 
		\multicolumn{1}{c}{CG} &
		\multicolumn{1}{c}{PA} & \multicolumn{1}{c}{CPU} & 
		\multicolumn{1}{c}{GE} & 
		\multicolumn{1}{c}{CG} &
		\multicolumn{1}{c}{PA} \\
		\midrule 
		\multirow{5}*{\news}  & PG & $*$738.7 & 10000 & 0 & 0.877 & $*$728.9 & 10000 & 
0 & 
0.935 \\ 
  & APG & 151.7 & 1583 & 0 & 0.877 & 758.3 & 8428 & 0 & 0.923 \\ 
  & APG+ & \bftab 5.0 & 52 & 63 & 0.853 & \bftab 16.1 & 171 & 67 & 0.923 \\ 
  & PG-LL & 366.7 & 4682 & 0 & 0.873 & $*$1494.4 & 20000 & 0 & 0.922 \\ 
  & APG-LL+ & 6.6 & 152 & 88 & 0.854 & 29.2 & 417 & 89 & 0.920 \\ 
\midrule
		\multirow{5}*{\rcv}  & PG & $*$58.4 & 10000 & 0 & 0.937 & $*$72.7 & 10000 & 0 & 
0.951 \\ 
  & APG & 12.6 & 1120 & 0 & 0.935 & 82.4 & 6372 & 0 & 0.934 \\ 
  & APG+ & \bftab 0.3 & 21 & 42 & 0.931 & \bftab 2.4 & 192 & 138 & 0.940 \\ 
  & PG-LL & 22.2 & 3638 & 0 & 0.935 & 72.1 & 8738 & 0 & 0.929 \\ 
  & APG-LL+ & 0.6 & 99 & 49 & 0.930 & 4.9 & 626 & 236 & 0.939 \\ 
\midrule
		\multirow{5}*{\webspam}  & PG & $*$18660.1 & 10000 & 0 & 0.964 & $*$30776.2 & 
10000 & 0 & 0.978 \\ 
  & APG & 19683.4 & 7682 & 0 & 0.981 & 7722.4 & 2008 & 0 & 0.991 \\ 
  & APG+ & \bftab 248.3 & 75 & 88 & 0.969 & \bftab 695.4 & 164 & 57 & 0.991 \\ 
  & PG-LL & 9001.3 & 4720 & 0 & 0.972 & 10163.5 & 3098 & 0 & 0.990 \\ 
  & APG-LL+ & 447.3 & 264 & 92 & 0.965 & 837.3 & 294 & 90 & 0.992 \\ 
\midrule
		 & &
		\multicolumn{1}{c}{CPU} & \multicolumn{1}{c}{GE} & 
		\multicolumn{1}{c}{CG} &
		\multicolumn{1}{c}{MSE} &
		\multicolumn{1}{c}{CPU} & \multicolumn{1}{c}{GE} & 
		\multicolumn{1}{c}{CG} &
		\multicolumn{1}{c}{MSE}\\
		\midrule 
		\multirow{5}*{\El}  & PG & $*$2998.6 & 10000 & 0 & 0.167 & $*$3644.1 & 10000 & 
0 & 
0.161 \\ 
  & APG & 270.6 & 669 & 0 & 0.136 & 811.8 & 1757 & 0 & 0.133 \\ 
  & APG+ & \bftab 19.5 & 40 & 49 & 0.141 & \bftab 105.6 & 222 & 124 & 0.132 \\ 
  & PG-LL & $*$6049.8 & 20000 & 0 & 0.132 & 2696.0 & 7086 & 0 & 0.132 \\ 
  & APG-LL+ & 41.2 & 142 & 38 & 0.142 & 107.5 & 326 & 100 & 0.138 \\ 
\midrule
		\multirow{5}*{\Et}  & PG & $*$242.7 & 10000 & 0 & 0.152 & $*$666.9 & 10000 & 0 
& 
0.152 \\ 
  & APG & \bftab 1.3 & 14 & 0 & 0.154 & \bftab 3.3 & 33 & 0 & 0.153 \\ 
  & APG+ & \bftab 1.3 & 8 & 6 & 0.141 & \bftab 3.3 & 31 & 7 & 0.139 \\ 
  & PG-LL & 110.6 & 4440 & 0 & 0.152 & 304.8 & 4558 & 0 & 0.151 \\ 
  & APG-LL+ & 1.7 & 34 & 6 & 0.141 & 3.7 & 47 & 7 & 0.139 \\ 
\midrule

	\end{tabular} 
\end{table}

Evidently, the extrapolation procedure in APG provides a 
significant improvement in the running time compared with the base algorithm 
PG, and further
incorporating subspace identification as in APG+ results to a very fast algorithm 
that outperforms PG and APG by magnitudes.
Since the per-iteration cost of PG and APG are almost the same as
argued in \cref{sec:accelerate}, we note that the convergence of APG in terms
of iterations is also superior to that of PG.

We also report the required time and number of gradient evaluations (which
is the main computation at each iteration) for the algorithms
to drive \cref{eq:opt} below $\hat \epsilon = 10^{-6}$.
For PG, APG, and APG+, one gradient evaluation is needed per
iteration, so the number of gradient evaluations is equivalent to the
iteration count.
For PG-LL, two gradient evaluations are needed per iteration, so its
cost is twice of other methods.
We also report the prediction performance on the test data,
and we in particular use the test
accuracy for \cref{eq:lr} and the mean-squared error for
\cref{eq:ls_loss}.
Results for the two smaller $s$ are in \cref{tbl:exp} while 
that for
the largest $s$ is in \cref{sec:moreexp}.
It is clear from the results in
\cref{tbl:exp} that APG 
outperforms PG-LL for most of the test instances considered, while
APG+ is magnitudes faster than PG-LL.
When we equip PG-LL with our acceleration techniques
by replacing $\Tpg$ in \cref{alg:extrapolation,alg:identification}
with the algorithmic map defining PG-LL, we can further speed 
up PG-LL greatly as shown under the name APG-LL+ (see \cref{tbl:exp}).
We do not observe a method that consistently possesses the best
prediction performance, as this is mainly affected by which local
optima is found, while no algorithm is able to find the best local
optima among all candidates.
With no prediction performance degradation, we see that APG+ and
APG-LL+ reduce the time needed to solve \cref{eq:problem} to a level
significantly lower than that of the state of the art.

\Alert{In \cref{sec:varyingresiduals}, we demonstrate the effect on 
prediction performance when we vary the residual \cref{eq:opt} and illustrate 
that tight residual level is \Alertt{indeed} required to obtain better prediction. Comparisons 
with a greedy method is shown in \cref{sec:comparison_greedy}.}

\begin{figure}[tbh]
	\centering 
	\begin{tabular}{@{}c@{}c@{}c@{}c@{}}
		\multicolumn{4}{c}{Sparse regularized logistic regression}\\
		\begin{subfigure}[b]{0.24\textwidth} 
			\includegraphics[width=\linewidth]{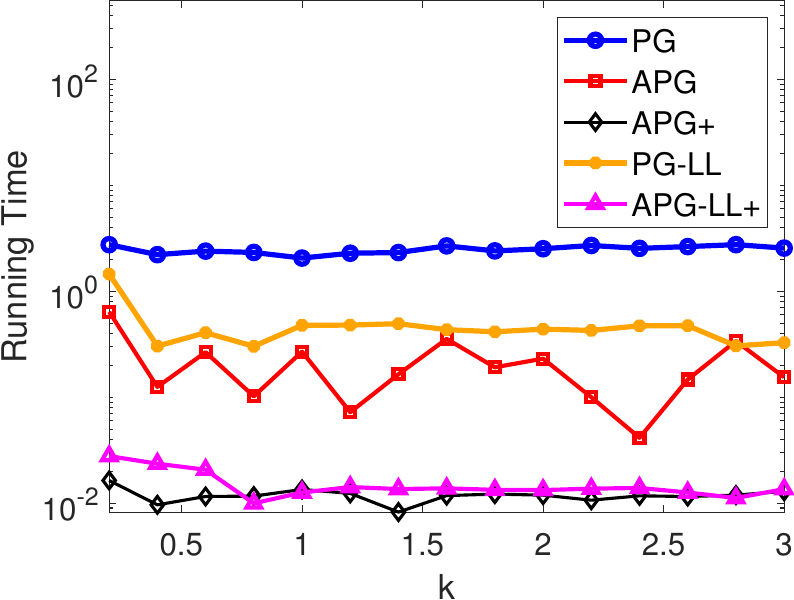}
			\caption{\colon} 
		\end{subfigure}& 
		\begin{subfigure}[b]{0.24\textwidth} 
			\includegraphics[width=\linewidth]{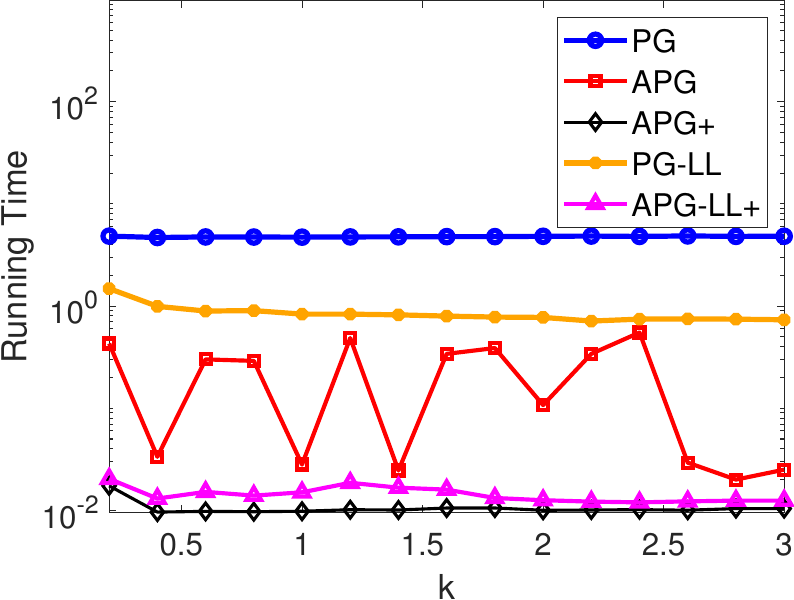}
			\caption{\duke} 
		\end{subfigure} &
		\begin{subfigure}[b]{0.24\textwidth} 
			\includegraphics[width=\linewidth]{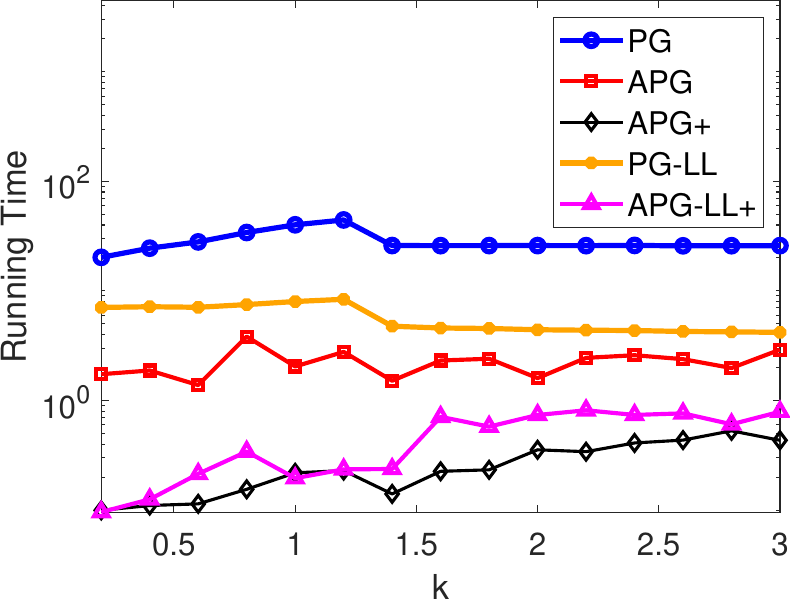}
			\caption{\gisette} 
		\end{subfigure}& 
		\begin{subfigure}[b]{0.24\textwidth} 
			\includegraphics[width=\linewidth]{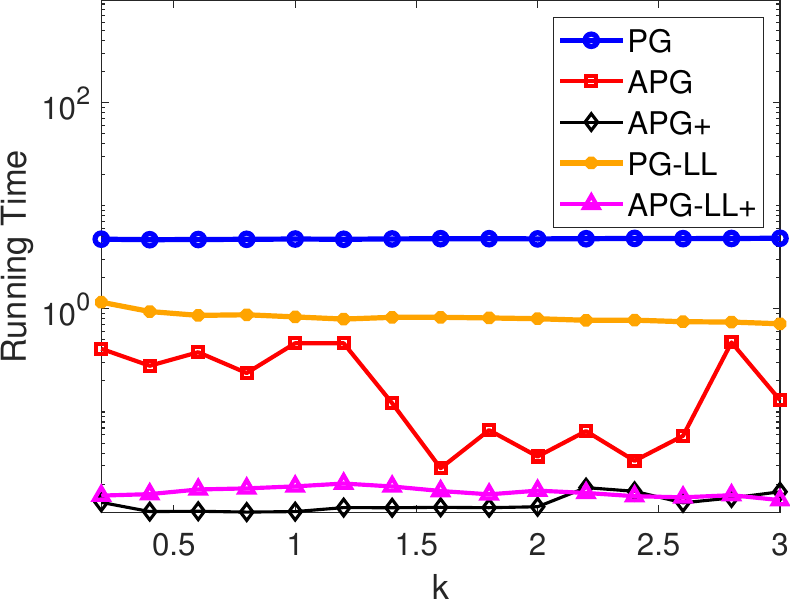}
			\caption{\leu} 
		\end{subfigure}\\
		\multicolumn{4}{c}{Sparse least squares regression}\\
		\begin{subfigure}[b]{0.24\textwidth} 
			\includegraphics[width=\linewidth]{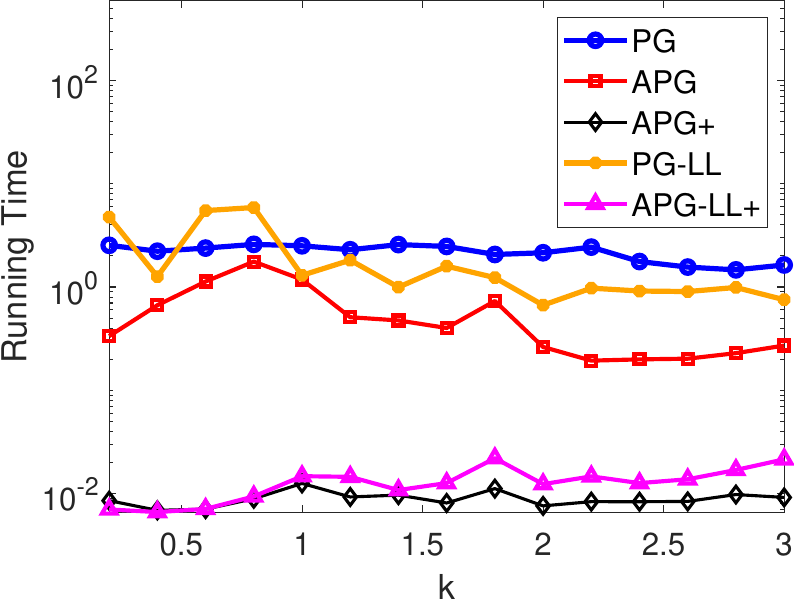}
			\caption{\colon} 
		\end{subfigure}& 
		\begin{subfigure}[b]{0.24\textwidth} 
			\includegraphics[width=\linewidth]{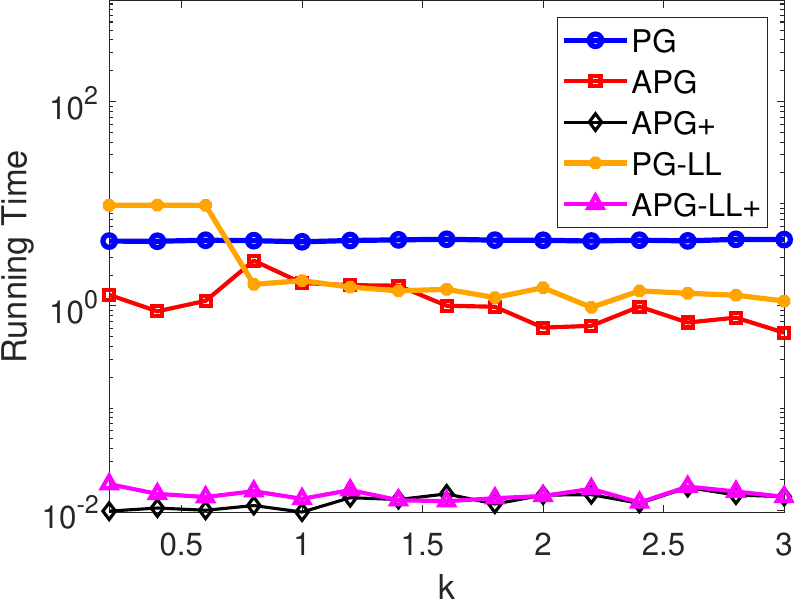}
			\caption{\duke} 
		\end{subfigure}  &
		\begin{subfigure}[b]{0.24\textwidth} 
			\includegraphics[width=\linewidth]{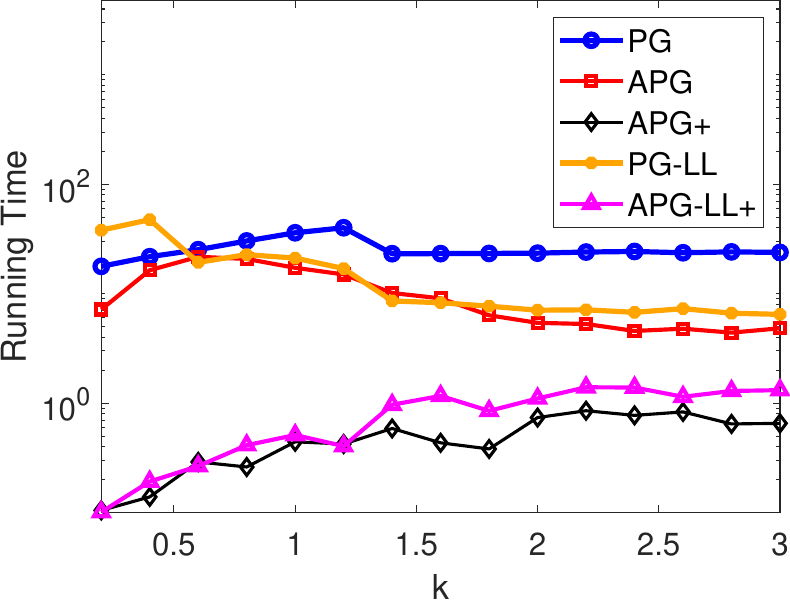}
			\caption{\gisette} 
		\end{subfigure}& 
		\begin{subfigure}[b]{0.24\textwidth} 
			\includegraphics[width=\linewidth]{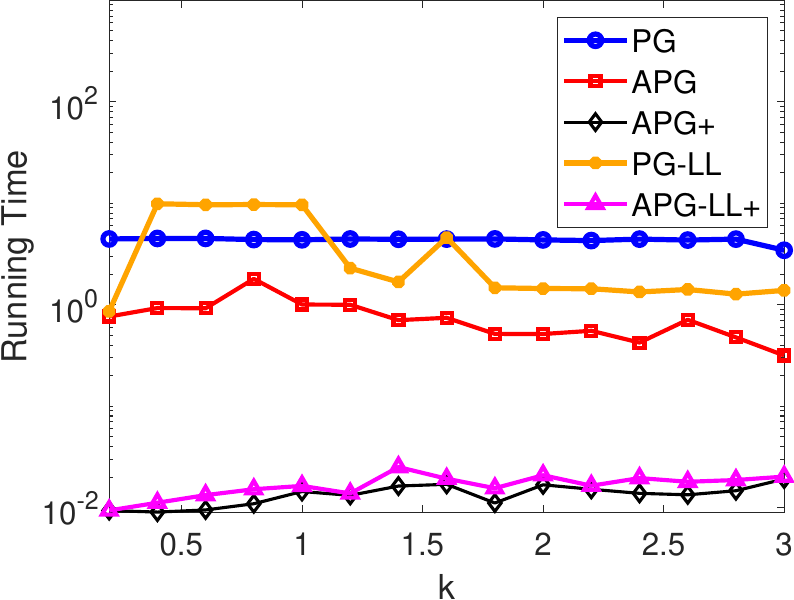}
			\caption{\leu} 
		\end{subfigure}
	\end{tabular} 
	\caption{\Alertt{Transition plots.  We present sparsity levels versus
	running time (in log scale). Top row: logistic loss. Bottom row:
least square loss.}}
	\label{fig:transitionplot} 
\end{figure}

\Alert{\paragraph{{Transition Plots.}} To demonstrate the behavior of the 
algorithm for increasing values of $s$, we fit the smaller datasets in \cref{tbl:data2} 
using logistic loss \cref{eq:lr} and least squares loss \cref{eq:ls_loss} for 
	varying 
	$s = \ceil{km}$, where $k=0.2,0.4,0.6,\dots,3$. \Alertt{The transition plots
	 are presented in
	 \cref{fig:transitionplot}.}} We note that the time is in log
	 scale.
	
	We can see clearly that APG+ and APG-LL+ are consistently magnitudes faster
	than the baseline PG method throughout all sparsity levels.
	On the other hand, the same-subspace extrapolation scheme of APG is
	consistently faster than PG and APG-LL and slower than the two Newton
	acceleration schemes, although the performance is sometimes closer to
	APG+/APG-LL+ while sometimes closer to PG.
	APG-LL tends to outperform PG in most situations as well, but in
	several cases when solving the least square problem, especially when
	$s$ is small, it can sometimes be slower than PG.
	Overall speaking, the results in the transition plots show that our
	proposed acceleration schemes are indeed effective for all sparsity
	levels tested.

\section{Conclusions}
\label{sec:conclusion}
In this work, we revisited the projected gradient algorithm for solving 
$\ell_0$-norm constrained 
optimization problems. Through a natural decomposition of the constraint set into subspaces
and the proven ability of the projected gradient method to identify a subspace that contains a
solution,
we further proposed effective acceleration schemes with provable
convergence speed improvements.
Experiments showed that our acceleration strategies improve
significantly both the convergence speed and the running time of the
original projected gradient algorithm,
and outperform the state of the art for $\ell_0$-norm constrained
problems by a huge margin.
We plan to extend our analysis and algorithm to the setting of a
nonconvex objective in the near future.

\Alert{\section*{Acknowledgments}
This work was supported in part by Academia Sinica Grand Challenge
Program Seed Grant No. AS-GCS-111-M05 and 
NSTC of R.O.C. grants 109-2222-E-001-003 and 111-2628-E-001-003. }

\bibliographystyle{plainnat}
\bibliography{../bibfile}

\ifdefined\neurips
\section*{Checklist}

\begin{enumerate}
	
	\item For all authors...
	\begin{enumerate}
		\item Do the main claims made in the abstract and introduction 
		accurately reflect the paper's contributions and scope?
		\answerYes{}
		\item Did you describe the limitations of your work?
		\answerYes{We only considered the case of $f$ convex, and we
			mentioned in the conclusions that we plan to extend our analysis
			to nonconvex $f$.}
		\item Did you discuss any potential negative societal impacts of your 
		work?
		\answerNA{This work has no negative societal impacts.}
		\item Have you read the ethics review guidelines and ensured that your 
		paper conforms to them?
		\answerYes{}
	\end{enumerate}
	
	\item If you are including theoretical results...
	\begin{enumerate}
		\item Did you state the full set of assumptions of all theoretical 
		results?
		\answerYes{}
		\item Did you include complete proofs of all theoretical results?
		\answerYes{}
	\end{enumerate}
	
	\item If you ran experiments...
	\begin{enumerate}
		\item Did you include the code, data, and instructions needed to 
		reproduce the main experimental results (either in the supplemental 
		material or as a URL)?
		\answerYes{In the supplemental materials.}
		\item Did you specify all the training details (e.g., data splits, 
		hyperparameters, how they were chosen)?
		\answerYes
		\item Did you report error bars (e.g., with respect to the random seed 
		after running experiments multiple times)?
		\answerNA{There is no randomness in our experiments or
			algorithms.}
		\item Did you include the total amount of compute and the type of 
		resources used (e.g., type of GPUs, internal cluster, or cloud 
		provider)?
		\answerYes
	\end{enumerate}
	
	\item If you are using existing assets (e.g., code, data, models) or 
	curating/releasing new assets...
	\begin{enumerate}
		\item If your work uses existing assets, did you cite the creators?
		\answerYes{The website of the source of data used in the experiments 
		are 
			provided.}
		\item Did you mention the license of the assets?
		\answerNA
		\item Did you include any new assets either in the supplemental 
		material or as a URL?
		\answerYes{Our codes are provided in the supplementary materials.}
		\item Did you discuss whether and how consent was obtained from people 
		whose data you're using/curating?
		\answerYes{The data used are publicly available.}
		\item Did you discuss whether the data you are using/curating contains 
		personally identifiable information or offensive content?
		\answerNA{The public data used do not contain such information/content.}
	\end{enumerate}
	
\end{enumerate}
\fi

\newpage
\appendix
\startcontents
\begin{spacing}{0}
\printcontents{atoc}{0}{\section*{Appendices}}
\end{spacing}

\section{Implementation Details for \cref{sec:identify}}
\label{subsec:ssn}
We first discuss our implementation for obtaining inexact SSN steps
described in \cref{sec:identify}.
Given any $w \in A_J$,
We use the notation $f_{J}(w_J) \coloneqq f(w)$ to denote the function
of considering only the coordinates of $w$ in $J$ as variables and
treating the remaining as constants equal to zero. For any $p \in \Re^s$, we
use $P^{-1}_{A_J}$ to denote the vector
$\hat p \in \Re^n$ with $\hat p_{J} = p$ and $\hat p_i = 0$ for $i
\notin J$.
Note that since $f_{J}$ is Lipschitz-continuously differentiable,
a generalized Hessian $\nabla^2 f_{J}$ always exists \cite{HUJ84a}.
When the set of generalized Hessian is not a singleton, we can pick
any element in the set.

In large-scale problems often faced in modern machine learning tasks,
$s$ can be large even if $s \ll n$, and thus forming the generalized
Hessian explicitly and inverting it could still be prohibitively
expensive even if we only consider the generalized Hessian in the
$s$-dimensional subspace.
Therefore, we resort to PCG that, given a preconditioner $M$,
iteratively uses the matrix-vector products $\nabla^2 f_J(w_J) v$
and $M^{-1}u$ for given vectors $u,v \in \Re^s$, which can be of much
lower cost especially if $M$ has certain structures to facilitate the
inverse.
Details of PCG can be found in, for instance,
\citet[Chapter~7]{NoceWrig06}.
The PCG approach provides an approximate solution to
\[
	p \approx \nabla^2 f_J(w_J)^{-1} \nabla f_J(w_J),
\]
or equivalently,
\begin{equation}
	p \approx \argmin_{\bar p} \left(Q_J(\bar p; w_J) \coloneqq
	\inner{\nabla f_J(w_J)}{\bar p} +
	\frac12 \inner{\bar p}{\nabla^2 f_J(w_J) \bar p}\right).
	\label{eq:pcg}
\end{equation}
In our implementation, inspired by the approach of \cite{CYH18a}, we
select the diagonal entries of $\nabla^2 f_J$ as our preconditioner
$M$, which provides better performance in our preliminary test over
using no preconditioner (or equivalently, taking $M$ as the identity
matrix).
As this choice of $M$ is a diagonal matrix, its inverse can be
computed efficiently in $O(s)$ time.

After obtaining $p$, given parameters $\beta, \sigma_2 \in (0,1)$,
we conduct a backtracking line search procedure to find the
largest nonnegative integer $i$ such that
\begin{equation}
	\label{eq:armijo}
	f_{J} \left( w_J + \beta^i p \right) \leq f_J \left( w_J
	\right) + \sigma_2 \beta^i \inner{\nabla f_J \left( w_J
	\right)}{p}
\end{equation}
and set the step size to $\alpha = \beta^i$.
Finally, the iterate is updated by
\[w_J \leftarrow w_J + \alpha p.
\]
If $\alpha$ is too small, or this decrease condition cannot be
satisfied even when $\beta^i$ is already extremely small, we
discard this SSN step and declare that this smooth optimization part
has failed in \cref{alg:identification}.

For the approximation criterion in \cref{eq:pcg}, let the $i$-th
iterate of PCG be $p^{(i)}$ and $Q_i \coloneqq Q_J(p^{(i)};w_J)$, we
follow \cite{GalL21a} to terminate PCG either when it reaches $s$
iterations (at which point theoretically it should have found the
exact solution of the right-hand side of \cref{eq:pcg}) or when the
$i$-th iterate satisfies $i \geq 1$ and
\begin{equation}
\frac{Q_i - Q_{i-1}}{\frac{Q_i}{i}} \leq \min\left\{ 0.5,
	\sqrt{\inner{\nabla f_J(w_J)}{M^{-1} \nabla f_J(w_J)}}\right\},
	\label{eq:terminate}
\end{equation}
where $Q_0 \coloneqq Q(0;w_J) = 0$.
It has been shown in \cite{GalL21a} that such a stopping condition
leads to $Q$-superlinear convergence to an optimum of $f_{J}$ when
$\nabla f_J$ is semismooth and $f$ is strongly convex.
In our case that alternates between such an SSN step and a PG step, we
will show that with \cref{eq:terminate}, the overall procedure will
enjoy superlinear convergence to $w^*$ if $\nabla f$ is semismooth
around $x^*$; see \cref{thm:superlinear} for more details.

One concern is that PCG only works when $\nabla^2 f_J$ is positive
definite, but our problem class only guarantees that it is positive
semidefinite.
To safeguard this issue, one can add a multiple of the identity to $\nabla^2
f_J$ as a damping term to make sure the quadratic term is always
positive definite.
A particularly useful way is to use $c \norm{\nabla f_J(w_J)}^\rho I$
as
the damping term for some $c>0$ and $\rho \in (0,1]$ in
\cref{eq:terminate}.
When $f_{J}$ satisfies a $q$-metric subregularity condition or an
error-bound condition,
this damping is known to produce a superlinear convergence rate of
order $(1+\rho)$ for a range of $q$ following the analysis in
\cite{MCY19a,MorYZZ20a}.
In \cref{thm:superlinear}, we do not consider any specific scenarios,
but just assume that the smooth optimization subroutine involved
itself has a superlinear convergence rate, and show that such a rate is
still retained when this subroutine is combined with our algorithm.
Therefore, discussions of various schemes including truncated Newton,
semismooth Newton, and damping, are all compatible with our general
framework to obtain superlinear convergence rates.

\section{Experimental settings}
\label{app:expdetails}
All experiments are conducted on a machine with 64GB memory and an Intel Xeon 
Silver 4208 CPU with 8 cores and 2.1GHz.
For all algorithms and all experiments, all cores are utilized.
The experiment environment runs Ubuntu 20.04 and MATLAB 2021b.
For experiments in \cref{sec:exp}, we use public data listed in \Alert{
\cref{tbl:data,tbl:data2}}.
\footnote{Downloaded from
	\url{http://www.csie.ntu.edu.tw/~cjlin/libsvmtools/datasets/}.} For the 
	datasets 
	that do not come with a test set, we manually do a $80/20$ split to obtain 
	a test set.

	\begin{table}
		\caption{Data statistics.}
		\label{tbl:data}
		\centering
		\setlength{\tabcolsep}{2pt}
		\begin{tabular}{@{}lrrrr@{}}
			\toprule
			Dataset & Loss & \#training & \#features  & \#test
			\\
			& & instances ($m$) & ($n$) & instances
			\\
			\midrule
			\news & \cref{eq:lr} & 15,997 & 1,355,191 & 3,999\\
			\rcv & \cref{eq:lr} & 20,242 & 47,236 & 677,399\\
			\webspam & \cref{eq:lr} & 280,000 & 16,609,143 & 70,000\\
			\midrule
			\El & \cref{eq:ls_loss} & 16,087 & 4,272,227 & 3,308\\
			\Et & \cref{eq:ls_loss} & 16,087 & 150,360 & 3,308\\
			\bottomrule
		\end{tabular}
	\end{table}

	\begin{table}[tbh]
	\caption{Data statistics for small datasets.}
	\label{tbl:data2}
	\centering
	\begin{tabular}{@{}lrrrr@{}}
		\toprule
		Dataset & Loss & \#training instances ($m$) 
		& \#features ($n$) & \#test instances
		\\
		\midrule
		\colon & \cref{eq:ls_loss}\&\cref{eq:lr} & 50 & 2,000 & 12\\
		\duke & \cref{eq:ls_loss}\&\cref{eq:lr} & 38 & 7,129 & 6\\
		\gisette & \cref{eq:ls_loss}\&\cref{eq:lr} & 1,000 & 5,000 &
		6,000 \\
		\leu & \cref{eq:ls_loss}\&\cref{eq:lr} & 38 & 7,129 &  34\\
		\bottomrule
	\end{tabular}
	\end{table}
The parameters used in our implementation are as follows.
We use $\mu = 10^{-10}$ in \cref{eq:lr}.
For \cref{alg:extrapolation}, $\sigma = 0.05$, $\eta = 0.5$, ${\epsilon
	= 10^{-20}}$, $\alpha_{\min}
= 1$, $\alpha_{\max} = 100$, $L$ is estimated using MATLAB's
\texttt{eigs} function to approximate the largest eigenvalue of $AA^\T$ with 
tolerance 
$10^{-3}$, and $\lambda = 0.999/L$.
In \cref{alg:identification}, we set $t = 1$ and $S = 5$, while for the
PCG and SSN subroutines, we set $\beta = 0.5$ and $\sigma_2 = 0.001$.

\section{Additional Experiments}
\label{sec:moreexp}
This section provides two sets of additional experiments.
We first present results of the datasets in \cref{sec:exp} with
different settings of $s$.
The second set of additional experiments are on some smaller datasets
that are often considered in existing works for the best subset
selection problem like \cite{BKM16}.

\subsection{Other settings of $s$}
\begin{figure}[tbh] 
	\centering 
	\begin{tabular}{cc}
		\begin{subfigure}[b]{0.31\textwidth} 
			\includegraphics[width=\linewidth]{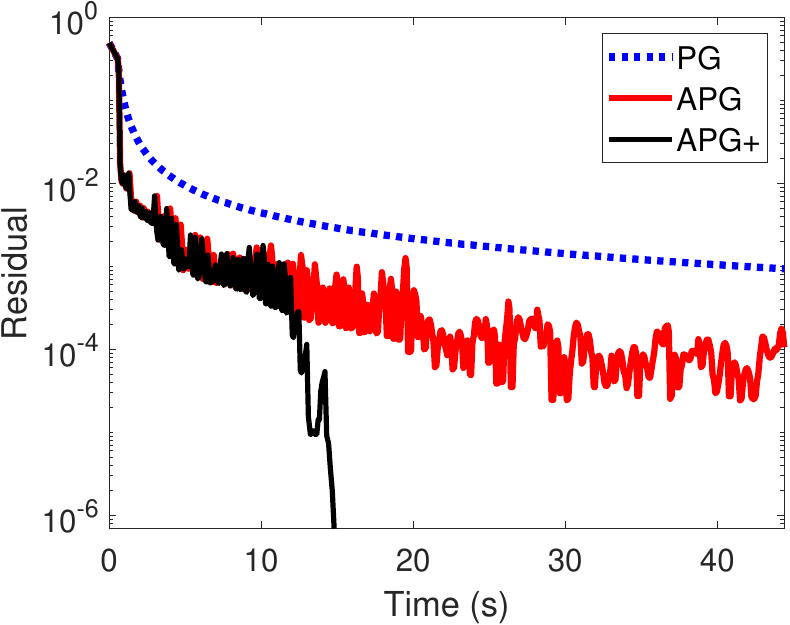}
			 
			\caption{\news, $s=\ceil {0.05m}$} 
		\end{subfigure}& 
		\begin{subfigure}[b]{0.31\textwidth} 
		\includegraphics[width=\linewidth]{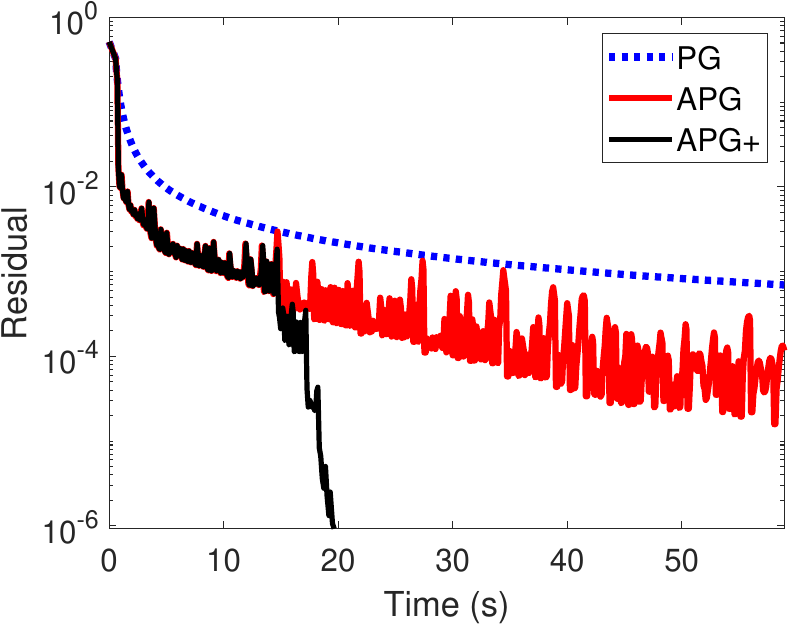}
		 
		\caption{\news, $s=\ceil {0.1m}$} 
		\end{subfigure}  \\
		\begin{subfigure}[b]{0.31\textwidth} 
			\includegraphics[width=\linewidth]{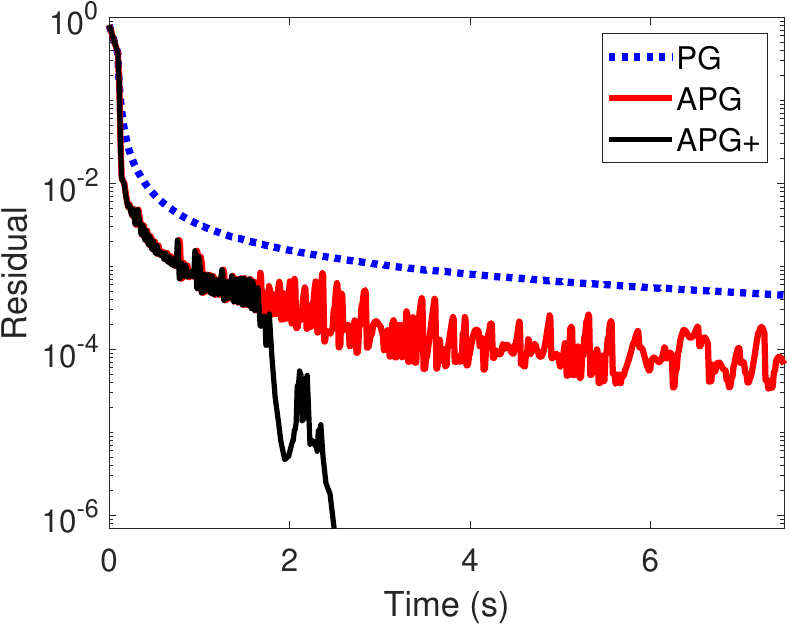}
			 
			\caption{\rcv, $s=\ceil {0.05m}$} 
		\end{subfigure}& 
		\begin{subfigure}[b]{0.31\textwidth} 
			\includegraphics[width=\linewidth]{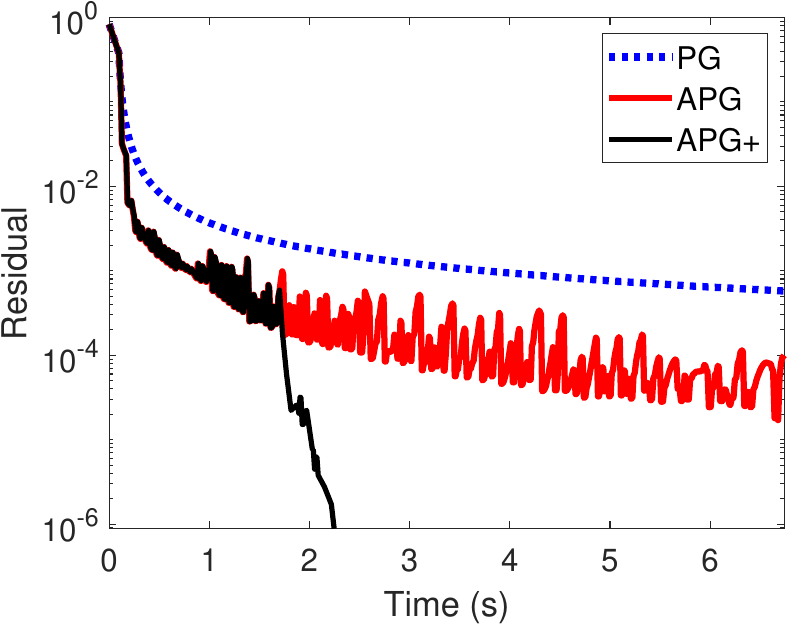}
			 
			\caption{\rcv, $s=\ceil {0.1m}$} 
		\end{subfigure} \\
		\begin{subfigure}[b]{0.31\textwidth} 
		\includegraphics[width=\linewidth]{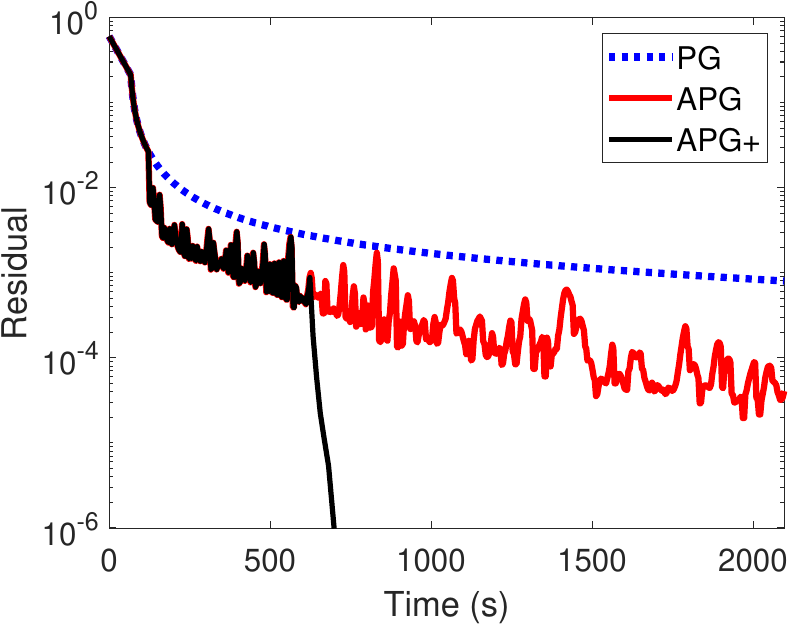}
		 
		\caption{\webspam, $s=\ceil {0.005m}$} 
		\end{subfigure}& 
		\begin{subfigure}[b]{0.31\textwidth} 
		\includegraphics[width=\linewidth]{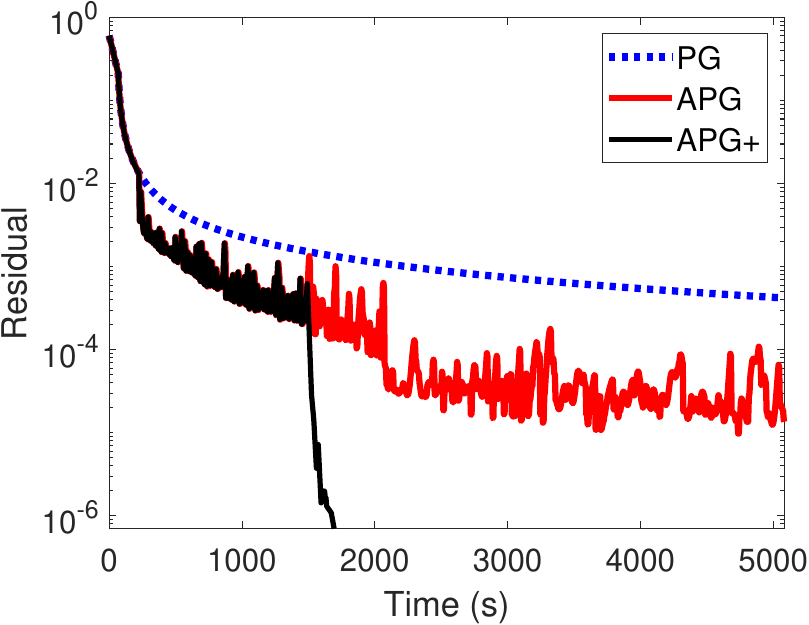}
		 
		\caption{\webspam, $s=\ceil {0.01m}$} 
		\end{subfigure} 
	\end{tabular} 
	\caption{Sparse regularized logistic loss regression.}
	\label{fig:logistic} 
\end{figure}

\begin{figure}[tbh]
	\centering 
	\begin{tabular}{cc}
		\begin{subfigure}[b]{0.31\textwidth} 
			\includegraphics[width=\linewidth]{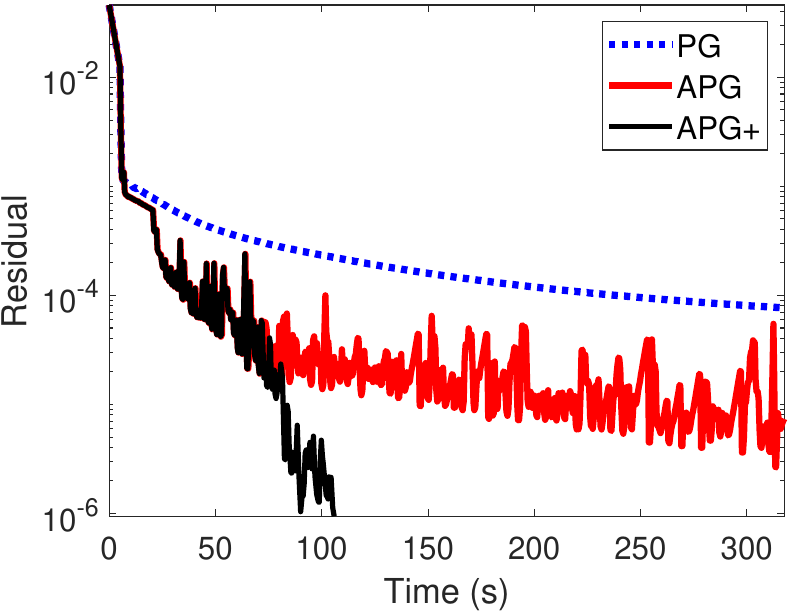}
			\caption{\El, $s=\ceil {0.05m}$} 
		\end{subfigure}& 
		\begin{subfigure}[b]{0.31\textwidth} 
			\includegraphics[width=\linewidth]{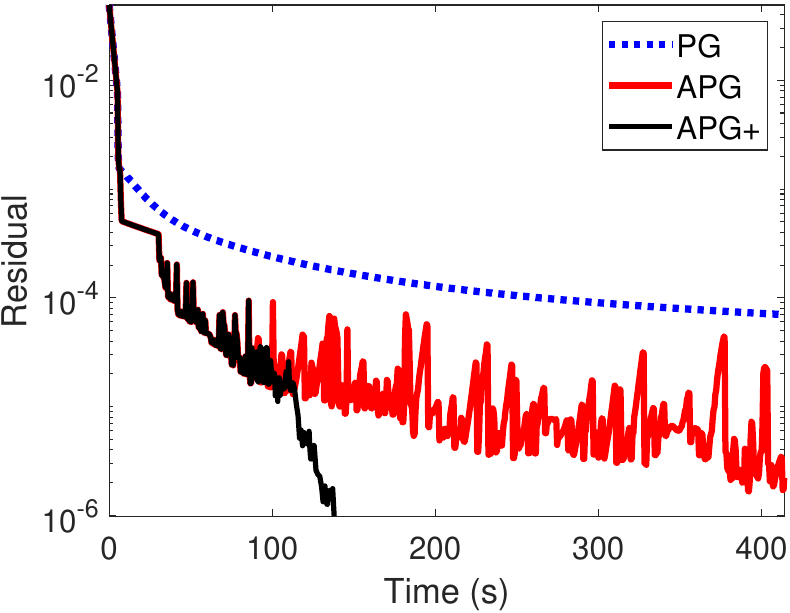}
			\caption{\El, $s=\ceil {0.1m}$} 
		\end{subfigure}  \\
		\begin{subfigure}[b]{0.31\textwidth} 
			\includegraphics[width=\linewidth]{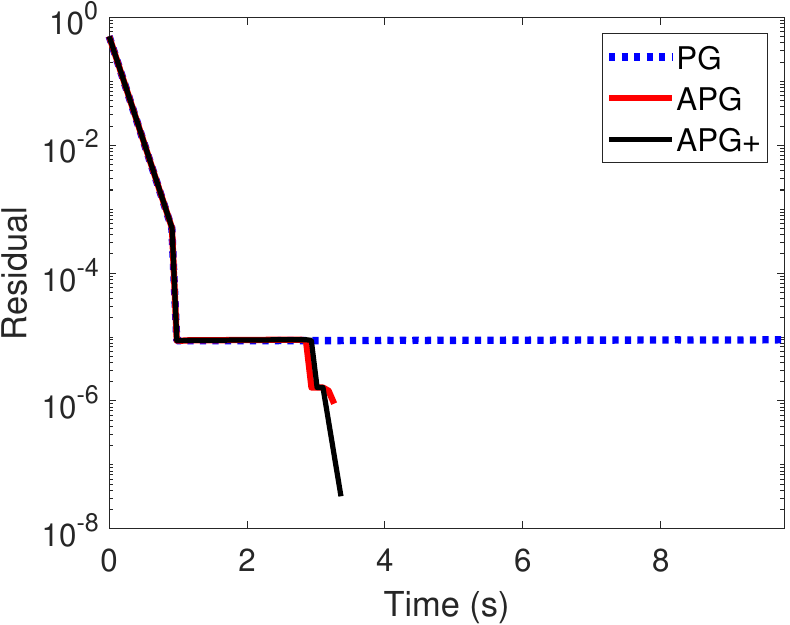}
			\caption{\Et, $s=\ceil {0.05m}$} 
		\end{subfigure}& 
		\begin{subfigure}[b]{0.31\textwidth} 
			\includegraphics[width=\linewidth]{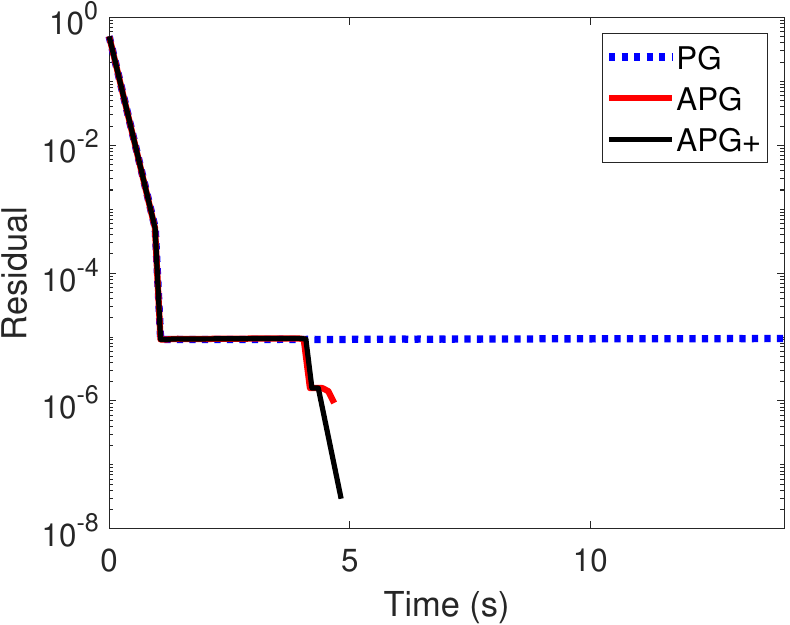}
			\caption{\Et, $s=\ceil {0.1m}$} 
		\end{subfigure}
	\end{tabular} 
	\caption{Sparse least squares regression.}
	\label{fig:leastsquare} 
\end{figure}

We present the other two settings of $s$ described in \cref{sec:exp}
in \cref{fig:logistic,fig:leastsquare}, and the continuation of \cref{tbl:exp} 
is presented in \cref{tbl:exp_continued}. Additional experiments with the 
setting of $s>m$ are presented in \cref{tbl:exp2,tbl:exp3}, which further 
exemplifies the benefits of our proposed acceleration strategies. 

Clearly, for the setting of $s_3$ as well as $s > m$, our acceleration
techniques continue to greatly improve upon existing methods in almost
all cases, with the only excepion being \webspam with $s > m$.
After a thorough check, we found that the reason is that in this
setting, due to the high dimensionality of $n$ and that many pieces of
$J \in \J_s$ can lead to a very low objective value, the subspaces in
which each $w^k$ lie change very frequently, so our extrapolation
barely take place.
This is a potential limit of our method, although in practice we
observe that for such easier datasets we probably can avoid this
problem by setting $s < m$, which would also make the problem much
easier to solve in general (note that with $s < m$, the prediction
performance on \webspam is not improving at all, suggesting that
indeed we do not need to consider the more difficult situation of $s <
m$).

We also observe that all for $s \geq m$ on \El, all accelerated
methods experience significantly larger MSE than the base PG method.
After a close examination, we find out that all such acceleration
methods provide much lower objective value than PG for the
minimization problem, indicating that this is merely due to
overfitting of the training data, and indeed PG is alway terminated
without reaching the prespecified stopping condition for these cases.
This indicates that the accelerated methods are actually performing well
from the optimization angle, and this overfitting issue is just a
matter of parameter selection.

For \Et, we see that for all settings of $s$,
identification does not show any additional time improvement in the
tables, while the figures clearly show that this is due to that this
step kicks in at a very late stage when the residual is already very
close to $\hat \epsilon$, and if we set $\hat \epsilon$ to a  smaller
value, we can expect observable running time difference between APG
and APG+.

\begin{table}[tbh] 
	\centering 
	\caption{Comparison of algorithms for \cref{eq:problem} to meet
		\cref{eq:opt} with $\hat \epsilon < 10^{-6}$, with
		\cref{eq:lr} and \cref{eq:ls_loss} and with sparsity level 
		$s_3 = 
		\ceil{0.1m}$ for all datasets 
		except \webspam where $s_3 = \ceil{0.01m}$.
		CPU: CPU time in seconds. GE: number of gradient
		evaluations. In one iteration, PG, APG, and APG+ needs one
		gradient evaluation , while PG-LL and PG-LL+ needs two.
		CG: number of Hessian-vector products in the PCG
		procedure for obtaining SSN steps.
		PA: prediction accuracy (for \cref{eq:lr}).
		MSE: mean-squared error (for \cref{eq:ls_loss}).
		Time with $*$ indicates that the algorithm is terminated after
		running $10000$ iterations without satisfying \cref{eq:opt}.
	}
	\label{tbl:exp_continued}
	\begin{tabular}{llrrrr}
		\toprule 
		\multirow{2}*{Dataset} & \multirow{3}*{Method} & 
		\multicolumn{4}{c}{$s_3$}   \\
		\cmidrule(lr){3-6} 
		& & \multicolumn{1}{c}{CPU} & \multicolumn{1}{c}{GE} & 
		\multicolumn{1}{c}{CG} &
		\multicolumn{1}{c}{PA}  \\
		\midrule 
		\multirow{5}*{\news}  & PG & $*$806.8 & 10000 & 0 & 0.947 \\ 
  & APG & 562.7 & 5972 & 0 & 0.927 \\ 
  & APG+ & \bftab 19.8 & 209 & 133 & 0.918 \\ 
  & PG-LL & 356.6 & 4578 & 0 & 0.930 \\ 
  & APG-LL+ & 23.2 & 463 & 223 & 0.918 \\ 
\midrule
		\multirow{5}*{\rcv}  & PG & $*$81.2 & 10000 & 0 & 0.953 \\ 
  & APG & 33.6 & 2556 & 0 & 0.943 \\ 
  & APG+ & \bftab 2.2 & 173 & 93 & 0.936 \\ 
  & PG-LL & 21.0 & 2292 & 0 & 0.940 \\ 
  & APG-LL+ & 4.5 & 542 & 106 & 0.933 \\ 
\midrule
		\multirow{5}*{\webspam}  & PG & $*$42487.5 & 10000 & 0 & 0.980 \\ 
  & APG & 11215.1 & 2242 & 0 & 0.993 \\ 
  & APG+ & 1664.7 & 313 & 83 & 0.994 \\ 
  & PG-LL & 14203.9 & 3176 & 0 & 0.992 \\ 
  & APG-LL+ & \bftab 1565.3 & 367 & 61 & 0.994 \\ 
\midrule
%
		\multirow{2}*{Dataset} & \multirow{3}*{Method} & 
		\multicolumn{4}{c}{$s_3$} \\
		\cmidrule(lr){3-6} 
		& & \multicolumn{1}{c}{CPU} & \multicolumn{1}{c}{GE} & 
		\multicolumn{1}{c}{CG} &
		\multicolumn{1}{c}{MSE} \\
		\midrule  
		\multirow{5}*{\El}  & PG & $*$4162.7 & 10000 & 0 & 0.160 \\ 
  & APG & 559.9 & 1084 & 0 & 0.142 \\ 
  & APG+ & \bftab 138.8 & 252 & 122 & 0.141 \\ 
  & PG-LL & 1996.5 & 4532 & 0 & 0.141 \\ 
  & APG-LL+ & 262.5 & 601 & 81 & 0.139 \\ 
\midrule
		\multirow{5}*{\Et}  & PG & $*$1086.3 & 10000 & 0 & 0.152 \\ 
  & APG & \bftab 4.7 & 33 & 0 & 0.153 \\ 
  & APG+ & \bftab 4.7 & 31 & 7 & 0.139 \\ 
  & PG-LL & 512.3 & 4602 & 0 & 0.151 \\ 
  & APG-LL+ & 8.6 & 75 & 7 & 0.139 \\ 
\midrule

	\end{tabular} 
\end{table}

\begin{table}[tbh] 
	\caption{Comparison of algorithms for \cref{eq:problem} to meet
		\cref{eq:opt} with $\hat \epsilon = 10^{-6}$, with
		\cref{eq:lr} and \cref{eq:ls_loss} and with sparsity levels $s\in \{ 
		m, \ceil{1.1m}\}$, i.e. $s\geq m$.
		CPU: CPU time in seconds. GE: number of gradient
		evaluations. In one iteration, PG, APG, and APG+ needs one
		gradient evaluation , while PG-LL and PG-LL+ needs two.
		CG: number of Hessian-vector products in the PCG
		procedure for obtaining SSN steps.
		PA: prediction accuracy (for \cref{eq:lr}).
		MSE: mean-squared error (for \cref{eq:ls_loss}).
		Time with $*$ indicates that the algorithm is terminated after
		running $10000$ iterations without satisfying \cref{eq:opt}. Time with 
		$\dagger$ indicates that the algorithm is terminated after exceeding 12 
		hours of running time without satisfying \cref{eq:opt}.}
	\label{tbl:exp2}
	\begin{tabular}{@{}l@{}lrrrrrrrr@{}}
		\toprule 
		\multirow{3}*{Dataset} & \multirow{3}*{Method} & 
		\multicolumn{4}{c}{$s=m$} & 
		\multicolumn{4}{c}{$s=\ceil{1.1m}$}\\
		\cmidrule(lr){3-6} 
		\cmidrule(lr){7-10} 
		& & \multicolumn{1}{c}{CPU} & \multicolumn{1}{c}{GE} & 
		\multicolumn{1}{c}{CG} &
		\multicolumn{1}{c}{PA} & \multicolumn{1}{c}{CPU} & 
		\multicolumn{1}{c}{GE} & 
		\multicolumn{1}{c}{CG} &
		\multicolumn{1}{c}{PA} \\
		\midrule 
		\multirow{5}*{\news}  & PG & $*$871.2 & 10000 & 0 & 0.963 & $*$869.0 & 10000 & 
0 & 0.963 \\ 
  & APG & 142.9 & 1482 & 0 & 0.962 & 191.6 & 1964 & 0 & 0.963 \\ 
  & APG+ & 58.3 & 619 & 24 & 0.966 & \bftab 64.8 & 684 & 26 & 0.969 \\ 
  & PG-LL & 156.5 & 1804 & 0 & 0.961 & 151.8 & 1778 & 0 & 0.961 \\ 
  & APG-LL+ & \bftab 47.7 & 555 & 9 & 0.958 & 66.2 & 743 & 9 & 0.955 \\ 
\midrule
		\multirow{5}*{\rcv}  & PG & $*$81.9 & 10000 & 0 & 0.959 & $*$82.6 & 10000 & 0 & 
0.959 \\ 
  & APG & 22.6 & 1859 & 0 & 0.956 & 18.5 & 1575 & 0 & 0.956 \\ 
  & APG+ & 5.1 & 468 & 47 & 0.952 & 5.4 & 524 & 47 & 0.953 \\ 
  & PG-LL & 15.7 & 1780 & 0 & 0.955 & 16.0 & 1784 & 0 & 0.955 \\ 
  & APG-LL+ & \bftab 4.9 & 539 & 29 & 0.951 & \bftab 4.6 & 490 & 38 & 0.951 \\ 
\midrule
		\multirow{5}*{\webspam}  & PG & $\dagger$43206.6 & 3902 & 0 & 0.977 & 
$\dagger$43203.1 & 
3870 & 0 & 0.977 \\ 
  & APG & $\dagger$43207.8 & 3866 & 0 & 0.985 & $\dagger$43202.0 & 3852 & 0 & 
  0.982 \\ 
  & APG+ & $\dagger$43207.5 & 3846 & 0 & 0.986 & $\dagger$43210.9 & 3879 & 0 & 
  0.983 \\ 
  & PG-LL & \bftab 35753.0 & 3190 & 0 & 0.995 & \bftab 35776.4 & 3190 & 0 & 0.995 \\ 
  & APG-LL+ & 36561.7 & 3190 & 0 & 0.995 & 36494.1 & 3190 & 0 & 0.995 \\ 
\midrule
		& &
		\multicolumn{1}{c}{CPU} & \multicolumn{1}{c}{GE} & 
		\multicolumn{1}{c}{CG} &
		\multicolumn{1}{c}{MSE} &
		\multicolumn{1}{c}{CPU} & \multicolumn{1}{c}{GE} & 
		\multicolumn{1}{c}{CG} &
		\multicolumn{1}{c}{MSE}\\
		\midrule 
		\multirow{5}*{\El}  & PG & $*$7039.7 & 10000 & 0 & 0.155 & $*$7172.7 & 10000 & 
0 & 
0.155 \\ 
  & APG & 4588.4 & 5819 & 0 & 0.207 & 5011.5 & 6275 & 0 & 0.213 \\ 
  & APG+ & \bftab 1050.0 & 1362 & 169 & 0.344 & 1320.0 & 1696 & 172 & 0.375 \\ 
  & PG-LL & 2261.4 & 3046 & 0 & 0.238 & 2292.0 & 3040 & 0 & 0.238 \\ 
  & APG-LL+ & 1220.1 & 1725 & 171 & 0.380 & \bftab 1282.3 & 1751 & 111 & 0.340 \\ 
\midrule
		\multirow{5}*{\Et}  & PG & $*$1821.0 & 10000 & 0 & 0.152 & $*$1819.9 & 10000 & 
0 & 
0.152 \\ 
  & APG & \bftab 67.4 & 353 & 0 & 0.155 & \bftab 69.0 & 363 & 0 & 0.155 \\ 
  & APG+ & 67.8 & 351 & 8 & 0.151 & 69.4 & 361 & 8 & 0.151 \\ 
  & PG-LL & 906.8 & 4832 & 0 & 0.151 & 909.6 & 4836 & 0 & 0.151 \\ 
  & APG-LL+ & 69.3 & 370 & 8 & 0.148 & 72.8 & 384 & 0 & 0.154 \\ 
\midrule
	\end{tabular} 
\end{table}

\begin{table}[tbh] 
	\caption{Comparison of algorithms for \cref{eq:problem} to meet
		\cref{eq:opt} with $\hat \epsilon = 10^{-6}$, with
		\cref{eq:lr} and \cref{eq:ls_loss} and with sparsity levels $s\in \{ 
		 \ceil{1.5m}, \ceil{2m}\}$, i.e. $s> m$.
		CPU: CPU time in seconds. GE: number of gradient
		evaluations. In one iteration, PG, APG, and APG+ needs one
		gradient evaluation , while PG-LL and PG-LL+ needs two.
		CG: number of Hessian-vector products in the PCG
		procedure for obtaining SSN steps.
		PA: prediction accuracy (for \cref{eq:lr}).
		MSE: mean-squared error (for \cref{eq:ls_loss}).
		Time with $*$ indicates that the algorithm is terminated after
		running $10000$ iterations without satisfying \cref{eq:opt}.
		Time with 
		$\dagger$ indicates that the algorithm is terminated after exceeding 12 
		hours of running time without satisfying \cref{eq:opt}.}
	\label{tbl:exp3}
	\begin{tabular}{@{}l@{}lrrrrrrrr@{}}
		\toprule 
		\multirow{3}*{Dataset} & \multirow{3}*{Method} & 
		\multicolumn{4}{c}{$s=\ceil{1.5m}$} & 
		\multicolumn{4}{c}{$s=\ceil{2m}$}\\
		\cmidrule(lr){3-6} 
		\cmidrule(lr){7-10} 
		& & \multicolumn{1}{c}{CPU} & \multicolumn{1}{c}{GE} & 
		\multicolumn{1}{c}{CG} &
		\multicolumn{1}{c}{PA} & \multicolumn{1}{c}{CPU} & 
		\multicolumn{1}{c}{GE} & 
		\multicolumn{1}{c}{CG} &
		\multicolumn{1}{c}{PA} \\
		\midrule 
		\multirow{5}*{\news}  & PG & $*$885.0 & 10000 & 0 & 0.964 & $*$904.5 & 10000 & 
0 & 
0.966 \\ 
  & APG & 208.2 & 2072 & 0 & 0.964 & 217.1 & 2170 & 0 & 0.964 \\ 
  & APG+ & 78.3 & 826 & 17 & 0.967 & 86.4 & 875 & 18 & 0.967 \\ 
  & PG-LL & 155.8 & 1736 & 0 & 0.963 & 153.1 & 1700 & 0 & 0.963 \\ 
  & APG-LL+ & \bftab 64.6 & 690 & 6 & 0.962 & \bftab 80.4 & 846 & 2 & 0.962 \\ 
\midrule
		\multirow{5}*{\rcv}  & PG & $*$84.8 & 10000 & 0 & 0.959 & $*$87.8 & 10000 & 0 & 
0.959 \\ 
  & APG & 19.6 & 1554 & 0 & 0.954 & 15.9 & 1317 & 0 & 0.955 \\ 
  & APG+ & \bftab 4.4 & 412 & 42 & 0.953 & \bftab 4.7 & 442 & 59 & 0.949 \\ 
  & PG-LL & 16.2 & 1784 & 0 & 0.956 & 16.6 & 1786 & 0 & 0.956 \\ 
  & APG-LL+ & 4.7 & 512 & 22 & 0.952 & 5.5 & 562 & 14 & 0.952 \\ 
\midrule
		\multirow{5}*{\webspam}  & PG & $\dagger$43201.3 & 3809 & 0 & 0.977 & 
$\dagger$43207.9 & 3807 & 0 & 0.977 \\ 
  & APG & $\dagger$43203.8 & 3815 & 0 & 0.978 & $\dagger$43201.7 & 3792 & 0 & 
  0.983 \\ 
  & APG+ & $\dagger$43202.7 & 3828 & 0 & 0.978 & $\dagger$43205.0 & 3783 & 0 & 
  0.983 \\ 
  & PG-LL & \bftab 36340.5 & 3190 & 0 & 0.995 & 36325.1 & 3190 & 0 & 0.995 \\ 
  & APG-LL+ & 36380.0 & 3190 & 0 & 0.995 & \bftab 31177.7 & 2716 & 24 & 0.995 \\ 
\midrule
		& &
		\multicolumn{1}{c}{CPU} & \multicolumn{1}{c}{GE} & 
		\multicolumn{1}{c}{CG} &
		\multicolumn{1}{c}{MSE} &
		\multicolumn{1}{c}{CPU} & \multicolumn{1}{c}{GE} & 
		\multicolumn{1}{c}{CG} &
		\multicolumn{1}{c}{MSE}\\
		\midrule 
		\multirow{5}*{\El}  & PG & $*$7617.3 & 10000 & 0 & 0.154 & $*$8003.0 & 10000 & 
0 & 
0.154 \\ 
  & APG & 5686.2 & 6781 & 0 & 0.209 & 6375.8 & 7300 & 0 & 0.201 \\ 
  & APG+ & 1697.5 & 2104 & 118 & 0.279 & 2098.5 & 2496 & 108 & 0.280 \\ 
  & PG-LL & 2398.6 & 3002 & 0 & 0.231 & 2460.1 & 2946 & 0 & 0.225 \\ 
  & APG-LL+ & \bftab 1362.7 & 1744 & 94 & 0.313 & \bftab 1672.1 & 2057 & 101 & 0.299 \\ 
\midrule
		\multirow{5}*{\Et}  & PG & $*$1870.2 & 10000 & 0 & 0.152 & $*$1894.5 & 10000 & 
0 & 
0.152 \\ 
  & APG & 89.0 & 465 & 0 & 0.155 & 97.0 & 500 & 0 & 0.155 \\ 
  & APG+ & 89.5 & 463 & 8 & 0.153 & 97.2 & 498 & 8 & 0.155 \\ 
  & PG-LL & 927.8 & 4848 & 0 & 0.151 & 948.6 & 4856 & 0 & 0.151 \\ 
  & APG-LL+ & \bftab 72.8 & 382 & 8 & 0.151 & \bftab 75.0 & 390 & 8 & 0.153 \\ 
\midrule
	\end{tabular} 
\end{table}

\subsection{Experiments with smaller datasets}
We now consider some other smaller datasets shown in \cref{tbl:data2},
which are also downloaded from the LIBSVM website.
Note that for \gisette, we interchanged the training and the test sets
to make $m < n$. For the setting of $s < m$, we consider $s \in \{
\ceil{0.01m},\ceil{0.05m},\ceil{0.1m},\ceil{0.5m}\}$,
while for the setting of $s \geq m$, we consider $s\in \{
m,\ceil{1.1m},\ceil{1.5m},2m\}$.
The results of least-square loss in \cref{eq:ls_loss} are shown in
\cref{tbl:exp_small_ls,tbl:exp_small_ls_others}, while the
results of the logistic loss in \cref{eq:lr} are shown in
\cref{tbl:exp_small_log,tbl:exp_small_log_others}.

We can clearly see from these results that our acceleration schemes
are also effective on smaller datasets to reduce the running time to
magnitudes shorter.
However, there are several cases that the running time is too short
such that the digits in the tables are unable to show difference
between APG and APG+.
We do not try to increase the number of digits in such cases, as the
running time is anyway already extremely short, and the difference
would not make much difference for problems that can be solved with
such high efficiency.

\begin{table}[tbh] 
	\centering 
	\caption{Comparison of algorithms for \cref{eq:problem} to meet
		\cref{eq:opt} with $\hat \epsilon = 10^{-6}$, with
		\cref{eq:ls_loss} with sparsity levels $s\in \{ 
			\ceil{0.01m},\ceil{0.05m},\ceil{0.1m},\ceil{0.5m}\}$.
		CPU: CPU time in seconds. GE: number of gradient
		evaluations. In one iteration, PG, APG, and APG+ needs one
		gradient evaluation , while PG-LL and PG-LL+ needs two.
		CG: number of Hessian-vector products in the PCG
		procedure for obtaining SSN steps.
		MSE: mean-squared error (for \cref{eq:ls_loss}).
		Time with $*$ indicates that the algorithm is terminated after
		running $10000$ iterations without satisfying \cref{eq:opt}.}
	\label{tbl:exp_small_ls}
	\begin{tabular}{@{}llrrrrrrrr@{}}
		\toprule 
		\multirow{3}*{Dataset} & \multirow{3}*{Method} & 
		\multicolumn{4}{c}{$s=\ceil{0.01m}$} & 
		\multicolumn{4}{c}{$s=\ceil{0.05m}$}\\
		\cmidrule(lr){3-6} 
		\cmidrule(lr){7-10} 
		& & \multicolumn{1}{c}{CPU} & \multicolumn{1}{c}{GE} & 
		\multicolumn{1}{c}{CG} &
		\multicolumn{1}{c}{MSE} & \multicolumn{1}{c}{CPU} & 
		\multicolumn{1}{c}{GE} & 
		\multicolumn{1}{c}{CG} &
		\multicolumn{1}{c}{MSE} \\
		\midrule 
		\multirow{5}*{\colon}  & PG & 0.44 & 2081 & 0 & 1.125 & 0.77 & 3162 & 0 & 0.645 \\ 
  & APG & \bftab 0.01 & 5 & 0 & 1.125 & 0.14 & 185 & 0 & 0.645 \\ 
  & APG+ & \bftab 0.01 & 5 & 0 & 1.125 & \bftab 0.01 & 8 & 3 & 0.645 \\ 
  & PG-LL & 0.21 & 664 & 0 & 1.125 & 0.19 & 768 & 0 & 0.645 \\ 
  & APG-LL+ & \bftab 0.01 & 10 & 0 & 1.125 & \bftab 0.01 & 16 & 2 & 0.646 \\ 
\midrule

		\multirow{5}*{\duke}  & PG & $*$4.07 & 10000 & 0 & 1.568 & $*$4.18 & 10000 & 0 
& 1.145 \\ 
  & APG & \bftab 0.01 & 5 & 0 & 1.581 & \bftab 0.01 & 8 & 0 & 1.140 \\ 
  & APG+ & \bftab 0.01 & 5 & 0 & 1.581 & \bftab 0.01 & 8 & 2 & 1.140 \\ 
  & PG-LL & 0.19 & 382 & 0 & 1.579 & 0.75 & 1514 & 0 & 1.141 \\ 
  & APG-LL+ & \bftab 0.01 & 10 & 0 & 1.581 & \bftab 0.01 & 16 & 2 & 1.140 \\ 
\midrule

		\multirow{5}*{\gisette}  & PG & $*$13.35 & 10000 & 0 & 0.465 & $*$14.70 & 10000 
& 0 & 0.304 \\ 
  & APG & 4.08 & 1526 & 0 & 0.464 & 3.63 & 1286 & 0 & 0.303 \\ 
  & APG+ & \bftab 0.07 & 11 & 18 & 0.464 & \bftab 0.08 & 13 & 38 & 0.334 \\ 
  & PG-LL & 2.45 & 1758 & 0 & 0.466 & $*$30.15 & 20000 & 0 & 0.268 \\ 
  & APG-LL+ & \bftab 0.07 & 32 & 18 & 0.464 & \bftab 0.08 & 37 & 23 & 0.337 \\ 
\midrule

		\multirow{5}*{\leu}  & PG & 3.67 & 8768 & 0 & 0.595 & 2.80 & 6726 & 0 & 0.566 \\ 
  & APG & \bftab 0.01 & 5 & 0 & 0.595 & \bftab 0.01 & 8 & 0 & 0.566 \\ 
  & APG+ & \bftab 0.01 & 5 & 0 & 0.595 & \bftab 0.01 & 8 & 2 & 0.566 \\ 
  & PG-LL & 0.14 & 302 & 0 & 0.595 & 0.77 & 1632 & 0 & 0.566 \\ 
  & APG-LL+ & \bftab 0.01 & 10 & 0 & 0.595 & \bftab 0.01 & 16 & 2 & 0.566 \\ 
\midrule

%
%
%
		\multirow{3}*{Dataset} & \multirow{3}*{Method} & 
		\multicolumn{4}{c}{$s=\ceil{0.1m}$} & 
		\multicolumn{4}{c}{$s=\ceil{0.5m}$}\\
		\cmidrule(lr){3-6} 
		\cmidrule(lr){7-10} 
		& & \multicolumn{1}{c}{CPU} & \multicolumn{1}{c}{GE} & 
		\multicolumn{1}{c}{CG} &
		\multicolumn{1}{c}{MSE} & \multicolumn{1}{c}{CPU} & 
		\multicolumn{1}{c}{GE} & 
		\multicolumn{1}{c}{CG} &
		\multicolumn{1}{c}{MSE} \\
		\midrule 
		\multirow{5}*{\colon}  & PG & 1.59 & 6951 & 0 & 0.652 & $*$2.51 & 10000 & 0 & 
1.855 \\ 
  & APG & 0.13 & 320 & 0 & 0.599 & 1.48 & 3268 & 0 & 2.461 \\ 
  & APG+ & \bftab 0.01 & 10 & 10 & 0.599 & \bftab 0.02 & 18 & 27 & 1.345 \\ 
  & PG-LL & 0.52 & 2990 & 0 & 0.656 & $*$6.26 & 20000 & 0 & 1.895 \\ 
  & APG-LL+ & \bftab 0.01 & 24 & 10 & 0.599 & \bftab 0.02 & 61 & 31 & 1.723 \\ 
\midrule

		\multirow{5}*{\duke}  & PG & $*$4.28 & 10000 & 0 & 0.860 & $*$4.18 & 10000 & 0 
& 
0.864 \\ 
  & APG & 0.02 & 23 & 0 & 0.882 & 1.26 & 1749 & 0 & 0.569 \\ 
  & APG+ & \bftab 0.01 & 9 & 5 & 0.882 & \bftab 0.01 & 11 & 14 & 1.060 \\ 
  & PG-LL & 0.84 & 1670 & 0 & 0.880 & 1.05 & 2284 & 0 & 1.089 \\ 
  & APG-LL+ & 0.02 & 19 & 5 & 0.882 & \bftab 0.01 & 28 & 14 & 1.060 \\ 
\midrule

		\multirow{5}*{\gisette}  & PG & $*$15.63 & 10000 & 0 & 0.243 & $*$23.63 & 10000 
& 0 & 0.220 \\ 
  & APG & 7.83 & 2697 & 0 & 0.212 & 14.97 & 3922 & 0 & 0.292 \\ 
  & APG+ & \bftab 0.08 & 13 & 41 & 0.260 & \bftab 0.22 & 43 & 108 & 0.253 \\ 
  & PG-LL & 5.37 & 3294 & 0 & 0.238 & $*$50.85 & 20000 & 0 & 0.364 \\ 
  & APG-LL+ & 0.09 & 54 & 40 & 0.259 & 0.48 & 271 & 149 & 0.258 \\ 
\midrule
		\multirow{5}*{\leu}  & PG & $*$4.35 & 10000 & 0 & 0.523 & $*$4.47 & 10000 & 0 & 
0.582 \\ 
  & APG & 0.99 & 1263 & 0 & 0.524 & 0.63 & 876 & 0 & 1.277 \\ 
  & APG+ & \bftab 0.01 & 9 & 5 & 0.524 & \bftab 0.01 & 10 & 16 & 1.676 \\ 
  & PG-LL & 0.93 & 1966 & 0 & 0.524 & $*$9.79 & 20000 & 0 & 1.226 \\ 
  & APG-LL+ & \bftab 0.01 & 19 & 5 & 0.524 & \bftab 0.01 & 30 & 16 & 1.676 \\ 
\midrule

	\end{tabular} 
\end{table}

\begin{table}[tbh] 
	\centering 
	\caption{Comparison of algorithms for \cref{eq:problem} to meet
		\cref{eq:opt} with $\hat \epsilon = 10^{-6}$, with
		\cref{eq:ls_loss} with sparsity levels $s\in \{
			m,\ceil{1.1m},\ceil{1.5m},2m\}$.
		CPU: CPU time in seconds. GE: number of gradient
		evaluations. In one iteration, PG, APG, and APG+ needs one
		gradient evaluation , while PG-LL and PG-LL+ needs two.
		CG: number of Hessian-vector products in the PCG
		procedure for obtaining SSN steps.
		MSE: mean-squared error (for \cref{eq:ls_loss}).
		Time with $*$ indicates that the algorithm is terminated after
		running $10000$ iterations without satisfying \cref{eq:opt}.}
	\label{tbl:exp_small_ls_others}
	\begin{tabular}{@{}llrrrrrrrr@{}}
		\toprule 
		\multirow{3}*{Dataset} & \multirow{3}*{Method} & 
		\multicolumn{4}{c}{$s=m$} & 
		\multicolumn{4}{c}{$s=\ceil{1.1m}$}\\
		\cmidrule(lr){3-6} 
		\cmidrule(lr){7-10} 
		& & \multicolumn{1}{c}{CPU} & \multicolumn{1}{c}{GE} & 
		\multicolumn{1}{c}{CG} &
		\multicolumn{1}{c}{MSE} & \multicolumn{1}{c}{CPU} & 
		\multicolumn{1}{c}{GE} & 
		\multicolumn{1}{c}{CG} &
		\multicolumn{1}{c}{MSE} \\
		\midrule 
		\multirow{5}*{\colon}  & PG & $*$2.45 & 10000 & 0 & 3.161 & $*$2.61 & 10000 & 0 
& 
3.048 \\ 
  & APG & 1.78 & 3272 & 0 & 2.994 & 0.52 & 1564 & 0 & 3.311 \\ 
  & APG+ & \bftab 0.03 & 39 & 169 & 10.643 & \bftab 0.02 & 26 & 100 & 6.470 \\ 
  & PG-LL & 1.08 & 5614 & 0 & 4.905 & 3.87 & 13116 & 0 & 4.510 \\ 
  & APG-LL+ & \bftab 0.03 & 192 & 128 & 12.519 & 0.03 & 137 & 101 & 5.062 \\ 
\midrule

		\multirow{5}*{\duke}  & PG & $*$4.18 & 10000 & 0 & 0.554 & $*$4.25 & 10000 & 0 
& 0.549 \\ 
  & APG & 1.42 & 1945 & 0 & 1.531 & 1.24 & 1631 & 0 & 0.505 \\ 
  & APG+ & \bftab 0.01 & 13 & 58 & 0.209 & \bftab 0.01 & 11 & 39 & 2.061 \\ 
  & PG-LL & 1.80 & 3782 & 0 & 0.635 & 1.47 & 3104 & 0 & 0.281 \\ 
  & APG-LL+ & 0.02 & 101 & 87 & 6.685 & \bftab 0.01 & 58 & 44 & 0.831 \\ 
\midrule

		\multirow{5}*{\gisette}  & PG & $*$35.83 & 10000 & 0 & 0.226 & $*$40.06 & 10000 
& 0 & 0.225 \\ 
  & APG & 16.50 & 3326 & 0 & 0.298 & 15.55 & 3000 & 0 & 0.308 \\ 
  & APG+ & \bftab 0.58 & 78 & 247 & 0.499 & \bftab 1.20 & 200 & 203 & 0.347 \\ 
  & PG-LL & 21.09 & 5364 & 0 & 0.367 & 18.91 & 4462 & 0 & 0.359 \\ 
  & APG-LL+ & 0.99 & 330 & 144 & 0.341 & 1.59 & 449 & 203 & 0.508 \\ 
\midrule
		\multirow{5}*{\leu}  & PG & $*$4.22 & 10000 & 0 & 0.649 & $*$4.30 & 10000 & 0 & 
0.703 \\ 
  & APG & 1.05 & 1202 & 0 & 0.967 & 1.06 & 1445 & 0 & 0.980 \\ 
  & APG+ & \bftab 0.02 & 15 & 104 & 4.722 & \bftab 0.02 & 17 & 98 & 9.202 \\ 
  & PG-LL & $*$9.64 & 20000 & 0 & 1.050 & 1.87 & 3968 & 0 & 1.718 \\ 
  & APG-LL+ & \bftab 0.02 & 118 & 104 & 4.722 & \bftab 0.02 & 112 & 98 & 9.202 \\ 
\midrule

%
%
%
		\multirow{3}*{Dataset} & \multirow{3}*{Method} & 
		\multicolumn{4}{c}{$s=\ceil{1.5m}$} & 
		\multicolumn{4}{c}{$s=2m$}\\
		\cmidrule(lr){3-6} 
		\cmidrule(lr){7-10} 
		& & \multicolumn{1}{c}{CPU} & \multicolumn{1}{c}{GE} & 
		\multicolumn{1}{c}{CG} &
		\multicolumn{1}{c}{MSE} & \multicolumn{1}{c}{CPU} & 
		\multicolumn{1}{c}{GE} & 
		\multicolumn{1}{c}{CG} &
		\multicolumn{1}{c}{MSE} \\
		\midrule 
		\multirow{5}*{\colon}  & PG & $*$2.71 & 10000 & 0 & 3.277 & 2.03 & 8195 & 0 & 
3.062 \\ 
  & APG & 0.27 & 805 & 0 & 2.836 & 0.22 & 614 & 0 & 3.016 \\ 
  & APG+ & \bftab 0.02 & 27 & 79 & 3.092 & \bftab 0.02 & 42 & 44 & 3.131 \\ 
  & PG-LL & 3.52 & 14734 & 0 & 3.016 & 0.66 & 3024 & 0 & 2.964 \\ 
  & APG-LL+ & 0.03 & 139 & 67 & 3.279 & 0.03 & 123 & 63 & 3.032 \\ 
\midrule
		\multirow{5}*{\duke}  & PG & $*$4.30 & 10000 & 0 & 0.392 & $*$4.28 & 10000 & 0 
& 0.333 \\ 
  & APG & 0.84 & 1162 & 0 & 0.271 & 0.51 & 687 & 0 & 0.301 \\ 
  & APG+ & \bftab 0.01 & 13 & 40 & 0.259 & \bftab 0.01 & 12 & 33 & 0.564 \\ 
  & PG-LL & 1.79 & 3740 & 0 & 0.373 & 1.49 & 2990 & 0 & 0.539 \\ 
  & APG-LL+ & 0.02 & 48 & 34 & 0.308 & 0.02 & 47 & 33 & 0.564 \\ 
\midrule
		\multirow{5}*{\gisette}  & PG & $*$23.09 & 10000 & 0 & 0.230 & $*$23.08 & 10000 
& 0 & 0.237 \\ 
  & APG & 7.40 & 2057 & 0 & 0.284 & 5.03 & 1443 & 0 & 0.282 \\ 
  & APG+ & \bftab 0.52 & 118 & 133 & 0.339 & 1.32 & 357 & 143 & 0.319 \\ 
  & PG-LL & 7.88 & 3242 & 0 & 0.347 & 6.98 & 2860 & 0 & 0.328 \\ 
  & APG-LL+ & 1.15 & 500 & 150 & 0.387 & \bftab 1.13 & 488 & 128 & 0.336 \\ 
\midrule
		\multirow{5}*{\leu}  & PG & $*$4.36 & 10000 & 0 & 0.730 & $*$4.30 & 10000 & 0 & 
0.704 \\ 
  & APG & 0.69 & 865 & 0 & 0.694 & 0.47 & 574 & 0 & 0.619 \\ 
  & APG+ & \bftab 0.02 & 19 & 53 & 1.369 & \bftab 0.01 & 12 & 42 & 0.944 \\ 
  & PG-LL & 1.92 & 3952 & 0 & 1.074 & 1.38 & 2862 & 0 & 0.714 \\ 
  & APG-LL+ & 0.03 & 81 & 49 & 1.230 & \bftab 0.01 & 43 & 29 & 0.926 \\ 
\midrule

	\end{tabular} 
\end{table}

\begin{table}[tbh] 
	\centering 
	\caption{Comparison of algorithms for \cref{eq:problem} to meet
		\cref{eq:opt} with $\hat \epsilon = 10^{-6}$, with
		\cref{eq:lr} with sparsity levels $s\in \{
			\ceil{0.01m},\ceil{0.05m}, \ceil{0.1m}, \ceil{0.5m}\}$.
		CPU: CPU time in seconds. GE: number of gradient
		evaluations. In one iteration, PG, APG, and APG+ needs one
		gradient evaluation , while PG-LL and PG-LL+ needs two.
		CG: number of Hessian-vector products in the PCG
		procedure for obtaining SSN steps.
		PA: prediction accuracy (for \cref{eq:lr}).
		Time with $*$ indicates that the algorithm is terminated after
		running $10000$ iterations without satisfying \cref{eq:opt}.}
	\label{tbl:exp_small_log}
	\begin{tabular}{@{}llrrrrrrrr@{}}
		\toprule 
		\multirow{3}*{Dataset} & \multirow{3}*{Method} & 
		\multicolumn{4}{c}{$s=\ceil{0.01m}$} & 
		\multicolumn{4}{c}{$s=\ceil{0.05m}$}\\
		\cmidrule(lr){3-6} 
		\cmidrule(lr){7-10} 
		& & \multicolumn{1}{c}{CPU} & \multicolumn{1}{c}{GE} & 
		\multicolumn{1}{c}{CG} &
		\multicolumn{1}{c}{PA} & \multicolumn{1}{c}{CPU} & 
		\multicolumn{1}{c}{GE} & 
		\multicolumn{1}{c}{CG} &
		\multicolumn{1}{c}{PA} \\
		\midrule 
		\multirow{5}*{\colon}  & PG & 1.69 & 7560 & 0 & 0.667 & $*$2.64 & 10000 & 0 & 
0.833 \\ 
  & APG & \bftab 0.01 & 12 & 0 & 0.667 & 0.93 & 1526 & 0 & 0.833 \\ 
  & APG+ & \bftab 0.01 & 9 & 2 & 0.667 & \bftab 0.01 & 10 & 8 & 0.833 \\ 
  & PG-LL & 0.06 & 226 & 0 & 0.667 & 0.40 & 1728 & 0 & 0.833 \\ 
  & APG-LL+ & \bftab 0.01 & 16 & 2 & 0.667 & 0.02 & 22 & 8 & 0.833 \\ 
\midrule

		\multirow{5}*{\duke}  & PG & $*$4.52 & 10000 & 0 & 0.000 & $*$4.72 & 10000 & 0 
& 0.500 \\ 
  & APG & \bftab 0.01 & 11 & 0 & 0.000 & \bftab 0.01 & 12 & 0 & 0.500 \\ 
  & APG+ & \bftab 0.01 & 8 & 1 & 0.000 & \bftab 0.01 & 8 & 2 & 0.500 \\ 
  & PG-LL & 0.47 & 958 & 0 & 0.000 & 0.58 & 1038 & 0 & 0.500 \\ 
  & APG-LL+ & \bftab 0.01 & 15 & 1 & 0.000 & \bftab 0.01 & 16 & 2 & 0.500 \\ 
\midrule

		\multirow{5}*{\gisette}  & PG & $*$16.11 & 10000 & 0 & 0.839 & $*$17.50 & 10000 
& 0 & 0.912 \\ 
  & APG & 11.62 & 3600 & 0 & 0.888 & 6.81 & 1980 & 0 & 0.916 \\ 
  & APG+ & \bftab 0.07 & 10 & 15 & 0.851 & \bftab 0.08 & 11 & 23 & 0.900 \\ 
  & PG-LL & $*$34.11 & 20000 & 0 & 0.863 & 5.34 & 2780 & 0 & 0.916 \\ 
  & APG-LL+ & 0.08 & 29 & 15 & 0.851 & 0.09 & 38 & 24 & 0.898 \\ 
\midrule
		\multirow{5}*{\leu}  & PG & $*$4.60 & 10000 & 0 & 0.824 & $*$4.81 & 10000 & 0 & 
0.882 \\ 
  & APG & \bftab 0.01 & 9 & 0 & 0.824 & 0.05 & 46 & 0 & 0.853 \\ 
  & APG+ & \bftab 0.01 & 9 & 2 & 0.824 & \bftab 0.01 & 10 & 6 & 0.853 \\ 
  & PG-LL & 0.81 & 1612 & 0 & 0.824 & 1.23 & 2278 & 0 & 0.853 \\ 
  & APG-LL+ & \bftab 0.01 & 16 & 2 & 0.824 & \bftab 0.01 & 20 & 6 & 0.853 \\ 
\midrule

%
%
%
		\multirow{3}*{Dataset} & \multirow{3}*{Method} & 
		\multicolumn{4}{c}{$s=\ceil{0.1m}$} & 
		\multicolumn{4}{c}{$s=\ceil{0.5m}$}\\
		\cmidrule(lr){3-6} 
		\cmidrule(lr){7-10} 
		& & \multicolumn{1}{c}{CPU} & \multicolumn{1}{c}{GE} & 
		\multicolumn{1}{c}{CG} &
		\multicolumn{1}{c}{PA} & \multicolumn{1}{c}{CPU} & 
		\multicolumn{1}{c}{GE} & 
		\multicolumn{1}{c}{CG} &
		\multicolumn{1}{c}{PA} \\
		\midrule 
		\multirow{5}*{\colon}  & PG & $*$2.48 & 10000 & 0 & 0.833 & $*$2.21 & 10000 & 0 
& 0.833 \\ 
  & APG & 2.51 & 3651 & 0 & 0.833 & 0.20 & 255 & 0 & 0.833 \\ 
  & APG+ & \bftab 0.02 & 11 & 12 & 0.833 & \bftab 0.02 & 14 & 55 & 0.833 \\ 
  & PG-LL & $*$5.79 & 20000 & 0 & 0.833 & 0.29 & 1376 & 0 & 0.833 \\ 
  & APG-LL+ & \bftab 0.02 & 29 & 15 & 0.833 & \bftab 0.02 & 58 & 44 & 0.833 \\ 
\midrule

		\multirow{5}*{\duke}  & PG & $*$4.71 & 10000 & 0 & 0.750 & $*$4.81 & 10000 & 0 
& 0.750 \\ 
  & APG & 0.03 & 33 & 0 & 0.500 & 0.43 & 431 & 0 & 0.750 \\ 
  & APG+ & \bftab 0.01 & 10 & 8 & 0.500 & \bftab 0.01 & 11 & 11 & 0.750 \\ 
  & PG-LL & 1.22 & 2352 & 0 & 0.500 & 0.99 & 1800 & 0 & 0.750 \\ 
  & APG-LL+ & \bftab 0.01 & 22 & 8 & 0.500 & 0.02 & 25 & 11 & 0.750 \\ 
\midrule

		\multirow{5}*{\gisette}  & PG & $*$18.31 & 10000 & 0 & 0.928 & $*$27.24 & 10000 
& 0 & 0.951 \\ 
  & APG & 4.79 & 1329 & 0 & 0.934 & 2.62 & 608 & 0 & 0.960 \\ 
  & APG+ & \bftab 0.09 & 12 & 42 & 0.923 & 0.23 & 48 & 34 & 0.956 \\ 
  & PG-LL & 6.40 & 3190 & 0 & 0.936 & 7.18 & 2324 & 0 & 0.955 \\ 
  & APG-LL+ & \bftab 0.09 & 56 & 42 & 0.923 & \bftab 0.14 & 59 & 45 & 0.955 \\ 
\midrule

		\multirow{5}*{\leu}  & PG & $*$4.79 & 10000 & 0 & 0.912 & $*$4.81 & 10000 & 0 & 
0.912 \\ 
  & APG & 0.54 & 522 & 0 & 0.824 & 0.29 & 304 & 0 & 0.882 \\ 
  & APG+ & \bftab 0.01 & 11 & 10 & 0.853 & \bftab 0.01 & 11 & 9 & 0.912 \\ 
  & PG-LL & 1.50 & 2726 & 0 & 0.912 & 0.92 & 1698 & 0 & 0.912 \\ 
  & APG-LL+ & 0.02 & 21 & 7 & 0.853 & 0.02 & 24 & 10 & 0.912 \\ 
\midrule
	\end{tabular} 
\end{table}

\begin{table}[tbh] 
	\centering 
	\caption{Comparison of algorithms for \cref{eq:problem} to meet
		\cref{eq:opt} with $\hat \epsilon = 10^{-6}$, with
		\cref{eq:lr} with sparsity levels $s\in \{ m,\ceil{1.1m},
		\ceil{1.5m}, 2m\}$.
		CPU: CPU time in seconds. GE: number of gradient
		evaluations. In one iteration, PG, APG, and APG+ needs one
		gradient evaluation , while PG-LL and PG-LL+ needs two.
		CG: number of Hessian-vector products in the PCG
		procedure for obtaining SSN steps.
		PA: prediction accuracy (for \cref{eq:lr}).
		Time with $*$ indicates that the algorithm is terminated after
		running $10000$ iterations without satisfying \cref{eq:opt}.}
	\label{tbl:exp_small_log_others}
	\begin{tabular}{@{}llrrrrrrrr@{}}
		\toprule 
		\multirow{3}*{Dataset} & \multirow{3}*{Method} & 
		\multicolumn{4}{c}{$s=m$} & 
		\multicolumn{4}{c}{$s=\ceil{1.1m}$}\\
		\cmidrule(lr){3-6} 
		\cmidrule(lr){7-10} 
		& & \multicolumn{1}{c}{CPU} & \multicolumn{1}{c}{GE} & 
		\multicolumn{1}{c}{CG} &
		\multicolumn{1}{c}{PA} & \multicolumn{1}{c}{CPU} & 
		\multicolumn{1}{c}{GE} & 
		\multicolumn{1}{c}{CG} &
		\multicolumn{1}{c}{PA} \\
		\midrule 
		\multirow{5}*{\colon}  & PG & $*$2.46 & 10000 & 0 & 0.833 & $*$2.33 & 10000 & 0 
& 0.833 \\ 
  & APG & 0.17 & 369 & 0 & 0.833 & 0.09 & 118 & 0 & 0.833 \\ 
  & APG+ & \bftab 0.02 & 14 & 34 & 0.833 & \bftab 0.02 & 19 & 35 & 0.833 \\ 
  & PG-LL & 0.32 & 1336 & 0 & 0.750 & 0.34 & 1364 & 0 & 0.750 \\ 
  & APG-LL+ & \bftab 0.02 & 73 & 37 & 0.833 & \bftab 0.02 & 53 & 39 & 0.833 \\ 
\midrule

		\multirow{5}*{\duke}  & PG & $*$4.82 & 10000 & 0 & 0.750 & $*$4.81 & 10000 & 0 
& 0.750 \\ 
  & APG & 0.03 & 31 & 0 & 0.750 & 0.33 & 351 & 0 & 0.500 \\ 
  & APG+ & \bftab 0.01 & 11 & 9 & 0.750 & \bftab 0.01 & 11 & 9 & 0.750 \\ 
  & PG-LL & 0.81 & 1558 & 0 & 0.750 & 0.80 & 1540 & 0 & 0.750 \\ 
  & APG-LL+ & 0.02 & 23 & 9 & 0.750 & 0.02 & 23 & 9 & 0.750 \\ 
\midrule
		\multirow{5}*{\gisette}  & PG & $*$40.08 & 10000 & 0 & 0.956 & $*$42.05 & 10000 
& 0 & 0.957 \\ 
  & APG & 1.96 & 366 & 0 & 0.960 & 3.94 & 668 & 0 & 0.957 \\ 
  & APG+ & 0.61 & 102 & 28 & 0.962 & \bftab 0.45 & 69 & 35 & 0.962 \\ 
  & PG-LL & 8.03 & 1862 & 0 & 0.962 & 8.19 & 1794 & 0 & 0.961 \\ 
  & APG-LL+ & \bftab 0.25 & 57 & 33 & 0.959 & 0.92 & 192 & 28 & 0.962 \\ 
\midrule
		\multirow{5}*{\leu}  & PG & $*$4.79 & 10000 & 0 & 0.912 & $*$4.75 & 10000 & 0 & 
0.912 \\ 
  & APG & 0.43 & 439 & 0 & 0.971 & 0.25 & 265 & 0 & 0.912 \\ 
  & APG+ & \bftab 0.01 & 11 & 8 & 0.912 & \bftab 0.02 & 11 & 8 & 0.941 \\ 
  & PG-LL & 0.83 & 1538 & 0 & 0.912 & 0.80 & 1516 & 0 & 0.912 \\ 
  & APG-LL+ & 0.02 & 22 & 8 & 0.912 & \bftab 0.02 & 22 & 8 & 0.941 \\ 
\midrule

%
%
%
		\multirow{3}*{Dataset} & \multirow{3}*{Method} & 
		\multicolumn{4}{c}{$s=\ceil{1.5m}$} & 
		\multicolumn{4}{c}{$s=2m$}\\
		\cmidrule(lr){3-6} 
		\cmidrule(lr){7-10} 
		& & \multicolumn{1}{c}{CPU} & \multicolumn{1}{c}{GE} & 
		\multicolumn{1}{c}{CG} &
		\multicolumn{1}{c}{PA} & \multicolumn{1}{c}{CPU} & 
		\multicolumn{1}{c}{GE} & 
		\multicolumn{1}{c}{CG} &
		\multicolumn{1}{c}{PA} \\
		\midrule 
		\multirow{5}*{\colon}  & PG & $*$2.69 & 10000 & 0 & 0.833 & $*$2.39 & 10000 & 0 
& 0.833 \\ 
  & APG & 0.09 & 219 & 0 & 0.667 & 0.23 & 304 & 0 & 0.833 \\ 
  & APG+ & \bftab 0.02 & 28 & 25 & 0.833 & 0.03 & 40 & 20 & 0.833 \\ 
  & PG-LL & 0.37 & 1282 & 0 & 0.833 & 0.30 & 1298 & 0 & 0.833 \\ 
  & APG-LL+ & \bftab 0.02 & 47 & 33 & 0.750 & \bftab 0.02 & 55 & 25 & 0.750 \\ 
\midrule

		\multirow{5}*{\duke}  & PG & $*$4.80 & 10000 & 0 & 0.750 & $*$4.83 & 10000 & 0 
& 0.750 \\ 
  & APG & 0.04 & 43 & 0 & 0.750 & 0.11 & 104 & 0 & 0.750 \\ 
  & APG+ & \bftab 0.01 & 11 & 9 & 0.750 & \bftab 0.01 & 12 & 11 & 0.750 \\ 
  & PG-LL & 0.84 & 1496 & 0 & 0.750 & 0.81 & 1416 & 0 & 0.750 \\ 
  & APG-LL+ & 0.02 & 30 & 6 & 0.750 & 0.02 & 25 & 11 & 0.750 \\ 
\midrule
		\multirow{5}*{\gisette}  & PG & $*$26.10 & 10000 & 0 & 0.956 & $*$26.29 & 10000 
& 0 & 0.956 \\ 
  & APG & 2.55 & 583 & 0 & 0.958 & 1.53 & 365 & 0 & 0.944 \\ 
  & APG+ & \bftab 0.48 & 107 & 25 & 0.960 & \bftab 0.69 & 159 & 23 & 0.956 \\ 
  & PG-LL & 4.66 & 1658 & 0 & 0.958 & 4.42 & 1568 & 0 & 0.959 \\ 
  & APG-LL+ & 0.63 & 193 & 15 & 0.960 & 0.74 & 242 & 12 & 0.958 \\ 
\midrule
		\multirow{5}*{\leu}  & PG & $*$4.67 & 10000 & 0 & 0.912 & $*$4.75 & 10000 & 0 & 
0.912 \\ 
  & APG & 0.44 & 435 & 0 & 0.971 & 0.03 & 36 & 0 & 0.941 \\ 
  & APG+ & \bftab 0.02 & 12 & 11 & 0.941 & \bftab 0.01 & 12 & 11 & 0.941 \\ 
  & PG-LL & 0.81 & 1480 & 0 & 0.912 & 0.78 & 1406 & 0 & 0.941 \\ 
  & APG-LL+ & \bftab 0.02 & 25 & 11 & 0.941 & 0.02 & 25 & 11 & 0.941 \\ 
\midrule

	\end{tabular} 
\end{table}

\clearpage


\subsection{\alert{Prediction accuracy for varying residuals}}
\label{sec:varyingresiduals}
\begin{figure}[tbh] 
	\centering 
%
%
%
		\begin{tabular}{cc}
		\begin{subfigure}[b]{0.45\textwidth} 
			
\includegraphics[width=.86\linewidth]{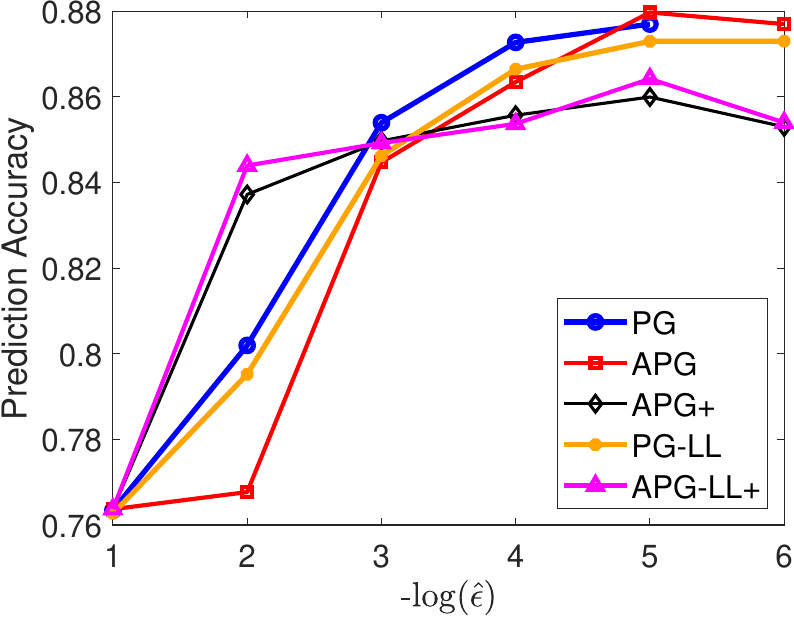}
			\caption{\news} 
		\end{subfigure}& 
		\begin{subfigure}[b]{0.45\textwidth} 
			
\includegraphics[width=.86\linewidth]{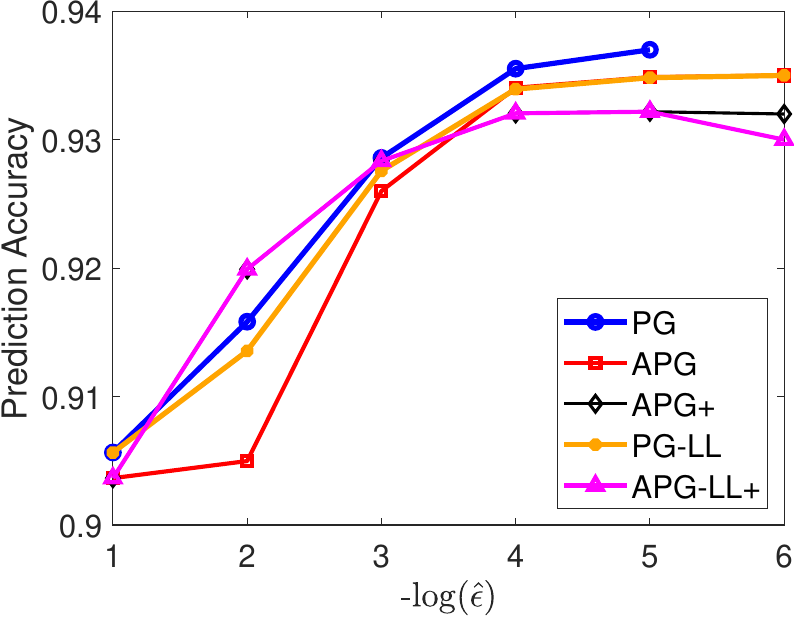}
			\caption{\rcv} 
		\end{subfigure}\\
		\begin{subfigure}[b]{0.45\textwidth} 
			
\includegraphics[width=.86\linewidth]{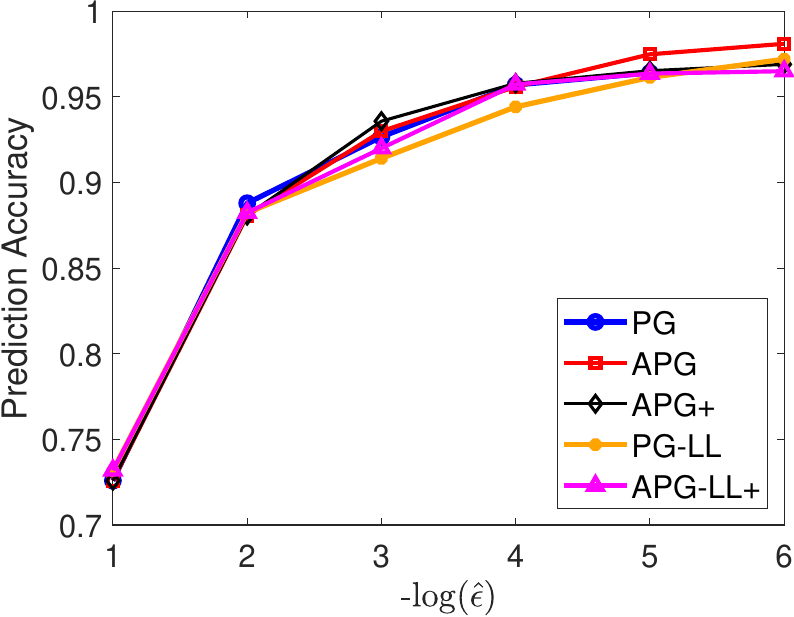}
			\caption{\webspam} 
		\end{subfigure}&
		\begin{subfigure}[b]{0.45\textwidth} 
\includegraphics[width=.86\linewidth]{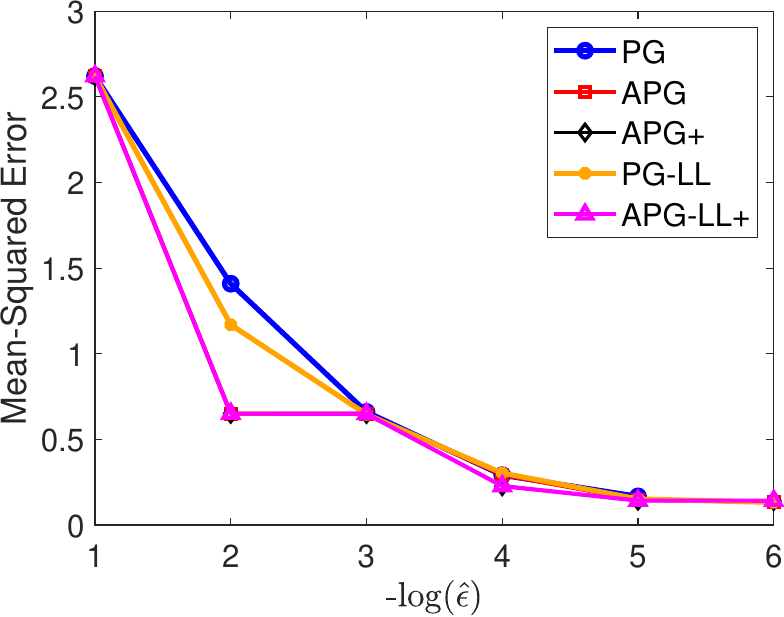}
			\caption{\El} 
			
		\end{subfigure}
	\end{tabular} 
	\caption{Prediction performance of different methods for sparse 
	logistic regression (\news, \rcv, \webspam) and least squares regression
	(\El) across varying levels of 
	residuals $\epsilon = 10^{-k}$, with $k=1,2,\dots, 6$. Generated plots 
	correspond to sparsity level of $s=\ceil{0.01m}$. }
	\label{fig:epsilon_sensitivity} 
\end{figure}

\alert{We present in \cref{fig:epsilon_sensitivity} the effect of varying the 
tolerance level $\hat{\epsilon}$ for the residual \cref{eq:opt}.
We can clearly see that in all cases, the prediction performance of
all methods keeps improving up to $\hat{\epsilon} = 10^{-5}$, which
indicates that our choice of a rather tight stopping condition is
indeed a suitable one for getting better prediction performance.
Note that in terms of comparison between different algorithms, 
the results in \cref{fig:epsilon_sensitivity} are consistent with that
in \cref{tbl:exp}.}

\subsection{\alert{Numerical comparison with a greedy method}}
\label{sec:comparison_greedy}
\begin{figure}[tbh] 
	\centering 
	\begin{tabular}{ccc}
		\begin{subfigure}[b]{0.32\textwidth} 
			\includegraphics[width=.86\linewidth]{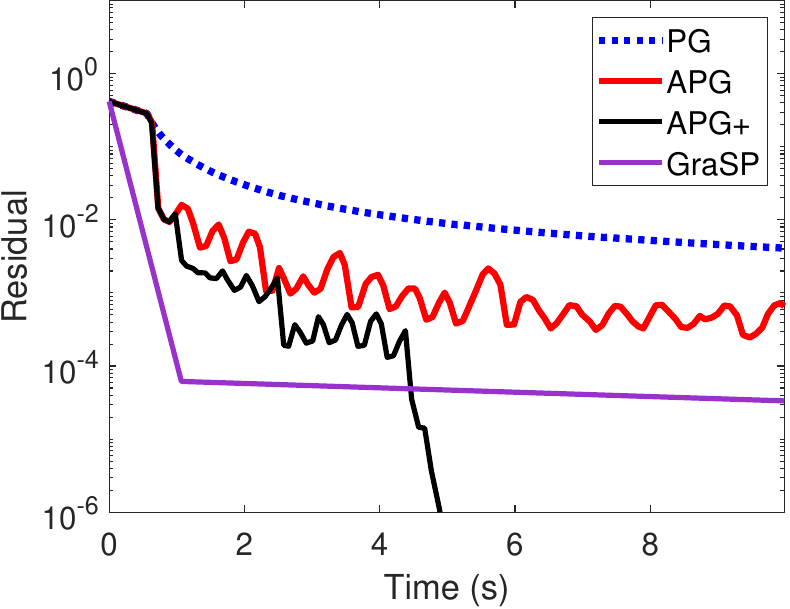}
			\caption{\news, $s=\ceil {0.01m}$} 
		\end{subfigure}& 
		\begin{subfigure}[b]{0.32\textwidth} 
			\includegraphics[width=.86\linewidth]{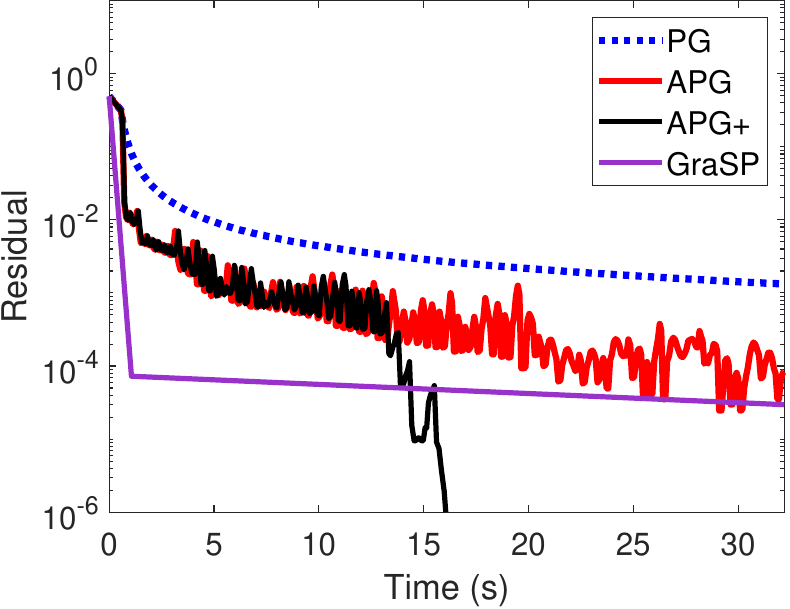}
			
			\caption{\news, $s=\ceil {0.05m}$} 
		\end{subfigure}& 
		\begin{subfigure}[b]{0.32\textwidth} 
			\includegraphics[width=.86\linewidth]{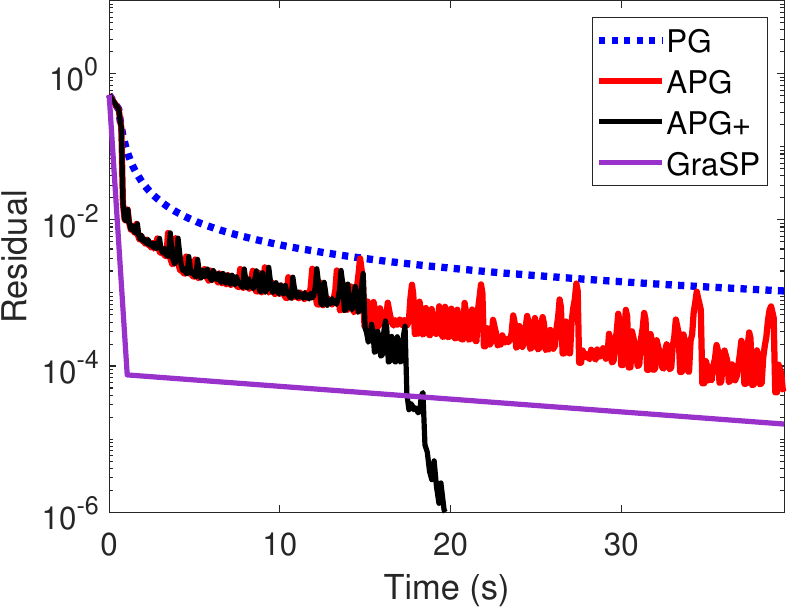}
			
			\caption{\news, $s=\ceil {0.1m}$} 
		\end{subfigure}
	\end{tabular} 
	\smallskip
		\begin{tabular}{ccc}
		\begin{subfigure}[b]{0.32\textwidth} 
			\includegraphics[width=.86\linewidth]{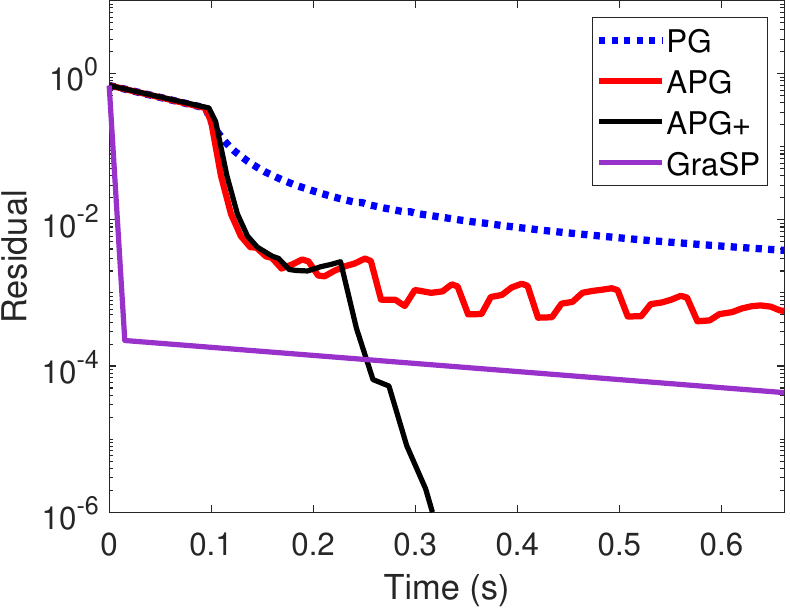}
			\caption{\rcv, $s=\ceil {0.01m}$} 
		\end{subfigure}& 
		\begin{subfigure}[b]{0.32\textwidth} 
			\includegraphics[width=.86\linewidth]{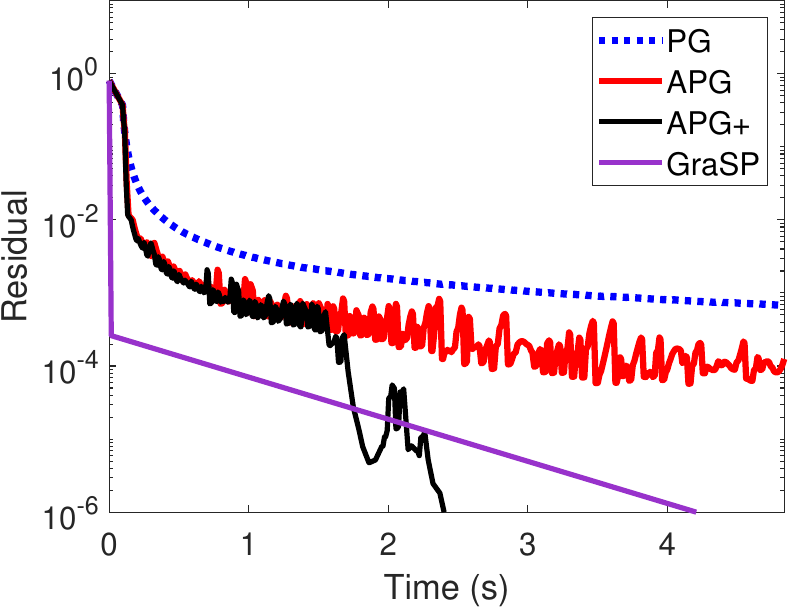}
			
			\caption{\rcv, $s=\ceil {0.05m}$} 
		\end{subfigure}& 
		\begin{subfigure}[b]{0.32\textwidth} 
			\includegraphics[width=.86\linewidth]{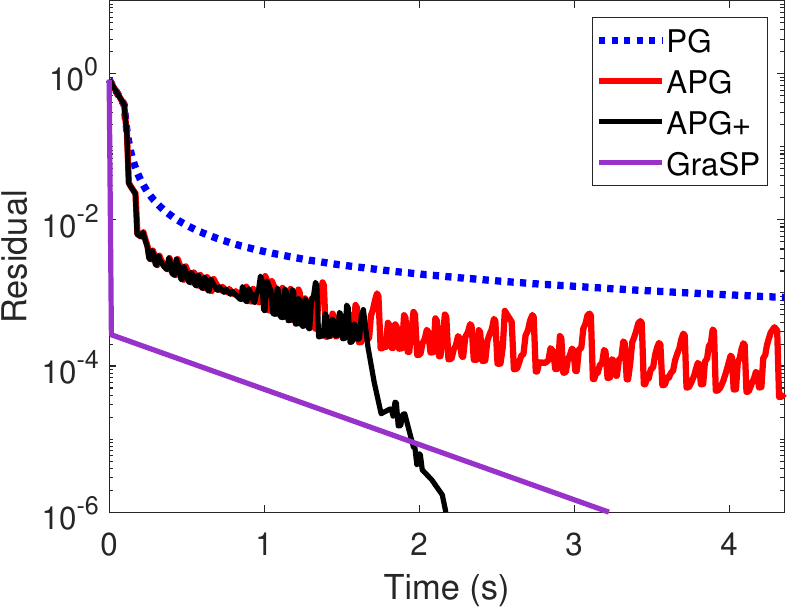}
			
			\caption{\rcv, $s=\ceil {0.1m}$} 
		\end{subfigure}
	\end{tabular} 
	\caption{Comparison of non-accelerated and accelerated projected gradient 
	methods with GraSP for solving regularized sparse logistic regression.}
	\label{fig:comparison_grasp} 
\end{figure}

\alert{
We present in \cref{fig:comparison_grasp} results of numerical
comparisons of our methods with the GraSP algorithm of
\citet{BahRB13a}, which is designed to solve our target problem \cref{eq:problem} 
for general loss functions $f$. Each iteration of GraSP involves a restricted 
minimization problem along a subspace of dimension at most $3s$, which is 
solved by a quasi-Newton approach.\footnote{We use their code for regularized
	sparse logistic regression downloaded
	from \url{https://sbahmani.ece.gatech.edu/GraSP.html}.} Note that
	GraSP is not 
ideal for large-scale datasets for its prohibitive memory consumption. For instance, it failed to fit the \webspam 
dataset even with the smallest sparsity level $s=\ceil{0.001m}$ \alertt{on our
machine with 64GB memory with an out of memory error}, whereas our 
proposed algorithm performs quite well for this instance. For a medium-sized 
dataset such as \news, we see from \cref{fig:comparison_grasp} that GraSP 
performs significantly slower than our proposed accelerated algorithm.
\alertt{In particular, its initial convergence is extremely fast, but it then
becomes stagnant after reaching a low-to-medium precision.}}

\clearpage 
\section{Proofs of Results in \cref{sec:pgm}}
\subsection{Proof of \cref{thm:pgm}}
\label{app:thmpgm}
\begin{proof}[Proof of part (a)]
	We note that the iterates of \cref{eq:pgm} are confined to $A_s$.
	Consider any accumulation point $w^*$ of $\{w^k\}$.
	For any given convergent subsequence $\{w^{k_r}\}$ with 
	$w^{k_r}\to w^*$, we obtain from the finiteness of $\J_s$ that
	there is $J \in \J_s$ such that
	\begin{equation}
	\label{eq:picked_subseq}
	w^{k_r+1} = 
	P_{A_J}(w^{k_r} - 
	\lambda \nabla f(w^{k_r}))
	\end{equation}
	 for infinitely many $r$. By taking subsequences if necessary, we assume 
	 that \cref{eq:picked_subseq} holds for all $r$ without loss of 
	 generality. Meanwhile, the PG 
	iterates \cref{eq:pgm} 
	can alternatively be written as
	\begin{equation}
	\label{eq:pgm_alternative}
	w^{k+1} \in \argmin _{y\in A_s} 
	~\inner{\nabla f(w^k)}{y-w^k} + 
	\frac{1}{2\lambda}\|y-w^k\|^2.
	\end{equation}
	We also use the 
	fact that \cref{eq:L-lipschitz} implies the well-known descent lemma 
	\cite[see for example,][Lemma 5.7]{Beck17}
	\begin{equation}
	\label{eq:descentlemma}
	f(w) \leq f(w') + \inner{\nabla f(w')}{w - w'} + 
	\frac{L}{2}\|w 
	- w'\|^2, \quad \forall w, w' \in \Re^n. 
	\end{equation}
	Noting that $\lambda \in (0,L^{-1})$, we have
	\begin{eqnarray*}
	&&f(w^{k+1}) \\
	& \overset{\cref{eq:descentlemma}}{\leq} & f(w^k) + 
	\inner{\nabla f(w^k)}{w^{k+1} - w^k} + 
	\frac{L}{2}\|w^{k+1} 
	- w^k\|^2 \\ 
	& = & f(w^k) + 
	\inner{\nabla f(w^k)}{w^{k+1} - w^k} + 
	\frac{1}{2\lambda}\|w^{k+1} 
	- w^k\|^2
	+ \left( \frac{\lambda L - 1}{2\lambda} \right) \|w^{k+1} 
	- w^k\|^2\\
	& \overset{\cref{eq:pgm_alternative}}{\leq} & f(w^k)  + \left( 
	\frac{\lambda L - 1}{2\lambda} \right) \|w^{k+1} 
	- w^k\|^2.
	\end{eqnarray*}
	That is, 
	\begin{equation}
	\frac{1-\lambda L}{2\lambda}\|w^k - w^{k+1}\|^2 \leq f(w^k) - 
	f(w^{k+1}).
	\label{eq:descent}
	\end{equation}
	\alert{If $w^{N}=w^{N+1}$ for some $N$, then it is clear from 
	\cref{eq:descent} 
	that $w^k=w^N$ for all $k\geq N$. Otherwise, }$\{f(w^k)\}$ is strictly 
	decreasing, proving the first 
	claim of part (a). Noting that $\nabla f$ and $P_{A_J}$ are continuous as 
	$A_J$ is closed and convex, we obtain from 
	\cref{eq:picked_subseq} that $w^{k_r+1} \to P_{A_J} (w^*-\lambda \nabla 
	f(w^*))$ as $r\to\infty$. On the other hand, we see from \cref{eq:descent} 
	and the lower boundedness of $f$
	that $\|w^k - w^{k+1}\|\to 0$, which, by means of the triangle inequality, 
	implies that $w^{k_r+1} \to w^*$. 
	Hence, we have $w^* = P_{A_J}   (w^*-\lambda \nabla 
	f(w^*))$.
	To complete the proof of part (a), we only need to show that 
	$P_{A_J} (w^*-\lambda \nabla 
	f(w^*)) \in P_{A_s} (w^*-\lambda \nabla 
	f(w^*))$. But from \cref{eq:picked_subseq} and \cref{eq:pgm}, we have that 
	for all $r$,  
	$w^{k_r} \in D_J$, where 
	\begin{equation}
	\label{eq:DJ}
	D_J\coloneqq \{ z \in \Re^n : \dist (z-\lambda \nabla f(z), 
	A_J) = \dist (z-\lambda \nabla f(z), A_s)\},
	\end{equation}
	where for any point $x$ and any set $A$, $\dist(x,A)$ is the
	distance from $x$ to $A$, defined as
	\[
		\dist(x,A) = \inf_{y \in A} \norm{x - y}.
	\]
	Since $A_J$ is closed, it 
	follows that $D_J$ is a closed set as well, and therefore $w^*\in D_J$. 
	That is, $P_{A_J} (w^*-\lambda \nabla 
	f(w^*)) \in P_{A_s} (w^*-\lambda \nabla 
	f(w^*))$, as desired. 

	Finally, we note that by representing \cref{eq:problem} as
	\[
		\min_{w}\quad f(w) + \delta_{A_s}(w),
	\]
	where $\delta_{A_s}$ is the indicator function of $A_s$ that
	outputs $0$ when $w \in A_s$ and infinity otherwise,
	a point $w^*$ is called stationary for \cref{eq:problem} if
	\[
		0 \in \partial f(w^*) + \partial \delta_{A_s}(w),
	\]
	where $\partial g$ is the limiting subdifferential in the sense of
	Clarke.
	On the other hand, the optimality condition of
	$w^* \in P_{A_s}(w^* - \lambda
	\nabla f(w^*))$ implies that
	\[
		0 \in w^* - \left( w^* - \lambda \nabla f(w^*) \right) +
		\partial \delta_{A_s}(w^*) \quad \Leftrightarrow \quad
		0 \in \lambda \nabla f(w^*) + \partial \delta_{A_s}(w^*).
	\]
	Since $\lambda \delta_{A_s} = \delta_{A_s}$ for any $\lambda$, the
	result above further implies that
	\[
		\lambda^{-1} 0 = 0 \in \nabla f(w^*) + \partial \delta_{A_s}(w^*),
	\]
	showing that $w^*$ is indeed a stationary point of 
	\cref{eq:problem}.
\end{proof}

\begin{proof}[Proof of part (b)]
	Suppose that $w^k \to w^*$ and define $I_{w^*}$ as in 
	\cref{eq:identification}. The finiteness of $\J_s$ 
	implies that there exists $\delta>0$ such that
	\begin{equation}
		B_{\delta}(w^*)\cap A_J = \emptyset,\quad \forall J\notin
		I_{w^*},
		\label{eq:split}
	\end{equation}
	where $B_{\delta}(w^*) \coloneqq \{ 
	w\in \Re^n : \|w-w^*\| < \delta\}$. Since $w^k\to w^*$, we can find $N>0$ 
	such that $w^k \in B_{\delta}(w^*)$ for all $k\geq N$. Hence, 
	\cref{eq:identification} immediately follows. 
	
	Now, suppose that $T_{\rm PG}(w^*)$ is a singleton for some accumulation 
	point $w^*$. Together with \cref{thm:pgm} (a), we have 
	\begin{equation}
	\label{eq:singleton}
	P_{A_s} (w^* - 
	\lambda \nabla f(w^*)) = \{w^*\}
	\end{equation}
	It is easy to verify that \cref{eq:singleton} implies that 
	\begin{equation}
	\label{eq:singleton2}
	w^* = 
	P_{A_J}(w^*-\lambda \nabla f(w^*)) \quad \forall J\in I_{w^*}.
	\end{equation}
	That is, $w^*$ 
	is a global minimum of $f$ over $A_J$ for all $J\in I_{w^*}$ due
	to the convexity of $f$.
	Now let $\delta>0$ be as defined in the preceding paragraph and $z\in A_s  
	\cap B_{\delta}(w^*)$. It then follows from \cref{eq:split} that
	$z\in A_J$ for some $J\in I_{w^*}$ .
	By the global minimality of $w^*$ for $f|_{A_J}$, it follows that
	$f(w^*)\leq f(z)$ for all $z\in A_s \cap B_{\delta}(w^*)$. That
	is, $w^*$ is a local minimum of $f$ over $A_s$.  
	

	It remains to show that the full sequence $\{w^k\}$ converges to $w^*$. To 
	this end, choose $\nu>0$ sufficiently small such that $z\in D_{J}\cap 
	B_{\nu}(w^*)$ for some $J$ implies $w^*\in D_J$, where $D_{J}$ is defined as 
	in 
	\cref{eq:DJ}. Note that such a $\nu$ exists as the collection $\{D_J : 
	J\in\J_s\}$ is finite and consists of closed sets. Let 
	$\{w^{k_r}\}_{r=0}^{\infty}$ be a subsequence converging to 
	$w^*$. We may assume without loss of generality that $w^{k_r} \in 
	B_{\nu}(w^*)$ for all $r\geq 0$. First, we show that $w^{k_0+1}\in 
	B_{\nu}(w^*)$. Note that the convexity of $f$ and 
	\cref{eq:L-lipschitz} result to 
	nonexpansiveness of the mapping $w \mapsto w-\lambda \nabla f(w)$
	because $\lambda \leq L^{-1}$
	\cite[Theorem 5.8]{Beck17}. From \cref{eq:pgm},
	$w^{k_0+1} = P_{A_{J}}\left(w^{k_0} -\lambda \nabla f(w^{k_0}) \right)$ for 
	some $J$. By the choice of $\nu$ and the fact that $\Tpg (w^*)$ consists of 
	one element, we see that $w^* = P_{A_J}(w^*-\lambda \nabla f(w^*))$. With these, 
	we have
	\begin{align}
		\nonumber
	\|w^{k_0+1} - w^*\| & = \| P_{A_{J}}\left(w^{k_0} -\lambda \nabla 
	f(w^{k_0}) \right) - P_{A_{J}} (w^* - 
	\lambda \nabla f(w^*)) \| \\
	\label{eq:firstbound}
	& \leq  \| (w^{k_0} - \lambda \nabla f(w^{k_0}) - (w^* - \lambda \nabla 
	f(w^*))\| \\
		\nonumber
	& \leq \|w^{k_0} - w^*\| \\
		\nonumber
	& < \nu ,
	\end{align}
	where \cref{eq:firstbound} follows from the nonexpansiveness of projection 
	mappings onto closed convex sets, while the second inequality follows from 
	the nonexpansiveness of $w \mapsto w-\lambda \nabla f(w)$. Proceeding 
	inductively, we see that $w^k 
	\in B_{\nu}(w^*)$ for all $k\geq k_0$ and 
	$\{\|w^k-w^*\|\}_{k=k_0}^{\infty}$ 
	is a decreasing sequence. As its subsequence $\{\|w^{k_r}-w^*\|\}$ 
	converges to zero, it follows that $w^k\to w^*$, as desired. 
\end{proof}

\begin{proof}[Proof of part (c)]
	Note that $w \mapsto w-\lambda \nabla f(w)$ being a 
	contraction implies that there exists $\gamma \in [0,1)$ such that 
	\begin{equation}
	\| (w^{k_0} - \lambda \nabla f(w^{k_0}) - (w^* - \lambda \nabla 
	f(w^*))\| \leq \gamma \|w^{k_0} - w^*\|.
	\label{eq:contraction}
	\end{equation}
	We then obtain the desired inequality \cref{eq:Q-linear} by
	combining \cref{eq:contraction} and \cref{eq:firstbound}.
\end{proof}

\subsection{Proof of \cref{thm:sublinear}}
\begin{proof} 
We first prove \cref{eq:sublinear}.
We have from \cref{thm:pgm} (b) that there exists $N>0$ such that 
\cref{eq:identification} holds. Moreover, recall from 
the proof of \cref{thm:pgm} that 
	\begin{equation}
	\label{eq:xi=argmin}
	w^{k} \in \argmin _{y\in\Re^n}\quad \inner{ \nabla f(w^{k-1})}{y-w^{k-1}} 
	+ 
	\frac{1}{2\lambda}\|y-w^{k-1}\|^2 + \delta_{A_s} (y),
	\end{equation}
where $\delta_{A_s}$ is the indicator function of $A_s$. If
$J_k\in\J_s$ satisfies $w^k\in A_{J_k}$, we can alternatively write
\cref{eq:xi=argmin} as
	\begin{equation}
		\label{eq:xi=argmin_AJ}
	w^{k} \in \argmin _{y\in\Re^n}\quad Q_{J_k}^{k-1}(y) \coloneqq \inner{ \nabla 
	f(w^{k-1})}{y-w^{k-1}} 
	+ 
	\frac{1}{2\lambda}\|y-w^{k-1}\|^2 + \delta_{A_{J_k}} (y).
	\end{equation}
Recognizing the right-hand side of \cref{eq:xi=argmin_AJ} as a strongly convex 
function of $y$, we obtain
	\begin{equation}
	\label{eq:strongconvexity}
	Q_{J_k}^{k-1}(y) - Q_{J_k}^{k-1}(w^k) \geq 
	\frac{1}{2\lambda}\|y-w^k\|^2, \quad 
	\forall y\in A_{J_k}.
	\end{equation}
Consequently, for any $k\geq N$, we have from \cref{eq:identification} that 
${J_k}\in I_{w^*}$ and so $w^* \in A_{J_k}$, which together with
\cref{eq:strongconvexity} implies
	\begin{equation}
	\label{eq:strongconvexity_star}
	Q_{J_k}^{k-1}(w^*) - Q_{J_k}^{k-1}(w^k) \geq 
	\frac{1}{2\lambda}\|w^*-w^k\|^2 \quad \forall 
	k\geq N. 
	\end{equation}
Meanwhile, by the descent lemma \cref{eq:descentlemma} and by noting that 
$\lambda\in(0,L^{-1})$, we have $f(w^{k-1}) + 
Q_{J_k}^{k-1}(w^k) \geq f(w^k)$. Thus, for all $k\geq N$,
	\begin{eqnarray*}
	\frac{1}{2\lambda}\|w^*-w^k\|^2 &  \leq &
 	Q_{J_k}^{k-1}(w^*) + f(w^{k-1}) - f(w^k) \\
 	& \overset{\cref{eq:xi=argmin_AJ}}{=} & \inner{ \nabla 
 		f(w^{k-1})}{w^*-w^{k-1}} 
 	+ 
 	\frac{1}{2\lambda}\|w^*-w^{k-1}\|^2 + f(w^{k-1}) - f(w^k)  \\
 	& \leq & f(w^*) - f(w^k) + \frac{1}{2\lambda}\|w^*-w^{k-1}\|^2 
 	\end{eqnarray*}
 where the last inequality holds by the convexity of $f$. Thus, we have
 	\begin{equation}
 	\sum_{j=N+1}^k \left( f(w^j) - f(w^*) \right) \leq \frac{1}{2\lambda} 
 	\left( \|w^*-w^{N}\|^2  - \|w^*-w^{k}\|^2 \right).
 	\end{equation}
Hence, \cref{eq:sublinear} immediately follows by noting that $\{ f(w^k)\}$ 
is monotonically decreasing, as proved in \cref{thm:pgm} (a), and
applying \cite[Lemma~1]{LeeW19a}.

Now we turn to \cref{eq:linear}.
From the convexity of $f$, we have that
\begin{equation}
	f(w^k) - f(w^*) \leq \inner{\nabla f(w^k)}{w^k - w^*}.
	\label{eq:convexity}
\end{equation}
By \cref{eq:singleton2}, we can
easily conclude that 
\begin{equation}
\label{eq:nabla_on_J=0}
(\nabla f(w^*))_J = 0 \quad \forall J\in I_{w^*}.
\end{equation}
Meanwhile, through \cref{eq:identification}, there exists $N$ such that for all 
$k
\geq N$, we can find $J_k \in \I_{w^*}$ such that
$w^k \in A_{J_k}$. Thus, we see that
	\begin{equation}
	\label{eq:inner}
	\inner{\nabla f(w^*)}{w^k - w^*} = 0, \quad \forall k \geq N,
	\end{equation}
because only entries of $\nabla f(w^*)$ outside $J_k\in I_{w^*}$ could be 
nonzero,
but those entries are identically $0$ for both $w^k$ and $w^*$, that is, 
$w^k_i 
= w_i^* = 0$ for all $i\notin J_k$.
We therefore proceed on with \cref{eq:convexity} as follows:
\begin{align}
	\nonumber
f(w^k) - f(w^*) &\overset{\cref{eq:inner}}{\leq} \inner{\nabla f(w^k) - \nabla 
f(w^*)}{w^k -
w^*}\\
	\nonumber
&\leq \norm{\nabla f(w^k) - \nabla f(w^*)}\norm{w^k - w^*}\\
&\overset{\cref{eq:L-lipschitz}}{\leq}  L \norm{w^k-w^*}^2, \quad \forall k 
\geq N,
\label{eq:quadraticbound}
\end{align}
where the second inequality is from the Cauchy-Schwarz inequality.
Finally, \cref{eq:linear} is proven by inserting \cref{eq:Q-linear}
into \cref{eq:quadraticbound}.
\end{proof}

\section{Proof of Results in \cref{sec:accelerate}}
\subsection{Proof of \cref{thm:global_acce}}
\begin{proof}
	Note that for any $k\geq 0$, we have from \cref{eq:suffdecrease} that 
		\begin{equation}
		\label{eq:suffdecrease_gen}
		f(z^k) \leq f(w^k) - \frac{\sigma}{2}t_k^2 
		\|d^k\|^2, \quad z^k \coloneqq w^k+t_kd^k,
		\end{equation}
	where $t_k$ is defined to be zero if the condition in
	\cref{line:if} of \cref{alg:extrapolation} is not 
	satisfied. Analogous to \cref{eq:descent}, we have from 
	$w^{k+1} \in \Tpg (z^k)$ that
		\begin{equation}
		\label{eq:descent2}
		\frac{1-\lambda L}{2\lambda}\|z^k - w^{k+1}\|^2\leq f(z^k) -
		f(w^{k+1}).
		\end{equation}
	Using \cref{eq:descent2} together with \cref{eq:suffdecrease_gen}, we have 
		\begin{equation}
		\frac{\sigma}{2}t_k^2 \|d^k\|^2 \leq f(w^k) - f(w^{k+1}).
		\end{equation}
	Then $\{f(w^k)\}$ is decreasing, and since $f$ is bounded below over $A_s$, 
	we have 
		\begin{equation}
		\label{eq:t_theta_d}
		t_k^2\|d^k\|^2 \to 0.
		\end{equation}
	
	Now, assume that $\{w^{k_r}\}$ is a subsequence of a 
	sequence generated by \cref{alg:extrapolation} that converges to $w^*$, and 
	as in the proof of \cref{thm:pgm}, we assume
	that there exists $J\in\J_s$ such that
		\begin{equation}
		w^{k_r+1} = P_{A_J}(z^{k_r} - \lambda \nabla f(z^{k_r}))
		\end{equation}
	 for all $r$. Then from 
	 \cref{eq:t_theta_d}, we have that $z^{k_r} \to w^*$ so 
	 that $w^{k_r+1} \to P_{A_J}(w^* - \lambda \nabla f(w^*))$. Meanwhile, 
	 \cref{eq:descent2} gives $\|z^k-w^{k+1}\|\to 0$, and therefore 
	 $w^{k_r+1}\to w^*$. Hence, $w^* = P_{A_J}(w^* - \lambda \nabla f(w^*))$. 
	The rest now follows from exactly the same arguments used in the latter 
	part of the proof of \cref{thm:pgm} (a).
\end{proof}

\subsection{Proof of \cref{thm:rate2}}


\begin{proof}
	We consider the sequence $\{w^0, z^0, w^1, z^1, w^2,z^2, \dotsc,
	w^k,z^k , \dotsc \}$ and remove those $z^k$ with $z^k = w^k$ from
	the sequence, and call the resulting sequence 
	$\{\bar{w}^m\}$. We will show that actually $\bar{w}^m \rightarrow
	w^*$, and since $\{w^k\}$ is a subsequence of $\{\bar{w}^m\}$, its
	convergence to the same point will ensue.
	Note that if $n_k$ denotes the number of successful 
	extrapolation steps in the first $k$ iterations of
	\cref{alg:extrapolation}, then $w^k=\bar{w}^{k+n_k}$. Moreover, 
	it is easy to check that $w^*$ is likewise an accumulation point of 
	$\{\bar{w}^m\}$. 
	
	To prove 
	the desired result, we first show that the following properties are 
	satisfied:
	\begin{description}
		\item[(H1)] There exists $a>0$ such that 
			\begin{equation}
			\label{eq:H1}
			f(\bar{w}^{m}) \leq f(\bar{w}^{m-1}) - 
			a\norm{\bar{w}^{m}-\bar{w}^{m-1}}^2 \quad \forall m\in 
			\mathbb{N}.
			\end{equation}
		\item[(H2)] There exists $b>0$ such that for all $m\in \mathbb{N}$, 
		there is a vector $v^{m} \in \partial \delta_{A_s} (\bar{w}^{m})$ 
		satisfying 
			\begin{equation}
			\label{eq:H2}
			\norm{\nabla f(\bar{w}^m) + v^m} \leq b 
			\norm{\bar{w}^{m}-\bar{w}^{m-1}}.
			\end{equation}	
	\end{description}
	To this end, fix $m\in \mathbb{N}$, and we separately consider two cases:
	$\bar{w}^m = w^k$ for some $k$ and $\bar{w}^m = z^k$ for some $k$.

	\textbf{Case I: $\bar{w}^m = w^k$.}\\
	First, suppose 
	that $\bar{w}^m = w^k$ for some $k$. Then $\bar{w}^{m-1} \in \{ w^{k-1}, 
	z^{k-1}\}$. In either case, we have from \cref{eq:descent} or
	\cref{eq:descent2} that 
		\begin{equation}
		\label{eq:descent3}
		f(\bar{w}^{m}) \leq f(\bar{w}^{m-1}) - 
		\frac{1-\lambda L}{2\lambda}\norm{\bar{w}^{m}-\bar{w}^{m-1}}^2, 
		\end{equation}
	that is, \cref{eq:H1} is satisfied with $a = (1-\lambda L)/(2\lambda)$. 
	On the other hand, since $\bar{w}^m \in \Tpg (\bar{w}^{m-1})$, similar
	to  \cref{eq:xi=argmin_AJ}, we have
		\begin{equation}
		\label{eq:xi=argmin_As}
				\bar{w}^{m} = \argmin_{y\in \Re^n}\quad \inner{\nabla 
				f(\bar{w}^{m-1})}{y-\bar{w}^{m-1}} + 
		\frac{1}{2\lambda}\|y-\bar{w}^{m-1}\|^2 + \delta_{A_{J_k}}(y),
		\end{equation}
	where $J_k\in \J_s$ satisfies $\bar{w}^m = w^k \in A_{J_k}$. From the 
	optimality condition of \cref{eq:xi=argmin_As}, we have
		\begin{equation}
		\label{eq:optimality condition}
		0\in \nabla f(\bar{w}^{m-1} )+ \frac{1}{\lambda}\left( \bar{w}^m - 
		\bar{w}^{m-1}\right)  + \partial \delta_{A_{J_k}}(\bar{w}^m).
		\end{equation}
	That is, 
		\begin{equation*}
		v^m \coloneqq  - \nabla f(\bar{w}^{m-1} ) - \frac{1}{\lambda}\left( 
		\bar{w}^m - 
		\bar{w}^{m-1}\right) \in  \partial \delta_{A_{J_k}} (\bar{w}^m).
		\end{equation*}
	From \cite[Equation~(19)]{HLN14}, the above equation implies $v^m\in \partial
	\delta_{A_s}(\bar{w}^m)$. Moreover, 
		\begin{eqnarray}
		\norm{\nabla f(\bar{w}^m) + v^m}& = &  \norm{\nabla f(\bar{w}^m)- 
\nabla f(\bar{w}^{m-1} ) - \frac{1}{\lambda}\left( 
		\bar{w}^m - \bar{w}^{m-1})\right)} \nonumber \\
		& \overset{\cref{eq:L-lipschitz}}{\leq} & \left( L + 
		\frac{1}{\lambda}\right)\norm{\bar{w}^m - \bar{w}^{m-1}},
		\end{eqnarray}
		that is, \cref{eq:H2} is satisfied with $b = L+\lambda^{-1}$. 
	
	\textbf{Case II: $\bar{w}^m = z^k$.}\\
	We now consider the other possibility that $\bar{w}^m = z^k$, in which 
	case, we 
	necessarily have $\bar{w}^{m-1} = w^k$, $d^k \neq 0$, and $t_k>0$. By 
	\cref{eq:suffdecrease}, \cref{eq:H1} is already satisfied with $a = \sigma$. 
	To prove that (H2) holds, we first bound the final step size $t_k$
	in the line search procedure.
	Using \cref{eq:descentlemma}, we see that
	\cref{eq:suffdecrease} is satisfied if
	\begin{equation*}
	t_k \inner{\nabla f(w^k)}{
		d^k} + \frac{{t}_k^2 L}{2} \norm{d^k}^2 \leq -\sigma
	t_k^2 \norm{d^k}^2 
	\end{equation*}
	or equivalently, recalling that $t_k>0$, we have
	\begin{equation}
	-\inner{\nabla f(w^k)}{d^k} \geq \left(\frac{L}{2} +
	\sigma)\right) t \norm{d^k}^2.
	\label{eq:bound2}
	\end{equation}
	Therefore, \cref{eq:suffdecrease} is satisfied whenever
	\begin{eqnarray}
		\nonumber
	t_k &\stackrel{\cref{eq:bound2}}{\leq} &
	- \frac{\inner{\nabla 
	f(w^k)}{ d^k}}{\norm{d^k}^2 \left( \frac{L}{2} + \sigma
\right)}\\
		\nonumber
	&\stackrel{d^k \in A_J}{=}&
	- \frac{\inner{(\nabla  f(w^k))_{J}}{(
			d^k)_J}}{\norm{d^k}^2 \left( \frac{L}{2} + \sigma \right)}\\
		\label{eq:lastbound}
	&\stackrel{\cref{eq:cosine}}{=}& \frac{\zeta_k \norm{(\nabla  
			f(w^k))_{J}}}{\left( \frac{L}{2} + \sigma
		\right)\norm{d^k}}.
	\end{eqnarray}
	By applying the condition of $\zeta_k \geq \epsilon$ to
	\eqref{eq:lastbound},
	we get that \cref{eq:suffdecrease} is satisfied whenever
	\[
		t_k \leq
	\frac{\epsilon \norm{(\nabla  f(w^k))_{J}}}{\left( \frac{L}{2} +
	\sigma\right) \norm{d^k}}.
	\]
	Therefore, we see that $t_k$ is lower-bounded by
	\begin{align*}
		t_k &\geq \min\left\{c_k \alpha_{\min}, 
		\frac{\eta \epsilon \norm{(\nabla  f(w^k))_{J}}}{\left( \frac{L}{2} +
		\sigma\right) \norm{d^k}} \right\} \\
			&\stackrel{\cref{eq:cosine}}{=}
			\min\left\{\frac{\norm{(\nabla f(w^k))_{J} }}{\zeta_k
				\norm{d^k}} \alpha_{\min}, 
			\frac{\eta \epsilon \norm{(\nabla  f(w^k))_{J}}}{\left( \frac{L}{2} +
			\sigma\right) \norm{d^k}} \right\}\\
			&\geq \frac{\norm{(\nabla f(w^k))_{J} }}{\norm{d^k}}
			\min\left\{
				\alpha_{\min}, 
			\frac{\eta \epsilon }{\left( \frac{L}{2} +
			\sigma\right)} \right\},
	\end{align*}
	where the factor of $\eta$ is to consider the possibility of
	overshooting and the last inequality is from that $\zeta_k \in
	(0,1]$ in \cref{eq:cosine}.
	We thus conclude that for the final update $t_k d^k$, we have
	\begin{equation}
	\norm{t_k d^k} \geq \norm{(\nabla  f(w^k))_{J}}
	\underline{t}, \quad \underline{t} \coloneqq
			\min\left\{
				\alpha_{\min}, 
			\frac{\eta \epsilon }{\left( \frac{L}{2} +
			\sigma\right)} \right\}.
	\label{eq:subgradlike}
	\end{equation}
	We then get from \cref{eq:subgradlike,eq:L-lipschitz} that
	\begin{equation}
	\norm{(\nabla f(z^k))_{J} } \leq \norm{(\nabla f(z^k))_{J}  -
		(\nabla f(w^k))_{J} } + \norm{(\nabla f(w^k))_{J} }
	\leq (L + \underline{t}^{-1}) \norm{z^k -
		w^k}.
	\label{eq:suff}
	\end{equation}
	We now furnish the vector required in (H2). Let $v^m\in \partial 
	\delta_{A_s}(\bar{w}^m)$ such that $(v^m)_{J} = 0$ and $(v^m)_{J^c} 
	\coloneqq -(\nabla f(z^k))_{J^c}$. Then by \cref{eq:suff}, it is clear that
		\[\norm{ \nabla f(\bar{w}^m) + v^m} = \norm{(\nabla f(z^k))_{J} } \leq 
		(L + \underline{t}^{-1}) \norm{\bar{w}^m - \bar{w}^{m-1}}. \]
	Setting $a = \min\left\{(1-\lambda L)/(2\lambda) , \sigma \right\}$ and 
	$b = \max \left\{ L + \lambda^{-1}, L+\underline t ^{-1}
	\right\}$, we see that (H1) and (H2) are both satisfied for case I and case 
	II. 
	
	The rest of the proof for convergence to $w^*$ will follow from
	arguments analogous to those used in \cite{HA13a}, with the only
	deviation that our condition \cref{eq:KL2} is weaker than the KL
	condition assumed in \cite{HA13a}.
	Through a careful inspection of the proof of
	\cite[Lemma~2.6,Corollary~2.8]{HA13a}, we see that
		\begin{equation}
			\label{eq:upperbound}
			(f(\bar{w}^m)-f(w^*))^{\alerttt{ \theta}} \leq
		\alerttt{b\kappa}\norm{\bar{w}^m-\bar{w}^{m-1}} , \quad
		\forall m\in \mathbb{N} 
		\end{equation}
	is a key inequality for the proof.
	The above inequality clearly holds when the conventional KL
	condition holds, and here we will show how \cref{eq:upperbound} will
	still hold under \cref{eq:KL2} so
	all remaining arguments in the proof of
	\cite[Lemma~2.6,Corollary~2.8]{HA13a} ensue to be valid.
	In either case I or case II above, note that the 
	vector $v^m$ in (H2) has the property that $(v^m)_{J}=0$, where $J\in \J_s$ 
	is the index set satisfying $\bar{w}^m\in A_J$. It follows that $(\nabla f(\bar{w}^m))_J$ is 
	a subvector of $\nabla f(\bar{w}^m)+v^m$, so
	\begin{equation}
		\norm{(\nabla f(\bar{w}^m))_J} \leq \norm{\nabla
			f(\bar{w}^m)+v^m}.
			\label{eq:subgrad}
	\end{equation}
	Using \cref{eq:H2,eq:subgrad,eq:KL2},
	we immediately obtain \cref{eq:upperbound}, as desired. 
	The convergence of the
	iterates then follows from \cite[Theorem 2.9]{HA13a}.
	
	As for the rates, we have from \cref{eq:upperbound} and \cref{eq:H1} that 
	\[a(f(\bar{w}^m)-f(w^*))^{\alerttt{2\theta}} \leq
		\alerttt{b^2\kappa^2} \left( 
		f(\bar{w}^{m-1}) - f (\bar{w}^m) \right). \]
	That is, 
	\begin{equation}
		aD_m^{\alerttt{2\theta}} \leq \alerttt{b^2\kappa^2}(D_{m-1}-D_m), \quad D_m \coloneqq 
		f(\bar{w}^m)-f(w^*).
		\label{eq:rate}
	\end{equation}
		We now separately consider different values of $\theta$.
		One result that is being used repeatedly in our discussion
		below is that we have from $\bar{w}^m \rightarrow w^*$ and the
		continuity of $f$ that
		\begin{equation}
			D_m \downarrow 0.
			\label{eq:Dm}
		\end{equation}
		Our proof for $\theta \in \alerttt{(1/2,1)}$ is inspired by the proof of
		Lemma~6 in \cite[Chapter~2.2]{Pol87a}.
		\begin{enumerate}[leftmargin=*]
			\item[(a)] When $\theta \in \alerttt{(1/2,1)}$,
				\cref{eq:rate} implies
				\begin{equation}
					D_m\left(\frac{a }{\alerttt{\kappa^2} b^2}
					D_m^{\alerttt{2\theta-1}}
					+ 1 \right) \leq D_{m-1},
					\label{eq:ratesublinear}
				\end{equation}
				and since $\alerttt{2\theta-1} \in (0,1)$,
				\cref{eq:ratesublinear} leads to
				\begin{equation}
					\label{eq:base}
					D_m^{-(\alerttt{2\theta-1})}
					\left( 1 + \frac{a }{\alerttt{\kappa^2}b^2}
					D_m^{\alerttt{
					2\theta-1}}\right)^{-\alerttt{(2\theta-1)}}
					\geq D_{m-1}^{-(
					\alerttt{2\theta-1})},
				\end{equation}
				and we have from \cref{eq:Dm} that
				$D_m^{\alerttt{2\theta-1}}
				\downarrow 0$.
				Therefore, we can find $k_0 \geq 0$ such that
				\[
					\frac{a }{\alerttt{\kappa^2}b^2}
					D_m^{\alerttt{2\theta-1}} < 1,
					\quad \forall m \geq k_0.
				\]
				As $-(\alerttt{2\theta-1}) \in (-1,0)$, for $m
				\geq k_0$ we get
				\begin{equation}
					\left( 1 + \frac{a }{\alerttt{\kappa^2}b^2}
					D_m^{\alerttt{
					2\theta-1}}\right)^{-(\alerttt{2\theta-1})} \leq 1
					+ (2^{\alerttt{-
					2\theta+1}} - 1) \frac{a }{\alerttt{\kappa^2}b^2}D_m^{
					\alerttt{2\theta-1}}.
					\label{eq:LHS}
				\end{equation}
				By combining \cref{eq:LHS,eq:base}, we get that for $m
				\geq k_0$,
				\begin{equation}
					D_m^{-(\alerttt{2\theta-1})} - (1 - 2^{\alerttt{1-2\theta}})
					\frac{a }{\alerttt{\kappa^2}b^2} \geq D_{m-1}^{-\alerttt{(
					2\theta-1})}.
					\label{eq:torate}
				\end{equation}
				We note that for $\theta \in \alerttt{(1/2,1)}$,
				$2^{\alerttt{-2\theta+1}}
				\in (1/2,1)$, so 
				\[
					C_\theta \coloneqq (1 - 2^{\alerttt{-2\theta+1}})
					\frac{a }{ \alerttt{\kappa^2}b^2} > 0.
				\]
				Thus, by summing \cref{eq:torate} for $m=k_0+1,
				k_0+1,\dotsc, k+n_k$ and telescoping,
				we get
				\[
					D_{k+n_k} \leq \left((k+n_k - k_0)
					C_\theta +
					D_{k_0}^{-(\alerttt{2\theta-1})}\right)^{\frac{-1}{\alerttt{2\theta-1}}}=
					O\left((k+n_k)^{\frac{-1}{\alerttt{2\theta-1}}}\right),
				\]
				as desired.

			\item[(b)] When $\theta = 1/2$, we see that \cref{eq:rate}
				reduces to
				\begin{equation}
					\left(a + b^2\alerttt{\kappa^2}\right) D_m \leq {b^2}\alerttt{\kappa^2} D_{m-1}
					\quad \Leftrightarrow \quad
					D_m \leq \frac{{b^2}\alerttt{\kappa^2}}{\left(a +
						{b^2}\alerttt{\kappa^2}\right)} D_{m-1},
				\label{eq:Qlinear}
				\end{equation}
				which shows a $Q$-linear convergence rate (as $a > 0$)
				that directly implies the desired exponential bound.
				
				For $\theta \in \alerttt{(0,1/2)}$, we get
				$\alerttt{2\theta} < 1$, and
				Thus, by the monotonicity of $\{f(\bar{w}^m)\}$, we
				can find $k_0 > 0$ such that $D_m \leq 1$ for all $m
				\geq k_0$.
				For such $m$, \cref{eq:rate} gets us
				\[
					a D_m \leq a D_m^{\alerttt{2\theta}} \leq
					{b^2}\alerttt{\kappa^2} \left(D_{m-1} - D_m\right),
					\quad \forall m \geq k_0,
				\]
				and the same $Q$-linear rate and exponential bound
				then follow from \cref{eq:Qlinear} and the argument
				that followed it.
			\item[(c)] When $\theta = \alerttt{0}$, \cref{eq:rate} becomes
				\[
					\frac{a}{\alerttt{\kappa^2}b^2} \leq \left( D_{m-1} - D_m
					\right).
				\]
				Hence, noting that $D_{m-1} -  D_m \to 0$ by \cref{eq:Dm} and 
				$a/(\alerttt{\kappa^2}b^2)>0$,  
				there must be $k_0\geq 0$ such that
				$D_m = 0$ for all $m \geq k_0$.
		\end{enumerate}
\end{proof}

\subsection{Proof of \cref{thm:newton}}
\begin{proof}
We will first establish the quadratic convergence of
$\{w^{k,j}\}_j$ to $w^*$ when $t$ approaches infinity.
The overall quadratic convergence can then be obtained by showing that
the iterates will all stay within the same $A_J$ and applying
\cref{thm:superlinear} in \cref{sec:superlinear}.

For the part of $\{w^{k,j}\}_j$ for a given $k$,
we note that since $\nabla^2 f_{J}(w^*)$ is positive definite and
$w^{k,0} \in U \cap A_J$,
$w^*_J$ is an isolated global optimum of $f_J$ (as $f_J$ is convex).
Moreover, the algorithm in \cref{eq:newton} clearly treats coordinates not in
$J$ as nonvariables, and thus the whole sequence of $\{w^{k,j}\}_j$
stays in $A_J$.
Therefore, $\{w^{k,j}\}$ converges quadratically to $w^*$ following
standard analysis for Newton methods; see, for example,
\citep[Chapter~3]{NoceWrig06}.
To satisfy the conditions of \cref{thm:superlinear}, we just need to
notice that if we group $t \geq 1$ consecutive Newton iterations as
the operation $T_2$,
the convergence speed is $2t \geq 2$, so the quadratic convergence
assumption is still satisfied.
It is also clear that since $w^*$ is stationary for
\cref{eq:problem},
$\nabla f_J(w^*_J) = 0$ and thus $w^*$ is a fixed point for the
Newton steps.
For \cref{eq:pgm}, clearly these suffice for our usage of
\cref{thm:superlinear} to reach the conclusion.
%
\end{proof}

\section{Superlinear convergence of \cref{alg:identification}}
\label{sec:superlinear}
In this section, we state and prove a general result of a two-step superlinear
convergence of \cref{alg:identification} that is similar in spirit
to that in \cite{BarIM20a} to simply assume that we have a
superlinearly convergent subroutine.
We consider this abstract form to demonstrate the versatility of our
framework and to allow full flexibility to
accommodate different problem conditions of $f|_{A_J}$, and also to fit
various algorithms like inexact damped/regularized (semismooth) Newton or
quasi-Newton methods, instead of giving the impression that we are restricted
to a certain algorithm.

\begin{theorem}
\label{thm:superlinear}
Assume that we have a mapping $T_1(w)$ such that its generated
iterates $\{w^k\}$ with $w^{k+1} \in T_1(w^k)$ converge to a
stationary point $w^*$ of \cref{eq:problem} and
\begin{equation}
	\norm{\hat w - T_1(w^*)} \leq \norm{w - w^*}, \quad
	\forall \hat w \in T_1 \left( w \right)
	\label{eq:nonexp}
\end{equation}
for all $w$ in a neighborhood $U$ of $w^*$ and in some $A_J$ with $J$
satisfying $J \in \I_{w^*}$,
and that there is another mapping $T_2$ that, when given an initial
point $w^0 \in A_J$, generates iterates that are all in $A_J$ and
superlinearly convergent to $w^*$ within $U$ for each $J \in
\I_{w^*}$ with $T_2(w^*) = w^*$, then the iterates generated by
\begin{equation}
	w^{k+1} \in T_2\left(T_1\left(w^{k}\right)\right)
	\label{eq:seq}
\end{equation}
converge to $w^*$ at the same superlinear rate as that of $T_2$.
\end{theorem}
\begin{proof}
We assume without loss of generality that
\begin{equation}
	\norm{T_2(w) - w^*} \leq c \norm{w - w^*}^{1+\rho}
	\label{eq:superlinear}
\end{equation}
for some $c, \rho > 0$ for all $w \in A_J \cap U$ for all $J \in
\I_{w^*}$.
Then by \cref{eq:nonexp}, and by denoting
\[
	\hat w^{k+1} \in T_1(w^k),
\]
as the element in $T_1(w^k)$ leading to $w^{k+1}$,
we obtain
\begin{align*}
	\norm{w^{k+1} - w^*} &= \norm{T_2 \left(\hat w^k  \right) -
	w^*}\\
	&\leq  c \norm{\hat w^k - w^*}^{1+\rho}\\
	&= c \norm{T_1(w^k) - w^*}^{1+\rho} \\
	&= c \norm{w^k - w^*}^{1+\rho},
\end{align*}
where the the first inequality is from \cref{eq:superlinear}.
Therefore, the conclusion of the theorem is proven.	
\end{proof}

%

\end{document}